\documentstyle[amsfonts,amssymb]{article}

\newcommand{\ndt}{\noindent}

\newcommand{\beq}{\begin{eqnarray}}
\newcommand{\eeq}{\end{eqnarray}}
\newcommand{\beqst}{\begin{eqnarray*}}
\newcommand{\eeqst}{\end{eqnarray*}}

\newcommand{\ol}{\overline}

\newcommand{\sgn}{\mbox{\rm sign\,}}

\newcommand{\R}{{\mathbb R}}
\newcommand{\C}{{\mathbb C}}

\newcommand{\qed}{\hfill $\square$}

\newcommand{\dsp}{\displaystyle}

\newtheorem{theorem}{Theorem}[section]
\newtheorem{lemma}[theorem]{Lemma}
\newtheorem{corollary}[theorem]{Corollary}
\newtheorem{proposition}{Proposition}[section]

\newtheorem{remark}[theorem]{Remark}

\title{\bf 
Waves in cosmological background with static Schwarzschild radius in the expanding universe
}

\author{{\bf Karen Yagdjian   } }

\begin{document}

\date{}

\maketitle

\thispagestyle{empty}
\begin{center}
{\small {School of Mathematical and Statistical Sciences,\\
University of Texas RGV, 
1201 W.~University Drive,  \\
Edinburg, TX 78539,
USA \\
e-mail: karen.yagdjian@utrgv.edu}}
\end{center}
\medskip
\vspace*{-0.6cm}
  
\addtocounter{section}{0}
\renewcommand{\theequation}{\thesection.\arabic{equation}}
\setcounter{equation}{0}
\pagenumbering{arabic}
\setcounter{page}{1}
\thispagestyle{empty}
\begin{abstract}

In this paper, we prove the existence of global in time small data solutions of  semilinear Klein-Gordon equations in  space-time with a static Schwarzschild radius in the expanding universe.

\end{abstract}

\section{The black hole in expanding universe. The   model with static Schwarzschild radius. Main results}

\setcounter{equation}{0}
\renewcommand{\theequation}{\thesection.\arabic{equation}}

The propagation of waves in the space-time of a single black hole and the partial differential equations describing them have been studied for quite a long time, and exhaustive answers to many interesting aspects of the problems, such as the linear stability of Schwarzschild black holes, decay of small solutions, Price's law, the formal mode analysis of the linearized equations, black hole shadow, particle creation, the ``John problem'', and the Strauss conjecture, are known. (See, e.g., \cite{Andersson-Blue-Wang, Bachelot-Nicolas,Birrell,Catania-Georgiev,Christ-Klainerman, Choquet-Bruhat_book, Daf_Rod_2005b, Daf_Rod_2005, Dafermos-Holzegel-Rodnianski,Dimock,Giorgi,Hawking,Hintz_2022_CMP,Hintz_2022,Lai-Zhou,Metcalfe,Nicolas, Ohanian-Ruffini,Parker,Perlick} and references therein.) In most publications on the partial differential equations in cosmological backgrounds, the black hole is assumed to be eternal; that is, the space-time and the Schwarzschild radius are assumed to be static. Actually, the latest astrophysical observational data confirm that the universe is expanding with acceleration  and that black holes are overwhelmingly present in the universe. The masses of the black holes are changing over time.  
Evidently, the metric of an expanding universe populated with black holes is extremely complicated. This gives rise to the question of how the propagation of waves in the cosmological background with black holes in the expanding universe can be mathematically reflected in the solutions of the related partial differential equations. We are motivated by the significant importance of the qualitative description of the solutions of the partial differential equations arising in cosmological backgrounds for understanding fundamental particle physics and the structure of the universe. In this paper, we focus on the equations of propagation of waves because the waves emitted by cosmic objects are one of the principal sources of empirical data in astrophysics. More precisely, we restrict ourselves to the case of a single black hole with a static Schwarzschild radius in the expanding universe and to the study of solutions of the linear and semilinear Klein-Gordon and wave equations with real and imaginary mass. The imaginary mass term appears in the Higgs boson equation \cite{CPDE_2010} and in the equation of tachyons \cite{Epstein-Moschella}.
 
To embed the black hole (the Schwarzschild  space-time) that has the line element
\[ 
\dsp ds^2= - \left( 1-\frac{2G M_{bh}}{c^2r} \right) \, c^2dt^2
+   \left( 1-\frac{2G M_{bh}}{c^2r} \right)^{-1}dr^2 +   r^2(d\theta ^2 + \sin^2 \theta \, d\phi ^2)
\]
in an expanded universe, we add the cosmological scale factor $a(t) $ to every component that is    measured in spatial linear units. Correspondingly, we write   the line element of such space-time  as follows:
\begin{equation}
 \label{ds_NM}
\dsp ds^2= - \left( 1-\frac{2G M_{bh}}{c^2a(t)r} \right) \, c^2dt^2
+ \left( 1-\frac{2G M_{bh}}{c^2a(t)r} \right)^{-1}a^2(t) dr^2 + a^2(t) r^2(d\theta ^2 + \sin^2 \theta \, d\phi ^2) .
\end{equation}
For the constant $M_{bh} $, this metric tensor solves Einstein's field equations with the diagonal energy-momentum tensor. 
Next, we take into account that black hole models with realistic behavior at infinity predict that the gravitating mass $M_{bh} $ of a black hole can increase with the expansion of the universe \cite{Farrah}. There are other various reasons (in \cite{Kastor}, ``This makes sense physically; ordinary matter would tend to accrete around the black hole.'' or, e.g., Hawking radiation) to believe that the mass $M_{bh} $ of the black hole (BH) is changing in time, that is, $ M_{bh} = M_{bh} (t) $. According to \cite{Farrah} ``Realistic astrophysical BH models must become cosmological at a large distance from the BH. Non-singular cosmological BH models can couple to the expansion of the universe, gaining mass proportional to the scale factor raised to some power $k$.''

The important characteristic of the BH is the so-called ``Schwarzschild radius''  $ \frac{2G M_{bh} (t)}{c^2a(t)} $. 
It was suggested and discussed in \cite{Macao},  the BH with the static (independent of time)   Schwarzschild radius embedded in the expanded universe, that is, $  \frac{d}{dt} a(t)>0$, meanwhile  
\[
R_{Sch}:=\frac{2G M_{bh} (t)}{c^2a (t)}=\frac{2G M_{bh} }{c^2 }, \quad \mbox{\rm where}\quad  M_{bh}=constant\,.
\]
The line element of this model (see also \cite[(4.117)]{Faraoni}) is given by 
\begin{equation}
\label{1.2}
\dsp ds^2= - \left( 1-\frac{2G M_{bh} }{c^2r} \right) \, c^2dt^2
+  \left( 1-\frac{2G M_{bh} }{c^2 r} \right)^{-1}a^2(t) dr^2 + a^2(t)r^2(d\theta ^2 + \sin^2 \theta \, d\phi ^2) .
\end{equation}
In \cite{Macao}, the case of the de~Sitter model with $a (t)= e^{Ht}$, where $H$ is the Hubble parameter,  was considered. 
It is worth mentioning  that   this model   solves the Einstein equation with the cosmological constant, and its 
     energy-momentum tensor   is of Type II (\cite[p. 89]{Hawking}). For  this model   the weak energy
condition is satisfied on some conic set consisting of the time-like vectors. In \cite{Macao}, the dominant energy condition (see, e.g.,
\cite[p. 91]{Hawking},  \cite[p. 51]{Choquet-Bruhat_book})  was addressed, the asymptotically dominant energy condition was defined and verified. Then it was discovered that an asymptotically strong energy condition was violated. 

Another model of   space-time  describing a black hole or massive object immersed in an expanding cosmological space-time is given by McVittie \cite{McVittie}. (For comprehensive discussion  and generalizations of the McVittie solution, see \cite[Ch.4]{Faraoni} and references therein.)   The dynamical many-black-hole space-times with well-controlled
asymptotic behavior as solutions of the Einstein vacuum
equation with positive cosmological constant under certain balance conditions on the black hole
parameters  are given by Hintz in  \cite{Hintz_2021}. 
\medskip

In this paper, we consider the Klein-Gordon equation for the self-interacting waves, that is,
\begin{eqnarray}
\label{KGE0}
\label{setCP}
 &  &
    \frac{\partial^2 \psi }{\partial t^2}
+ 3H       \frac{\partial \psi }{\partial t}
- e^{-2Ht}{\mathcal{A}}(x,\partial_x)\psi +   \frac{m^2c^4}{h^2}\psi +V(x,t) \psi=c^2\left( 1-\frac{2G M_{bh}}{c^2r} \right)\Psi(x,\psi)\,,
\end{eqnarray} 
where ${\mathcal{A}}(x,\partial_x) $ written in spherical coordinates is the following operator:
\begin{eqnarray}
\label{OpA}
{\mathcal{A}}(x,\partial_x)  
&  :=   &
 c^2\Bigg\{
    \left( 1-\frac{2G M_{bh}}{c^2r} \right)^2 \frac{\partial^2   }{\partial r^2} \\
&  &
\quad +\frac{2}{  r   }    \left( 1-\frac{ G M_{bh}}{c^2r} \right) \left( 1-\frac{2G M_{bh}}{c^2r} \right)\frac{\partial   }{\partial r}
 + \left( 1-\frac{2G M_{bh}}{c^2r} \right) \frac{1}{r^2}\Delta _{{\mathbb S}^2}   \Bigg\}     \,, \nonumber
\end{eqnarray}
while $ \Delta _{{\mathbb S}^2}$ is the Laplace operator on the unit sphere $ {\mathbb S}^2 \subset \R^3$. In (\ref{setCP}),   $V(r,t) $ is the potential that, in particular,  includes the case of the  gravitational potential $V(r,t)=-\frac{m^2c^2}{h^2}  \frac{2G M_{bh}}{ r}$. 
The term $ \Psi(x,\psi)$ represents the self-interaction of the field and vanishes for waves without self-interaction.
\medskip

 Lemma~\ref{L_Cos_Pr} shows that the damping term of the covariant Klein-Gordon equation  is independent of time and location     in the background (\ref{ds_NM}) only if the mass $M_{bh}=M_{bh}(t)$ of BH is proportional to the scale factor. This is where the cosmological principle (see, e.g., \cite[Sec. 9.1]{Ohanian-Ruffini}) comes into play.  
\medskip

We analyze the waves by appealing to the integral transform approach developed in \cite{Macao,MN2015}. It turns out that,  due to those   integral transforms, it is possible to reduce    the problem with infinite time to the problem with finite time  and to apply an energy estimate for the finite time,   thus eliminating the  growth of the energy. Moreover, this allows us to avoid the severity of the global construction of the phase function, which is one of the challenges of  micro-local analysis. 
We must    emphasize  that this is possible since the de~Sitter space-time has a permanently bounded domain of influence. 
\medskip

The  
covariant  Klein-Gordon equation in the black hole with static  Schwarzschild radius  
embedded in the de~Sitter universe, that is, in the metric (\ref{1.2}) with $a(t)=e^{Ht}$, 
can be written in the Cauchy-Kowalewski  form as follows:
\begin{eqnarray}
\label{KGE}
 &  &
    \frac{\partial^2 \psi }{\partial t^2}
+ 3 H      \frac{\partial \psi }{\partial t}
-e^{-2Ht}{\mathcal{A}}(x,\partial_x)\psi+   \frac{m^2c^4}{h^2} \left( 1-\frac{2G M_{bh}}{c^2r} \right)\psi =c^2\left( 1-\frac{2G M_{bh}}{c^2r} \right)\Psi(\psi)\,.
\end{eqnarray}
The term $\frac{m^2c^4}{h^2} \left( 1-\frac{2G M_{bh}}{c^2r} \right) $ can be split  into the ``rest mass'' 
term $ \frac{m^2c^4}{h^2} $ and the gravitational (Newtonian) potential part $-\frac{m^2c^2}{h^2}  \frac{2G M_{bh}}{ r} $.
The equation (\ref{KGE}) can be  regarded as an addition of the gravitational  potential  
\begin{equation}
\label{Potential}
 V(r)=   -   \frac{m^2c^2}{h^2} \frac{2G M_{bh}}{ r}
\end{equation}
to the equation 
\begin{eqnarray}
\label{1.3}
 &  &
    \frac{\partial^2 \psi }{\partial t^2}
+ 3H        \frac{\partial \psi }{\partial t}
-e^{-2Ht}{\mathcal{A}}(x,\partial_x)\psi+   \frac{m^2c^4}{h^2} \psi=c^2\left( 1-\frac{2G M_{bh}}{c^2r} \right)\Psi(\psi)\,.
\end{eqnarray}

The Klein-Gordon equation in the cosmological background with a static Schwarzschild radius in the expanding 
universe is (\ref{KGE0}). 
In this paper, we consider waves (solutions of the equations) in the exterior of the black hole denoted $B^{ext}_{Sch}:=\{x \in \R^3\,|\, |x| > R_{Sch} \}$. 
Bearing in mind the gravitational potential, we relate the properties of the potential $V(x,t)\in C^2 (B^{ext}_{Sch}\times[0, \infty))$ to the setting of the Cauchy problem for semilinear equation, more precisely, with   the support of the initial functions.  Denote $H_{(s )}:= H_{(s )}( \R^3)$ the Sobolev space, while ${\mathcal B}^\infty(B^{ext}_{Sch}\times[0,\infty)) $ is a set of all smooth functions with uniformly bounded  derivatives of any order. We also denote  $\pi_x$ the projection operator on $\R^3$.
\medskip

\noindent
{\bf The Cauchy Problem.} Let the number $R_{ID} $ be such that $R_{ID} > R_{Sch}$. 
For every given $\psi_0  ,  \psi_1 \in H_{(s)}   $ with supp\,$\psi_0\, \cup$ 
supp\,$\psi_1  \subseteq \{ x\in  \R^3\,|\, |x| \geq  R_{ID}\} \subset B^{ext}_{Sch}$,  find 
 {\sl a global  in time solution}  $\psi \in C^1([0,\infty); H_{(s)})$  of the equation (\ref{setCP}), such that supp~$\psi (t) \subseteq B^{ext}_{Sch}$ for all $t >0$,  
and which takes the  initial values
\begin{eqnarray}
\label{CP}
&   &
\psi (x,0) =  \psi_0 (x),\quad \psi_t (x,0) =  \psi_1 (x), \quad \mbox{\rm for all} \quad x \in \R^3.
\end{eqnarray} 

\ndt
 {\bf Condition ($\mathcal V$)} on the potential: {\it $V(x,t)\in {\mathcal B}^\infty (B^{ext}_{Sch}\times[0, \infty))$ and for given $s$ there is $\varepsilon_0>0 $ such that  
$
 \|V (x,t) \Phi (x) \| _{H_{(s)}}\leq \varepsilon_0 \| \Phi   \| _{H_{(s)}}$   for all $t \in [0,\infty)$ and all $ \Phi \in H_{(s)}  \, , 
$
such that  compact supp~$\Phi  \subset B^{ext}_{Sch}$.} 
\medskip
 
 \ndt
 The non-linear term  is supposed to satisfy the following condition. \medskip

 \ndt
{\bf Condition ($\mathcal L$).} {\it The smooth in $x \in B^{ext}_{Sch}$ function $\Psi =\Psi (x,\psi )$ is
said to be Lipschitz continuous with exponent $\alpha \geq 0 $ in the space $H_{(s)} $   if supp~$\Psi (x,\psi )\subseteq $ supp~$ \psi $ and there is a constant 
$C \geq 0$  such that}
\[
 \|  \Psi (x,\psi  _1 (x))- \Psi (x,\psi  _2(x) ) \|_{H_{(s)} }  \leq C
\| \psi  _1 -  \psi _2   \|_{H_{(s)} }
\left( \|  \psi _1  \|^{\alpha} _{H_{(s)} }
+ \|  \psi _2   \|^{\alpha} _{H_{(s)} } \right) \,\, \mbox{\it for all} \,\,  \psi  _1, \psi  _2 \in  H_{(s)}\,.
\]
The interesting cases of the semilinear term are $\Psi(\psi)=|\psi|^{1+\alpha}$ and $ \Psi(\psi)=\psi|\psi|^{\alpha}$ in (\ref{setCP}).

We say that {\it an equation has a large mass} if $m^2 \geq \frac{9H^2h^2}{4 c^4} $. First, we consider the case of  large mass and   
the Cauchy problem in the Sobolev space $ H_{(s)}$ with $ s > 3/2$, which is an algebra.
Define  the metric space
\[
X({R,H_{(s)},\gamma})  := \left\{ \psi  \in C([0,\infty) ; H_{(s)} ) \; \Big|  \;
 \parallel  \psi   \parallel _X := \sup_{t \in [0,\infty) } e^{\gamma t}  \parallel
\psi  (x ,t) \parallel _{H_{(s)}}
\le R \right\}\,,
\]
where $\gamma \in {\mathbb R}$, with the metric
$d(\psi _1,\psi _2) := \sup_{t \in [0,\infty) }  e^{\gamma t}  \parallel  \psi _1 (x , t) - \psi _2 (x ,t) \parallel _{H_{(s)}}$.

\begin{theorem}
\label{T_large_mass}
Consider the Cauchy problem for the equation (\ref{setCP}) 
in $  \R^3\times [0,\infty) $  with the initial conditions
\begin{equation}
\label{13.12}
\psi (x,0)=\psi _0(x)\in H_{(s)},\quad \partial_t \psi (x,0)=\psi _1(x)\in H_{(s)}    \,,
\end{equation}
where 
\begin{equation}
\label{13.13}
\mbox{\rm supp} \,\psi _0, \mbox{\rm supp} \,\psi _1 \subseteq \left\{ x \in \R^3\,\left|\, |x|\geq R_{ID}> \frac{c}{H}+ R_{Sch} \right\} \subset  B^{ext}_{Sch}.  \right.
\end{equation}
Assume that   the potential $V(x,t)\in {\mathcal B}^\infty (B^{ext}_{Sch}\times[0, \infty))$ satisfies the condition 
($\mathcal V$). 
Assume also that the physical mass $m $ of the field is large, that is, 
$
 \frac{m^2c^4}{h^2}  \geq \frac{9H^2}{4 }\,, 
$
and the nonlinear term $ \Psi(x,\psi)$ satisfies condition $({\mathcal L})$ with $\alpha >0$ and $\Psi(x,0) =0$. 

If $\varepsilon_0$ and the norms 
 $
 \|\psi _0\|_{H_{(s)}} , \|\psi _1\|_{H_{(s)}}  $ 
with $s> 3/2 $ are sufficiently small, then the Cauchy problem (\ref{setCP})\&(\ref{13.12})\&(\ref{13.13})  has a global solution
\[
\psi \in C^2([0,\infty);H_{(s)})\,.
\] 
Moreover, the solution $ \psi$ belongs to the space $X({R,H_{(s)},\gamma})  $, $ \gamma \in (0,H)$,  that is, the solution  $ \psi$ decays according to 
\[
  \parallel
\psi  (x ,t) \parallel _{H_{(s)}}
\le R e^{-\gamma t} , \quad  t \in [0,\infty)\,.
\] 
 If $\frac{m^2c^4}{h^2}  > \frac{9H^2}{4 }$ or $\psi _0=0 $, then $\gamma=H$.
\end{theorem}

Next, we assume that $m \in \C$ and define
\[
M:=  \left(\frac{9H^2 }{4 }-\frac{m^2c^4}{h^2} \right)^{1/2} \,.
\]
Then we consider the case of {\it small mass}, that is,   $\Re M >0 $. 
Here, $iM $ can be regarded as a {\it curved mass}. The Higgs boson equation and the equation of tachions have small masses. 
\begin{theorem}
\label{TIE}  
Assume that  $\Psi (x,\psi)$  is Lipschitz continuous in the  space $H_{(s)} $, $ s > 3/2$, $\Psi(x,0)=0$, and that $\alpha >0 $. 
Assume that   the potential $V(x,t)\in {\mathcal B}^\infty  (B^{ext}_{Sch}\times[0, \infty))$ satisfies condition ($\mathcal V$).

\noindent 
 \mbox{\rm (i)} Suppose that  $0< \Re M<H/2$.
  Then for every given functions $  \psi _0(x  ), \psi _1(x  ) \in  H_{(s)}  $  such that    
\begin{equation}
\label{smalldata} 
   \| \psi _0   
\|_{H_{(s)}}
+   \|\psi _1  
\|_{H_{(s)}}  < \varepsilon\,,  
\end{equation}
and for sufficiently small  $\varepsilon $, $\varepsilon_0 $,    the Cauchy problem  (\ref{setCP})\&(\ref{13.12})\&(\ref{13.13}) has a global solution  \\ $ \psi \in   C^2([0,\infty);H_{(s)}) $. For the solution  with $\gamma \in [0,H) $, one has  
\begin{equation}
\label{3.2}  
\sup_{t \in [0,\infty)}  e^{\gamma  t}  \| \psi  (\cdot ,t) \|_{H_{(s)} }  < 2\varepsilon \,.
\end{equation}
If $V(x,t)=0$, then $\gamma $ can be chosen as  $\gamma =H $.  
For $\Re M<\frac{H}{2} $   and $\gamma  < \frac{H}{1+\alpha}$,   the  norm of  $ \partial_t \psi  $ decays as follows,
\begin{eqnarray} 
\label{der_t_i}
\| \partial_t \psi(t,x)\|_{H_{(s-1)} }  
 & \leq   &
\cases{ C\varepsilon e ^{- \gamma t}     ,
\cr
   C\varepsilon e ^{- \gamma  (1+\alpha) t }     ,
  \quad \mbox{\rm if}  \quad V=0.
}
\end{eqnarray} \\
\noindent
 \mbox{\rm (ii)} Suppose that  $\Re M \in (H/2,3H/2)$ or $M= H/2$.  Then for every given functions $  \psi _0(x  ), \psi _1(x  ) \in  H_{(s)}  $  such that    (\ref{smalldata}),  
 for every $\gamma  $, $\gamma <(3H/2-\Re M)/(\alpha +1) $, and for sufficiently small  $\varepsilon_0 $, $\varepsilon $,  the Cauchy problem  (\ref{setCP})\&(\ref{13.12})\&(\ref{13.13})   has a global solution \, $  \psi (x ,t) \in C^2([0,\infty);H_{(s)}) $. For the solution, the inequality (\ref{3.2}) is fulfilled. 
 
For $\Re M \in (H/2,3H/2)$ or $M= H/2$,    $\gamma <(3H/2-\Re M)/(\alpha +1) $,   and $ \delta>  \Re M- H/2+ \gamma  (1+\alpha)  $, the norm of  $ \partial_t \psi  $  satisfies 
\begin{eqnarray} 
\label{der_t_ii}  
\| \partial_t \psi(t,x)\|_{H_{(s-1)} }   
 &  \leq & 
\cases{C \varepsilon e ^{(\delta- \gamma) t}  \cr
 C\varepsilon e ^{(\delta-\gamma (1+\alpha))  t}        ,
  \quad \mbox{\rm if} \quad V=0  . } 
\end{eqnarray} 
\noindent
\mbox{\rm (iii)} Suppose that  $\Re M>3H/2$. Then for  every given functions $  \psi _0(x  ), \psi _1(x  ) \in  H_{(s)}  $  such that    (\ref{smalldata}),  and every $\gamma$,  
$  \gamma < ( {3H}/{2} - \Re M)/(\alpha +1)  $,    
the  solution \, $ \psi (x ,t)   $ of the problem  (\ref{setCP})\&(\ref{13.12})\&(\ref{13.13})  has   
the lifespan $T_{ls}$ that can be estimated from below by 
\[
   T_{ls}  
 \geq  -\frac{1}{\gamma}\ln \left(    \varepsilon  \right)- C(\alpha,\gamma,\varepsilon _0,H, M )
\] 
 with some number $C(\alpha,\gamma,\varepsilon _0,H, M )$.
\end{theorem}
\medskip

\noindent
{\bf Examples of the potential.} 
(i) For the case of gravitational potential (\ref{Potential}),  
the conditions of the theorems  imply that    
$\frac{m^2c^2}{h^2} \frac{2G M_{bh}}{R_{ID}} $ is sufficiently small. \\

\noindent
(ii) In the case of general time-dependent potential,  we can assume that 
$\sup_{r\geq 2G M_{bh}/c^2 ,\, t \in [0,\infty)}|V(x,t)|$  is sufficiently small. \\

\noindent
(iii) If we consider the case of time-dependent potential
$
V(x,t)= 
-   m_H^2c^4 h^{-2} e^{- 2Ht }$, $m_H=const > 0$,
 then the equation  (\ref{KGE})  leads to
\begin{eqnarray*}
 &  &
    \frac{\partial^2 \psi }{\partial t^2}
+ 3H       \frac{\partial \psi }{\partial t}
- e^{-2Ht} 
    \left( {\mathcal{A}}(x,\partial_x) +  \frac{m_H^2c^4}{h^2} \right) \psi  +   \frac{m^2c^4}{h^2} \psi =c^2\left( 1-\frac{2G M_{bh}}{c^2r} \right)\Psi(\psi)\,.
\end{eqnarray*}
In this case, in the application of the integral transform approach, one can appeal to the results of \cite{Bachelot-Nicolas,Dimock,Nicolas}.

\noindent
(iv) The Yukawa potential \cite{Griffiths} $
V(r,t)= -g^2r^{-1} e^{-\alpha m r}  $, 
 which is a  model for the binding force 
in an atomic nucleus.  
 Here $ m$ is the mass of the particle, $g$   is the amplitude of potential, $ \alpha$  is a  scaling constant,  and $1/(\alpha m)$ its range.
\\
\noindent
(v)  
If we consider the equation with the distributed mass term, that is, replace  (see \cite[p. 51]{Frankel})
$
m^2 \longmapsto  {m^2}/{\left( 1-\frac{2G M_{bh}}{c^2r} \right)} 
$, 
then we  arrive at equation (\ref{1.3}) without potential.\\

Another interesting and important model with $M_{bh}(t) \sim a^3(t)$ (see  \cite{Farrah}) will be discussed in the forthcoming paper.
\begin{remark}
Due to the scale factor $a(t)=e^{Ht}$,  the equation (\ref{setCP}) has multiple characteristics at $t=\infty$; this is reflected in the choice of initial functions   $  \psi _0(x  ), \psi _1(x  ) \in  H_{(s)}  $ with the same $s$ for both functions.  In fact, the structures of Fourier integral operators in (\ref{psi}) generating solutions to the linear equation via each initial function  coincide.
\end{remark}

\begin{remark} 
In  light  of Corollary~\ref{C2.3}, it will be interesting to relax the condition $ R_{ID}> {c}/{H}+ R_{Sch}$   
 of (\ref{13.13}). 
\end{remark}

\begin{remark}
The decay of energy in the case of a large mass can be considered by classical methods (see, e.g., \cite{NA2015,Nakamura_JMAA2014})
and will be done in the forthcoming paper.
\end{remark}

\begin{remark}
It will also be interesting  to combine the integral transform approach with the results on the Cauchy problem 
for the  linear Klein-Gordon equation on Schwarzschild-like metric
$
    \frac{\partial^2 v }{\partial t^2}
-      {\mathcal A}(x,\partial_x)v +  \frac{m_H^2c^4}{h^2}v  =0,
$ 
obtained in \cite{Bachelot-Nicolas,Nicolas}. 
   The term with ${m_H^2c^4}/{h^2}>0$ can be 
regarded as  potential due to the expansion. (Compare with  the condition  $\xi>0$ of \cite{Bachelot-Nicolas,Nicolas}, that is,  $m_H>0 $, that is crucial for the results of \cite{Bachelot-Nicolas,Nicolas}.)
\end{remark}

\begin{remark}
If we carry out the above-mentioned derivation  in the case of the Majumdar-Papapetrou  multi-black-hole solutions (see, e.g., \cite{Hintz_2021,Kastor}) of the Einstein equation and assume, in accordance with \cite{Farrah},   that every black hole has a static Schwarzschild radius, then we arrive at a  similar picture of the propagation of the waves in the expanding universe. This will be done in the forthcoming paper. 
\end{remark} 
\medskip

\ndt
{\bf Outline of the proof.}
The integral transform approach \cite{MN2015} applied to the initial value  problem (\ref{CP}) for the equation 
\[ 
    \frac{\partial^2 \psi }{\partial t^2}
+ 3H        \frac{\partial \psi }{\partial t}
-e^{-2Ht}{\mathcal{A}}(x,\partial_x)\psi +   \frac{m^2c^4}{h^2} \psi=f\,,
\] 
leads to the following formula (the kernels $ E(r,t;0,b;M)$, $K_0( s,t;M) $, and $K_1( s,t;M) $, are defined in (\ref{0.6}), (\ref{K0MH}), and (\ref{K1MH}), respectively) for the solution
\begin{eqnarray}
\label{psi} 
\psi (x,t) 
&  =  &
e ^{-\frac{3}{2}Ht}  2   \int_{ 0}^{t} db
  \int_{ 0}^{\phi (t)- \phi (b)} e ^{\frac{3H}{2}b}  E(r,t;0,b;M)  v_f(x,r ;b) \, dr  
+ e ^{-Ht} v_{\psi_0}  (x, \phi (t))   \\
&  &
+ \, e ^{-\frac{3}{2}Ht}\int_{ 0}^{\phi (t)}  \left[ 2K_0( s,t;M)    
+   3HK_1( s,t;M)  \right] v_{\psi_0 } (x,  s)
  ds  \nonumber \\
&  &
+\, 2e ^{-\frac{3}{2}Ht}\int_{0}^{\phi (t) }  v_{\psi_1 } (x,  s)
  K_1( s,t;M)   ds,\quad  x \in \Omega \subseteq \R^3 \,,\,\, t \in I=[0,T] \subseteq [0,\infty)\,, \nonumber 
\end{eqnarray}
where $0<T\leq\infty$, $\phi (t):= (1-e^{-Ht} )/H$, $M^2 =  \frac{9 H^2}{4 } -\frac{m^2c^4}{h^2}$, and  
$v_{f} (x,  s)$ is a solution of 
\begin{equation}
\label{vf}
\cases{ 
 v_{tt} -   {\mathcal A}(x,\partial_x)  v  =  0 \,, \quad x \in \Omega \,,\quad t \in [0,(1-e^{-HT})/H]\,,\cr
v(x,0;b)=f(x,b)\,, \quad v_t(x,0;b)= 0\,, \quad x \in \Omega,\quad b \in I\,, 
}
\end{equation}
while the function \, $  v_\varphi(x, t) \in C_{x,t}^{m,2}(\Omega \times [0,(1-e^{-HT})/H])$ \, is a   solution   of the   problem
\begin{equation}
\label{vphi}
\cases{
 v_{tt}-   {\mathcal A}(x,\partial_x)  v =0, \quad x \in \Omega \,,\quad t \in [0,(1-e^{-HT})/H]\,, \cr
 v(x,0)= \varphi (x), \quad v_t(x,0)=0\,,\quad x \in \Omega\,.  
} 
\end{equation} 
One can regard that integral transform as an analytical mechanism  that,  from the
massless field in the static BH space-time, generates massive particles in the space-time of the BH in the expanding universe.
Considerations of geodesics in the black hole space-time (see, e.g., \cite[Sec.19,20]{Chandrasekhar} and  \cite[Ch.18]{Taylor}) show  that the  $\pi_x($supp~$v_{\psi  } (x,  s))$ is compact  for all $s \in [0,1/H] $   and on the positive distance from the event horizon $r=  R_{Sch} $ if the distance  dist((supp$\,\psi_0 \cup$  supp$\,\psi_1), \ol{S_{R_{Sch}}(0)} $), that is, $R_{ID}$,  
is  sufficiently large. 
But even for the initial data without the last restriction on the supports, the function (\ref{psi}) solves the equation as long as the functions $v_{\psi_0} $,  $v_{\psi_1} $,  and $ v_f$ are defined.

If one applies the Liouville transform  with  $u=e^{\frac{3H}{2 }t} \psi $  
to the   Klein-Gordon  equation,  then the covariant   Klein-Gordon  equation with the source   $f$ becomes  
\begin{eqnarray*}
&   &
     \frac{\partial^2 u }{\partial t^2}
- e^{-2Ht}{\mathcal{A}}(x,\partial_x)u
+   \frac{m^2c^4}{h^2} u-\frac{9H^2}{4 }        u  +V(r,t)u=g \,,
\end{eqnarray*}
where $ g=e^{\frac{3H}{2 }t}f $. 
 This is the non-covariant   Klein-Gordon equation with the ``imaginary mass''
\begin{eqnarray*}
&  &
u_{tt} -  e^{-2Ht}{\mathcal{A}}(x,\partial_x)  u - M^2 u  +V(r,t)u=  g \,,
\end{eqnarray*}
where the   mass term is
$
M^2 =  \frac{9 H^2}{4 } -\frac{m^2c^4}{h^2}$.

Thus, we are in a position to apply Theorem~2.1~\cite{Macao} and to reveal the properties of the black hole in the de~Sitter background.
The treatment of the semi-linear equation is based on Banach's fixed point theorem and on the estimates for the solution of the linear equation.

\section{Preliminaries. Linear equation}
\label{S_FPS}

The next lemma shows that  the  space-time of the BH  with a static Schwarzschild radius is the only space-time that dissipates   the waves   independently of time  and spatial coordinates. 
\begin{lemma}
\label{L_Cos_Pr}
Consider the d'Alembert operator in the metric  (\ref{ds_NM}) with depending on the time mass of BH, $M_{bh}=M_{bh}(t)$, and a positive function $a(t)\in C^1([0,\infty))$. The only function $M_{bh}(t)\in C^1([0,\infty))$ that makes the damping term of the d'Alembert operator independent of spatial variables   is $M_{bh}(t)= const \cdot a(t) $.  Moreover, for the de~Sitter  scale factor  $a(t)=e^{Ht}$, the damping term is independent of time as well. 
\end{lemma}

\noindent
{\bf Proof.} The damping term of d'Alembert operator is the ratio of the coefficients of $\psi_t$ and $\psi_{tt} $. The derivative of that ratio is 
\[
\frac{\partial}{\partial r} \left(\frac{3 c^2 r a(t) a'(t)-8 G M(t) a'(t)+2 G a(t) M'(t)}{a(t) \left(c^2 r a(t)-2 G M(t)\right)}\right)=\frac{2 c^2 G \left(M(t) a'(t)-a(t) M'(t)\right)}{\left(c^2 r a(t)-2 G M(t)\right)^2}\,,
\]
which vanishes only if $M_{bh}(t)= const \cdot a(t) $,     and the statement follows from (\ref{setCP}). The lemma is proved.
\qed

\subsection{Equation in Cartesian coordinates. Finite propagation speed}

\setcounter{equation}{0}
\renewcommand{\theequation}{\thesection.\arabic{equation}}

When no ambiguity arises, we will use the notations $\vec{ x}=(x_1,x_2,x_3):=(x,y,z)$ and $\vec{\xi }=(\xi _1,\xi _2,\xi _3)$. The scalar  product in $\R^3$ will be denoted $\vec{ x}\cdot \vec{\xi }$. 
In  the case of $H=0$, the linear Klein-Gordon equation without a source in Cartesian coordinates can be written as follows:
\begin{eqnarray}
\label{EqCart}
 &  &
   \psi_{tt}(\vec{x},t)
-c^2\Bigg\{
    F(|\vec{x}|)^2  \frac{1}{|\vec{x}|^2} \sum_{k ,\ell=1,2,3}x_k x_\ell   \psi _{x_k x_\ell} (\vec{x},t)  
+    \left( 1-\frac{ G M_{bh}}{c^2|\vec{x}|} \right) F(|\vec{x}|) \frac{2}{|\vec{x}|^2}  \sum_{k  =1,2,3}x_k \psi_{x_k}(\vec{x},t)   \nonumber \\
 &  & 
 + F(|\vec{x}|) \frac{1}{|\vec{x}|^2}\Big[ \sum_{k  =1,2,3}    |\vec{x}|^2   \psi _{x_kx _k} (\vec{x},t) - \sum_{k ,\ell =1,2,3 }x_k x_\ell   \psi _{x_k x_\ell} (\vec{x},t)  
-2\sum_{k  =1,2,3} x_k \psi _{x_k} (\vec{x},t) \Big]  \Bigg\} \nonumber\\
&  &
+    \frac{m^2 c^4 }{h^2}  \psi+V(\vec{x},t)\psi =0,  
\end{eqnarray}
where
\[
F(|\vec{x}|)=F(r):= 1-\frac{2G M_{bh}}{c^2|\vec{x}|}= 1-\frac{R_{Sch}}{ |\vec{x}|},\quad r:=|\vec{x}|:=\sqrt{x^2+y^2+z^2}, \quad r >R_{Sch}\,.
\]
Thus, the symbol $ {\mathcal A}(\vec{ x};\vec{\xi } )$ of the operator $ {\mathcal A}(\vec{ x},\partial_x)$ (\ref{OpA})  is given by
\begin{eqnarray}
\label{SymbA}
{\mathcal A}(\vec{ x};\vec{\xi } )
 & = &
  A_2(\vec{ x};\vec{\xi } )+A_1(\vec{ x};\vec{\xi } )\,,
\end{eqnarray}
where $ A_2(\vec{ x};\vec{\xi } )$ and $A_1(\vec{ x};\vec{\xi } )$ are the principal symbol and the low order 
symbol, respectively, and 
\[
 A_2(\vec{ x};\vec{\xi } ) 
  =   
  -c^{2}\left(1-\frac{R_{Sch}}{ |\vec{ x}|} \right)\left( | \vec{\xi }|^2 -\frac{ R_{Sch} \left(\vec{ x}\cdot \vec{\xi }\right)^2}{ |\vec{ x}|^3 }\right) ,\quad 
A_1(\vec{ x};\vec{\xi } ) 
  =  -c^{2}\left(1-\frac{R_{Sch}}{|\vec{ x}|} \right)\frac{ i R_{Sch}  \left(\vec{ x}\cdot \vec{\xi }\right)}{ |\vec{ x}|^3 }.
\]
Consider the zeros of the principal symbol of the equation, that is, solutions to 
\[ 
\tau^2 
-
c^{2}\left(1-\frac{R_{Sch}}{ |\vec{ x}|} \right) |\vec{ \xi }|^2\Bigg( 
 1  -\frac{R_{Sch}}{|\vec{ x}|} \frac{\left(\vec{ x}\cdot \vec{\xi }\right)^2}{ |\vec{ x}|^2|\vec{\xi} |^2} 
\Bigg)=0  \,.
\]  
It is evident that for $  |\vec{\xi} |=1$ and   $|\vec{x}|  > R_{Sch}$, we have
\[
|\tau|^2 
  \leq   
c^{2}\left(1-\frac{R_{Sch}}{|\vec{ x}|} \right)^2 \,.
\]
Thus, the equation (\ref{EqCart}) is strictly hyperbolic on every compact set in $B_{Sch}^{ext}\times \R $ 
and has multiple characteristics on the sphere  $r=R _{Sch}$. This indicates the behavior of the light cone approaching the event horizon $r=R_{Sch}$.
Since the operator has  multiple characteristics, the well-posedness of the Cauchy problem requires some kind of Levi condition. (See, for detail,
\cite{YagBook}.)  
In the interior of BH, the operator is not hyperbolic in  the direction of time, but it is hyperbolic in the radial direction.

Consider  the null radial geodesies of the space-time (\ref{1.2}) when it has the permanently restricted domain of influence, that is, 
\[
a(t) >0,  \quad \frac{d}{dt}a(t)>0 \quad A(t) :=\int_0^t \frac{1}{a(s)} ds \leq A(\infty)< \infty \quad \mbox{\rm for all} \quad t \in [0,\infty). 
\]
More exactly, consider the geodesic solving 
\begin{equation}
\label{r_geodesic}
\frac{dr}{dt}=-c   \frac{1}{a(t)}  \left( 1-\frac{R_{Sch}}{r} \right) 
\end{equation}
and starting at $R_{ID}$, that is $r(0)=R_{ID}$. The existence of global in time  geodesic is given by the following statement.
\begin{lemma}
For $ R_{ID}>R_{Sch} $ there is a positive number $\varepsilon$ such that the implicit function  $r=r(t)$ given by the equation
\begin{equation}
\label{r_geod_form}
 R_{ID}-r -R_{Sch}  \ln \left( 1- \frac{ R_{ID} - r}{R_{ID}-R_{Sch}}  \right)   =cA(t)
\end{equation}  
is well defined for all $t \in [0,\infty)$ and satisfies the inequality 
\begin{equation}
\label{r_support}
r(t) \geq R_{Sch}+\varepsilon \quad \mbox{ for all} \quad t >0\,.
\end{equation}
\end{lemma}

\ndt
{\bf Proof.}  By solving equation (\ref{r_geodesic}), we arrive at the formula (\ref{r_geod_form}). 
Consider an implicit function $z(\tau)$ given by the equation
\[
z -R_{Sch}  \ln \left( 1- \frac{ z}{R_{ID}-R_{Sch}}  \right)   =\tau 
\] 
such that $z(0)=0$, and $\tau \in [0,  \infty)  $. 
This function is well defined if $z(\tau) \in [0, R_{ID}-R_{Sch})$,  is positive, and is continuous. Indeed, 
\[
\frac{dz}{d\tau}\left(1 + R_{Sch}   \frac{ 1}{R_{ID}-R_{Sch}-z }   \right)   =1
\]
as long as $z(\tau)< R_{ID}-R_{Sch}$ implies $\frac{dz}{d\tau}>0 $ and  the function $z=z(\tau)$ is well defined. Denote $z_1=z(cA(\infty))$, then 
$z(\tau) \leq z_1$ for all $\tau \in [0,cA(\infty)]$. 
The number $z_1\in [0, R_{ID}-R_{Sch})$ exists since 
\[
\lim_{z  \nearrow  R_{ID}-R_{Sch} } \left(  z -R_{Sch}  \ln \left( 1- \frac{ z}{R_{ID}-R_{Sch}}  \right) \right)  = \infty \,.
\] 
Consequently, there is a positive number $\varepsilon$ such that 
$z(\tau) \leq R_{ID}-R_{Sch} -\varepsilon$   for all $\tau \in [0,cA(\infty)]$. 
The proper time $\tau$  is defined by
\[
 \tau (t):=cA(t)\quad \mbox{\rm and} \quad \frac{d\tau}{dt}= c\frac{1}{a(t)} >0 \quad \mbox{\rm for all} \quad t \in [0,\infty)\,.
\]
 Consider the function $r=r(t)$, which is  defined on $[0,\infty) $ and  
$r(t)= R_{ID}-z(\tau (t))$,   where $\tau (t)< cA(\infty)$   for all $t \in [0,\infty)$. 
The inequality (\ref{r_support}) is proved. \qed

\begin{corollary}
\label{C2.3}
If $H>0$, then for every compact $K\subset  B_{Sch}^{ext}\subset \R^3$, dist$(\partial K,B_{Sch})>c/H$,  if the initial functions have compact supports  in $K$, then there is $\varepsilon>0$ such that the solution  of (\ref{setCP}) with $\Psi=0$  has a compact support in $ \{x \in \R^3 \,|\, |x| > R_{Sch}+\varepsilon\} \subset B_{Sch}^{ext}$ for all   $t \in [0,\infty)$. 
\end{corollary}

\ndt
{\bf Proof.} This follows from the dependence domain  Theorem~4.10.1~\cite{YagBook} and Theorem~6.10~\cite{Mizohata}. \qed
\medskip

If $H=0$, then for every compact $K\subset  B_{Sch}^{ext}\subset \R^3$, if the initial functions have compact supports  in $K$   and $d:=$dist$(K;\{\sqrt{x^2+y^2+z^2} \leq R_{Sch}\})>0$, then the solution of (\ref{setCP}) with $\Psi=0$  has a compact support in $ B_{Sch}^{ext}$ for all time $t \in [0,d/c)$. 
\medskip
 
In order to derive $H_{(s)}$-estimates, we use some auxiliary  operator defined as follows. For every given fixed compact $K\subset  B_{Sch}^{ext}$ the coefficients can be continued smoothly outside of the small $\delta$-neighborhood of $K$ to be constants. This continuation does not affect   solutions with initial data supported by $K$ at least for time duration $d/c$. Thus, one can replace the operator ${\mathcal{A}}(x,\partial_x) $ with its continuation (auxiliary operator); we will do it without special notification.  More precisely, we define the
{\sl auxiliary operator} ${\mathcal A}_{\varepsilon }(\vec{ x};\vec{\xi } )$ 
as an operator   with the symbol
\begin{eqnarray}
\label{AUX} 
{\mathcal A}_{\varepsilon }(\vec{ x};\vec{\xi } )  
& := &
c^{2}\left(1-\chi_\varepsilon (\vec{ x})\frac{R_{Sch}}{  |\vec{ x}|}\right)\left(-|\xi |^2    +\chi_\varepsilon (\vec{ x})\frac{R_{Sch} \left(\vec{ x}\cdot \vec{\xi }\right)^2}{  |\vec{ x}|^{3} } 
-\chi_\varepsilon (\vec{ x})\frac{ i R_{Sch}  \left(\vec{ x}\cdot \vec{\xi }\right)}{    |\vec{ x}|^{3} }\right), 
\end{eqnarray}
where  $\chi_\varepsilon $ is a cutoff function vanishing if $|\vec{ x}| < R_{Sch}+\varepsilon/2 $ while $\chi_\varepsilon (\vec{ x})=1$ if $|\vec{ x}| > R_{Sch}+\varepsilon  $.

\subsection{Linear equation. Energy estimates and   energy conservation}

Let $\partial_t^2- A(x,\partial_x) =\partial_t^2- \sum_{|\alpha |\leq 2} a_{\alpha }(x)\partial _x^\alpha $ be a  second-order strictly hyperbolic   operator with coefficients $a_{\alpha } \in {\mathcal B}^\infty $, where   ${\mathcal B}^\infty $  is the space of all $C^\infty ({\mathbb R}^3)$ functions with uniformly bounded derivatives of all orders.  
 Let $v=v(x,t) $ be the solution of the problem
\[
\cases{
\partial_t^2 v - A(x,\partial_x)v =0, \quad x \in {\mathbb R}^n, \quad t \geq 0, \cr 
v(x,0) = v_0 (x),\quad   v_t (x,0)= v_1(x),\quad x \in {\mathbb R}^n\,.}
\]
 The following energy estimate is well known. (See, e.g., \cite{Taylor}.) For every $s \in {\mathbb R}$ and given $T>0$ there is  $ C_s(T)$ such that  
\begin{eqnarray}
\label{Blplq}
&   &
 \|v_{t}(t)\|_{H_{(s)}}   +  \|v (t)\|_{H_{(s+1)}} \leq  C_s(T) ( \|v_1\|_{H_{(s)}} +  \|v_0 \|_{H_{(s+1)}} )  , \quad 0 \leq  t \leq T  \,.
\end{eqnarray}

We note that although in this estimate the time interval is bounded, however, due to the integral transform approach given in \cite{MN2015}, it is possible to reduce    the problem with infinite time to the problem with  finite time and to apply (\ref{Blplq}). 
In fact,   this is possible since   de~Sitter space-time in FLRW coordinates has a permanently bounded domain of influence.

We are going to apply the estimate (\ref{Blplq}) to the problem     
\begin{eqnarray}
\label{2.1}
&   &
\partial_t^2 v - {\mathcal A}_{\varepsilon }(x,\partial_x)v =0, \quad x \in {\mathbb R}^3, \quad t \geq 0,  \\
\label{2.2}
&  &
v(x,0) = v_0 (x),\quad   v_t (x,0)= v_1(x),\quad x \in {\mathbb R}^3\,,
\end{eqnarray}
where the operator ${\mathcal A}_{\varepsilon }(x,\partial_x)$ has a symbol ${\mathcal A}_{\varepsilon }(x;\xi)$ of   (\ref{AUX}). In that case the constant $C_s(T) $ depends on $ \varepsilon $ as well.

The conservation of the energy of the solution of the equation 
\begin{eqnarray}
\label{Lin5.1}
 &  &
    \frac{\partial^2 \psi }{\partial t^2}
-{\mathcal{A}}(x,\partial_x)\psi +  F(r) \frac{m^2 c^4 }{h^2}  \psi  =0 
\end{eqnarray}
is known (see, e.g., \cite{Nicolas}). More exactly, 
for   initial data with the supports in $B^{ext}_{Sch} $ the energy
\begin{eqnarray*}
E(t) & = &
 \int_{R_{Sch}}^\infty \int_{{\mathbb S}^2} \left\{ \frac{1}{F(r)} |  \partial_t  \psi |^2 
+ 
    F(r) |\partial_r \psi  |^2  +  \frac{1}{r^2}|\nabla _{{\mathbb S}^2}\psi|^2 + \frac{m^2 c^4 }{h^2}  |\psi|^2\right\}r^2\, dr \, d \Omega_2
\end{eqnarray*}
is conserved as long as the solution exists, that is, for all time of the existence of the solution,   
\begin{eqnarray} 
\label{T_EC}
&  &
\frac{d}{dt}E(t) =0\,.
\end{eqnarray}
We write the energy 
in  Cartesian coordinates as follows
 \begin{eqnarray*}
E(t) 
& = &
 \int_{\R^3} \Bigg\{ \left(1-\frac{2 M_{bh}}{\sqrt{x^2+y^2+z^2}}\right)^{-1}  |   \psi_t (x,y,z)|^2
+  |\psi_x  (x,y,z)|^2  +  |\psi_y(x,y,z)|^2  + |\psi_z (x,y,z)|^2 \nonumber \\
&  &
-2 M_{bh}  \frac{1}{\left(x^2+y^2+z^2\right)^{3/2}}\Big| x  \psi_x (x,y,z)  + y      \psi_y (x,y,z)    + z   \psi_z (x,y,z) \Big|^2  
+ \frac{m^2 c^4 }{h^2}  |\psi|^2 \Bigg\} \,dx\,dy\,dz\,. 
\end{eqnarray*}

\subsection{Equation in self-adjoint Cauchy-Kowalewski form}

The semi-linear Klein-Gordon equation without potential is  
\begin{eqnarray}
\label{KGSLKF}
&  &
    \frac{\partial^2  \psi  }{\partial t^2}
+ 3H     \frac{\partial \psi }{\partial t} 
-e^{ -2Ht } {\mathcal{A}}(x,\partial_x)\psi +F (r) \frac{m^2 c^4 }{h^2}\psi = c^2F (r)\Psi(\psi)\,, 
\end{eqnarray}
where $ {\mathcal{A}}(x,\partial_x)$ is defined in (\ref{OpA}). 
The   operator $ {\mathcal{A}}(x,\partial_x)={\mathcal A}(\vec{ x},\partial_x)={\mathcal{A}}(x,y,z;D_x,D_y,D_z)$ 
has the symbol (\ref{SymbA}) 
in Cartesian coordinates.  
We apply the  Liouville transform
\[
\psi= e^{-\frac{3H}{2 }t}\sqrt{F(r)}u
\]
to the     covariant   Klein-Gordon  equation (\ref{KGSLKF}), which turns into the  
non-covariant equation
\[ 
   \frac{1}{c^2}\frac{\partial^2   u }{\partial t^2} 
- e^{-2Ht}{\mathcal A}_{3/2}(x,\partial_x)u 
+\left(   \frac{m^2 c^2 }{h^2}-\frac{9H^2}{4c^2}\right)u   -\frac{2G M_{bh}}{c^2r}   \frac{m^2 c^2 }{h^2}  u   
  =   e^{\frac{3H}{2 }t}\sqrt{F(r)}\Psi \left(e^{-\frac{3H}{2 }t}\sqrt{F(r)}u \right)\,, 
\]
where  the  operator ${\mathcal A}_{3/2}$ in the spherical coordinates is defined by
\[
{\mathcal A}_{3/2}(x,\partial_x)v
  :=   
    F (r)^{3/2} \frac{\partial^2    }{\partial r^2}\sqrt{F(r)}v
+\sqrt{F(r)}\frac{2}{  r   }    \left( 1-\frac{ G M_{bh}}{c^2r} \right)  \frac{\partial   }{\partial r}   \sqrt{F(r)}v
+ F(r) \frac{1}{  r^2 } \Delta_{S^2} v \,,
\]
while the term
$
-\frac{2G M_{bh}}{c^2r}   \frac{m^2 c^2 }{h^2}
$
 can be regarded as a potential $V=V(x,y,z)$.

If $H=0$, then we obtain from (\ref{KGSLKF}) the semi-linear Klein-Gordon equation in the static universe
\begin{eqnarray}
\label{KGSLKFH0}
&  &
    \frac{\partial^2  \psi  }{\partial t^2} 
- {\mathcal{A}}(x,\partial_x)\psi +F (r) \frac{m^2 c^4 }{h^2}\psi = c^2F (r)\Psi(\psi)\,.
\end{eqnarray}
The lemma below shows that the  Liouville transform   makes self-adjoint the spatial part of the operator of the left-hand side of (\ref{KGSLKFH0}),  
and that equation   reads
\[
 \frac{\partial^2   v }{\partial t^2} 
- c^2{\mathcal A}_{3/2}(x,\partial_x)v
+F (r) \frac{m^2 c^4 }{h^2}v = c^2\sqrt{F(r)}\Psi(\sqrt{F(r)}v)\,.
\]
The symbol ${\mathcal A}_{3/2}(x,\xi ) $  of operator ${\mathcal A}_{3/2} (x,\partial_x) $ in Cartesian coordinates   is
\[ 
{\mathcal A}_{3/2}(\vec{x},\vec{\xi}) 
  :=   
\left(1-\frac{2 G M_{bh}}{c^2 |\vec{x}|}\right) \left(-|\vec{\xi}|^2 +\frac{2 G M_{bh} \left(\vec{x}\cdot\vec{\xi}\right)^2}{c^2 |\vec{x}|^{3}}\right)
+\frac{G^2 M_{bh}^2}{c^4 |\vec{x}|^4} .
\]
\begin{lemma}
The operator ${\mathcal A}_{3/2}(x,\partial_x)$  is self-adjoint   on $C_0^\infty(B^{ext}_{Sch})$. On every closed subset in $B^{ext}_{Sch}$ the operator ${\mathcal A}_{3/2}(x,\partial_x)$ is an elliptic  operator that is non-positive on the subspace of functions with the supports in $B^{ext}_{Sch}$.
\end{lemma}

\noindent
{\bf Proof.}  
For the   vectors $\vec{x}=(x_1,x_2,x_3):=(x,y,z) \in \R^3 $ and $\vec{\xi}=(\xi_1,\xi_2,\xi_3) \in \R^3 $,  the  Cauchy inequality implies 
\[
 |\vec{\xi}|^2   
  >   
 \frac{2 G M_{bh} \left(\vec{x}\cdot \vec{\xi}  \right)^2}{c^2 | \vec{x}|^3 }\,.
\]
The direct calculations of the symbol of the  adjoint operator 
show  that ${\mathcal A}_{3/2}$ is self-adjoint. Indeed, 
\[
{\mathcal A}_{3/2}(\vec{x},\vec{\xi})=\sum_{|\alpha|=0,1,2}\frac{(-i)^{|\alpha|}}{\alpha !}\partial_x^\alpha \partial_\xi^\alpha {\ol {\mathcal A}_{3/2} (\vec{x},\vec{\xi})}.
\]
Moreover,
\[
\sum_{k=1,2,3} -i \frac{\partial ^2{\mathcal A}_{3/2}(\vec{x},\vec{\xi})}{\partial x_k\, \partial \xi _k}=0\,.
\]
The operator ${\mathcal A}_{3/2}(x,y,z;D_x,D_y,D_z)$  can be written as follows
\[
{\mathcal A}_{3/2}(x,y,z;D_x,D_y,D_z)v=\sum_{i,j=1}^3\frac{\partial}{\partial x_i}\left( a_{ij}(\vec{x})\frac{\partial}{\partial x_j}v \right)  +\frac{G^2 M^2}{c^4 |\vec{x}|^4 }v\,,
\]
since 
$\sum_{i=1}^3 \frac{\partial }{\partial x_i}a_{ij}(\vec{x}) =0$  for $j=1,2,3$,  
with the coefficients $ a_{ij}(\vec{x}) $   such that
\begin{eqnarray*}
a_{kk}(\vec{x}) 
& := &
-\left(1-\frac{2 G M_{bh}}{c^2 |\vec{x}|}\right) \left(\frac{2 G M_{bh} x_k^2}{c^2 |\vec{x}|^{3 }}-2\right), \quad k=1,2,3,\\
a_{k\ell}(\vec{x}) 
& = & 
a_{\ell k}(\vec{x})= 
-\left(1-\frac{2 G M_{bh}}{c^2 |\vec{x}|}\right)\frac{ 2  G M_{bh} x_k x_\ell  }{c^2 |\vec{x}|^{3 }}, \quad k,\ell=1,2,3,\,\, k\not=\ell \,.
\end{eqnarray*}
Thus, one can write
\begin{eqnarray*}
&  &
\int_{\R^3}({\mathcal A}_{3/2}(x,\partial_x)u(x,y,z))\ol{v(x,y,z)}\,dx dy dz\\
& = &
\int_{\R^3} \left\{ \sum_{i,j=1}^3a_{ij}(x,y,z)\frac{\partial^2  u(x,y,z)}{\partial x_i\partial x_j} \ol{v(x,y,z)}   
+\frac{G^2 M_{bh}^2}{c^4 \left(x^2+y^2+z^2\right)^2}u(x,y,z) \ol{v(x,y,z)}\right\}\,dx dy dz \\
& = &
\int_{\R^3} \left\{ \sum_{i,j=1}^3a_{ij}(x,y,z)u (x,y,z)\frac{\partial^2 \ol{v(x,y,z)}}{\partial x_i\partial x_j}    +\frac{G^2 M_{bh}^2}{c^4 \left(x^2+y^2+z^2\right)^2}u(x,y,z) \ol{v(x,y,z)}\right\}\,dx dy dz  
\end{eqnarray*}
 for every $u,v \in C_0^\infty(B_{Sch}^{ext})$ functions.

The bilinear form $ \sum_{i,j=1}^3  a_{ij}(x,y,z) \xi_{i }\xi_{ j} $ is positive in    $B^{ext}_{Sch}$   since the principal minors are 
\[
M_1=\left(1-\frac{2 G M_{bh}}{c^2 |\vec{x}|}\right) \left( 2-\frac{2 G M_{bh} x^2}{c^2 |\vec{x}|^{3}}\right), \quad
 M_2=4-\frac{4 G M_{bh} \left(x^2+y^2\right)}{c^2 |\vec{x}|},\quad M_3=8-\frac{8 G M_{bh}}{c^2 |\vec{x}|}\,.
\]
Next, consider for the real valued function $ v \in C_0^\infty(B_{Sch}^{ext})$  the inner product
\begin{eqnarray*}
&  &
 ({\mathcal A}_{3/2}(x,\partial_x)v,v)_{L^2(\R^3)}\\
&  = & 
   \int_0^\infty \int_{S^2} F (r) \left(\frac{\partial^2    }{\partial r^2}\sqrt{F(r)}v \right) \sqrt{F(r)} v\,r^2dr d\Omega_2\\
   &  &
+ \int_0^\infty \int_{S^2}   \frac{2}{  r   }    \left( 1-\frac{ G M_{bh}}{c^2r} \right)  \left(\frac{\partial   }{\partial r}   \sqrt{F(r)}v\right)\sqrt{F(r)}v\,r^2dr d\Omega_2
+  \int_0^\infty \int_{S^2}  \left( F(r) \frac{1}{  r^2 } \Delta_{S^2} v\right)v\,r^2dr d\Omega_2\,.
\end{eqnarray*}
Then we integrate by parts  the first term of the last identity
\begin{eqnarray*}
&  &
   \int_0^\infty \int_{S^2} F (r) \left(\frac{\partial^2    }{\partial r^2}\sqrt{F(r)}v \right) \sqrt{F(r)} v\,r^2dr d\Omega_2\\
   & = &
  - \int_0^\infty \int_{S^2} F (r) \left(\frac{\partial     }{\partial r }\sqrt{F(r)}v \right)^2  \,r^2dr d\Omega_2
- \int_0^\infty \int_{S^2} \left(\frac{ 2G M_{bh}}{c^2r^2} \right)\left(\frac{\partial     }{\partial r }\sqrt{F(r)}v \right)\sqrt{F(r)}v   \,r^2dr d\Omega_2\\
&  &
- \int_0^\infty \int_{S^2}   F (r)   \left(\frac{\partial     }{\partial r }\sqrt{F(r)}v \right)\sqrt{F(r)}v2r   \, dr d\Omega_2\,.
\end{eqnarray*}
Hence,
\begin{eqnarray*}
&  &
 ({\mathcal A}_{3/2}(x,\partial_x)v,v)_{L^2(\R^3)}\\
&  = & 
  - \int_0^\infty \int_{S^2} F (r) \left(\frac{\partial     }{\partial r }\sqrt{F(r)}v \right)^2  \,r^2dr d\Omega_2
- \int_0^\infty \int_{S^2} \left(\frac{ 2G M_{bh}}{c^2r^2} \right)\left(\frac{\partial     }{\partial r }\sqrt{F(r)}v \right)\sqrt{F(r)}v   \,r^2dr d\Omega_2\\
&  &
- \int_0^\infty \int_{S^2}   F (r)   \left(\frac{\partial     }{\partial r }\sqrt{F(r)}v \right)\sqrt{F(r)}v2r   \, dr d\Omega_2\\
   &  &
+ \int_0^\infty \int_{S^2}   \frac{2}{  r   }    \left( 1-\frac{ G M_{bh}}{c^2r} \right)  \left(\frac{\partial   }{\partial r}   \sqrt{F(r)}v\right)\sqrt{F(r)}v\,r^2dr d\Omega_2
+  \int_0^\infty \int_{S^2}  \left( F(r) \frac{1}{  r^2 } \Delta_{S^2} v\right)v\,r^2dr d\Omega_2\,.
\end{eqnarray*}
The first and last terms of the previous identity are non-positive. The sum of the remaining terms  vanishes 
\begin{eqnarray*}
&  & 
- \int_0^\infty \int_{S^2} \left(\frac{ 2G M_{bh}}{c^2r^2} \right)\left(\frac{\partial     }{\partial r }\sqrt{F(r)}v \right)\sqrt{F(r)}v   \,r^2dr d\Omega_2\\
&  &
- \int_0^\infty \int_{S^2}   F (r)   \left(\frac{\partial     }{\partial r }\sqrt{F(r)}v \right)\sqrt{F(r)}v2r   \, dr d\Omega_2\\
   &  &
+ \int_0^\infty \int_{S^2}   \frac{2}{  r   }    \left( 1-\frac{ G M_{bh}}{c^2r} \right)  \left(\frac{\partial   }{\partial r}   \sqrt{F(r)}v\right)\sqrt{F(r)}v\,r^2dr d\Omega_2 \\
& = &
\int_0^\infty \int_{S^2} \Bigg[ -   \frac{ 2G M_{bh}}{c^2 }   -   2rF (r)
+ 2r   \left( 1-\frac{ G M_{bh}}{c^2r} \right)\Bigg]\left(\frac{\partial     }{\partial r }\sqrt{F(r)}v \right)\sqrt{F(r)}v   \,dr d\Omega_2\\
& = &
0
\end{eqnarray*}
for the real valued smooth function with the compact support in $B^{ext}_{Sch}$. 
The lemma is proved. \qed

\section{Representation formula for the solution of generalized linear Klein-Gordon equation in de~Sitter space-time}

\setcounter{equation}{0}
\renewcommand{\theequation}{\thesection.\arabic{equation}}

In this section we set $c=1$. 
We recall (see \cite{Yag_Galst_CMP,MN2015})  the fundamental solutions  of the Klein-Gordon equation in the de~Sitter space-time. For $(x_0, t_0) \in {\mathbb R}^n\times {\mathbb R}$, $M \in {\mathbb C}$,  we define 
{\it chronological future} (``forward light cone'') $D_+(x_0, t_0)$  of the point $(x_0,t_0)  \in {\mathbb R}^4$ and the {\it chronological past} (``backward light cone'') $D_-(x_0, t_0)$.  The forward and backward light cones are defined as follows:
\[ 
D_{\pm}\left(x_{0}, t_{0}\right) :=\left\{(x, t) \in {\mathbb R}^{3+1} ;
\left|x-x_{0}\right| \leq \pm\left(\phi (t) -\phi (t_{0})  \right)\right\}\,,
\]
where $\phi (t):= (1-e^{-Ht} )/H$ is a distance function.  
In fact, any intersection of $D_-(x_0, t_0)$ with the hyperplane $t = const < t_0$ determines the
so-called {\it dependence domain} for the point $ (x_0, t_0)$, while the intersection of $D_+(x_0, t_0)$
with the hyperplane $t = const > t_0$ is the so-called {\it domain of influence} of the point
$ (x_0, t_0)$.
We define also the function
\begin{eqnarray*}
E(x,t;x_0,t_0;M)
& :=  &
4^{-\frac{M}{H}} e^{ M  (  t_0+ t)} \left(\left(e^{-H t_0}+e^{-H t}\right)^2-(x - x_0)^2\right)^{\frac{M}{H}-\frac{1}{2}}  \\
&  &
\times  F \left(\frac{1}{2}-\frac{M}{H},\frac{1}{2}-\frac{M}{H};1;\frac{\left(e^{-H t}-e^{-H t_0}\right)^2-(x - x_0)^2}{\left(e^{-H t}+e^{-Ht_0}\right)^2-(x - x_0)^2}\right)\,, \nonumber
\end{eqnarray*}
where   $(x,t) \in D_+ (x_0,t_0)\cup D_- (x_0,t_0) $  and $F\big(a, b;c; \zeta \big) $ is the hypergeometric function (see, e.g.,\cite{B-E}). 
When no ambiguity arises,  we use the notation $x^2:= |x|^2$ for $x \in {\mathbb R}^n $.
Thus, the function $E  $ depends on $r^2= (x- x_0 )^2/H^2$, and we will write $E  (r,t;0,t_0;M) $ for   $E(x,t;x_0,t_0;M)  $:
\begin{eqnarray}
\label{0.6} 
E(r,t;0,t_0;M)
& :=  &
4^{-\frac{M}{H}} e^{ M  (  t_0+ t)} \left(\left(e^{-H t_0}+e^{-H t}\right)^2-(H r)^2\right)^{\frac{M}{H}-\frac{1}{2}}  \nonumber \\
&  &
\times  F \left(\frac{1}{2}-\frac{M}{H},\frac{1}{2}-\frac{M}{H};1;\frac{\left(-e^{-H t}+e^{-H t_0}\right)^2-(r H)^2}{\left(e^{-H t}+e^{-Ht_0}\right)^2-(r H)^2}\right)  \,.
\end{eqnarray}
Additional to (\ref{0.6}) 
for $M \in {\mathbb C}$  we recall two more kernel functions from \cite{Yag_Galst_CMP,MN2015}
\begin{eqnarray}
\label{K0MH}
K_0(r,t;M)
&  := &
-\left[\frac{\partial}{\partial b} E(r, t ; 0, b ; M)\right]_{b=0}\,,   \\
\label{K1MH}
K_1(r,t;M)
& :=  &
E(r,t;0, 0;M)   \,.  
\end{eqnarray}
Then according to \cite{MN2015} the solution operator for the Cauchy problem for the scalar {\it generalized Klein-Gordon equation}
in the de~Sitter space-time  
\[
\left( \partial_t^2     
-  e^{-2Ht}    {\mathcal A}(x,\partial_x)     
 -  M ^2   \right)u  =f,\quad u (x,0)= u_0 (x), \quad u_t (x,0)= u _1 (x) \,,  
 \]
is given as follows
\[
u(x,t) = {\cal G} ( x,t,D_x;M)[f]+{\cal K}_0(x,t,D_x;M)[u_0 ]+ {\cal K}_1(x,t,D_x;M)[u_1 ] \,.
\]
Here ${\mathcal A}(x,\partial_x) $ is the   differential operator 
\[
{\mathcal A}(x,\partial_x)=\sum_{|\alpha| \leq m} a_\alpha (x)D_x^\alpha \,, \qquad a_\alpha \in C^\infty (\Omega )\,, 
\]
 and the coefficients $ a_\alpha (x)$ 
are $C^\infty$-functions in the open domain $\Omega \subseteq {\mathbb R}^n $. 
The kernels   $K_0(z,t;M) $   and  $K_1(z,t;M) $ can be written in the explicit form as follows
\begin{eqnarray*}
K_0(r,t;M) 
\!&\!=\!&\!
-4^{-\frac{M}{H}} \left((1+e^{  -H t})^2- H^2 r^2\right)^{\frac{M}{H}-\frac{5}{2}} e^{t (M-4 H)}
 \\
&  &
\times \Bigg\{ \! e^{  -2H t}\left( (1+e^{  -H t})^2- H^2 r^2  \right) \!\!\left(-e^{2 H t} \left(H  (H M r^2-1 )+M\right)+H e^{H t}+M\right)  \\
&  &
\times F \left(\frac{1}{2}-\frac{M}{H},\frac{1}{2}-\frac{M}{H};1;\frac{\left(1-e^{-H t}\right)^2-H^2 r^2}{\left(1+e^{-H t}\right)^2-H^2 r^2}\right)\\
&  &
+\frac{1}{H} (H-2 M)^2 e^{3t H}   \left( e^{-2 H t}- H^2 r^2-1 ) \right)  
 \\
&  &
\left. \times F \left(\frac{3}{2}-\frac{M}{H},\frac{3}{2}-\frac{M}{H};2;\frac{\left(1-e^{-H t}\right)^2-H^2 r^2}{\left(1+e^{-H t}\right)^2-H^2 r^2}\right)\right\},\\ 
K_1(r,t;M)
& =  &
4^{-\frac{M}{H}} e^{M t} \left(\left(1+e^{-H t}\right)^2-(H r)^2\right)^{\frac{M}{H}-\frac{1}{2}} \\
&  &
\times  F \left(\frac{1}{2}-\frac{M}{H},\frac{1}{2}-\frac{M}{H};1;\frac{\left(1-e^{-H t} \right)^2-(r H)^2}{\left(1+e^{-H t} \right)^2-(r H)^2}\right) \nonumber \,.  
\end{eqnarray*} 
Next we 
  recall the results of  Theorem~1.1~\cite{MN2015}.  
For $ f \in C^\infty (\Omega\times I  )$,\, $ I=[0,T]$, $0< T \leq \infty$, and \, $ \varphi_0 $,  $ \varphi_1 \in C (\Omega ) $,
let  the function\,
$v_f(x,t;b) \in C_{x,t,b}^{m,2,0}(\Omega \times [0,(1-e^{-HT})/H]\times I)$\,
be a solution to the   problem (\ref{vf})  
and the function \, $  v_\varphi(x, t) \in C_{x,t}^{m,2}(\Omega \times [0,(1-e^{-HT})/H])$ \, be a   solution   of the   problem
(\ref{vphi}). 
Then the function  $u= u(x,t)$   defined by
\begin{eqnarray}
\label{u}
u(x,t)
&  =  &
2   \int_{ 0}^{t} db
  \int_{ 0}^{\phi (t)- \phi (b)}  E(r,t;0,b;M)  v_f(x,r ;b) \, dr  
+ e ^{\frac{Ht}{2}} v_{u_0}  (x, \phi (t))   \\
&  &
+ \, 2\int_{ 0}^{\phi (t)}  K_0( s,t;M)v_{u_0}  (x, s)   ds  \nonumber 
+\, 2\int_{0}^{\phi (t) }  v_{u_1 } (x,  s)
  K_1( s,t;M)   ds
, \quad x \in \Omega  , \,  t \in I ,
\end{eqnarray}
where $\phi (t):= (1-e^{-Ht} )/H$,  
 solves the problem
\[
\cases{
u_{tt} - e^{-2Ht}{\mathcal A}(x,\partial_x)  u - M^2 u= f, \quad  x \in \Omega \,,\,\, t \in I,\cr
  u(x,0)= u_0 (x)\, , \quad u_t(x,0)=u_1 (x),\quad x \in \Omega\,.
}
\]
Here the kernels  $E$, $K_0$, and $K_1$ have been defined in (\ref{0.6}), (\ref{K0MH}), and (\ref{K1MH}), respectively.
Consequently, for $n=3$ we have the representation (\ref{psi}).

We need also the second  equivalent  form  of the kernel $K_0 $ given in the next statement.  
\begin{lemma} 
 \label{L8.1}
 The kernel $K_0 $ can be written as follows
\begin{eqnarray*}
&  &
K_0(r,t;M) \\
\!&\!=\!&\!
\frac{4^{-\frac{M}{H}} e^{M t} \left(\left(e^{-H t}+1\right)^2-H^2 r^2\right)^{\frac{M}{H}-\frac{1}{2}}}{\left(1-e^{-H t}\right)^2-H^2 r^2}\\
&  &
\times\Bigg[ \left(H e^{-H t}-H+M e^{-2 H t}-M-H^2 M r^2\right)  F \left(\frac{1}{2}-\frac{M}{H},\frac{1}{2}-\frac{M}{H};1;\frac{\left(-1+e^{-H t}\right)^2-H^2 r^2}{\left(1+e^{-H t}\right)^2-H^2 r^2}\right)\\
&  &
+\left(\frac{H}{2}+M\right) \left(H^2 r^2-e^{-2 H t}+1\right)  F \left(-\frac{1}{2}-\frac{M}{H},\frac{1}{2}-\frac{M}{H};1;\frac{\left(-1+e^{-H t}\right)^2-H^2 r^2}{\left(1+e^{-H t}\right)^2-H^2 r^2}\right)\Bigg]\,.  
\end{eqnarray*}
\end{lemma} 
\medskip

\ndt
{\bf Proof.}  Indeed, according to  \cite[(42), Sec.2.8]{B-E} 
\[
F\left(\frac{3}{2}-\frac{M}{H},\frac{3}{2}-\frac{M}{H};2;z\right)\\
  =  
\frac{1}{z \left(\frac{M}{H}-\frac{1}{2}\right) }F\left(\frac{3}{2}-\frac{M}{H},\frac{1}{2}-\frac{M}{H};1;z\right)-(1-z) F\left(\frac{3}{2}-\frac{M}{H},\frac{3}{2}-\frac{M}{H};1;z\right)\,.
 \]
 Next we apply \cite[(36), Sec.2.8]{B-E} and write
\begin{eqnarray*} 
 F \left(\frac{1}{2}-\frac{M}{H},\frac{3}{2}-\frac{M}{H};1;z\right)
& = &
-\frac{1}{(1-z) \left(\frac{1}{2}-\frac{M}{H}\right)}\Bigg[ \frac{2 M }{H} F \left(\frac{1}{2}-\frac{M}{H},\frac{1}{2}-\frac{M}{H};1;z\right)\\
&  &
-\left(\frac{M}{H}+\frac{1}{2}\right)    F \left(-\frac{M}{H}-\frac{1}{2},\frac{1}{2}-\frac{M}{H};1;z\right)\Bigg]\,.
\end{eqnarray*}
Then we    apply \cite[(36), Sec.2.8]{B-E} once again
\begin{eqnarray*} 
 F \left(\frac{3}{2}-\frac{M}{H},\frac{3}{2}-\frac{M}{H};1;z\right) 
& = &
-\frac{1}{(z-1)^2 (H-2 M)}\Bigg[(H-H z+ 2 M (z+3))  F \left(\frac{1}{2}-\frac{M}{H},\frac{1}{2}-\frac{M}{H};1;z\right)\\
&  &
-2 (H+2 M)  F \left(-\frac{1}{2}-\frac{M}{H},\frac{1}{2}-\frac{M}{H};1;z\right)\Bigg] \,. 
\end{eqnarray*}
Hence
\begin{eqnarray*}
F\left(\frac{3}{2}-\frac{M}{H},\frac{3}{2}-\frac{M}{H};2;z\right)
& = &
-\frac{ 2 H }{(z-1) z (H-2 M)^2}\Bigg[  
(H+2 M)  F \left(-\frac{1}{2}-\frac{M}{H},\frac{1}{2}-\frac{M}{H};1;z\right)\\
&  &
+(H (z-1)-2 M (z+1))  F \left(\frac{1}{2}-\frac{M}{H},\frac{1}{2}-\frac{M}{H};1;z\right)\Bigg]\,.
\end{eqnarray*}
For $z=  \left((e^{-H t}-1 )^2-H^2 r^2\right)/\left((e^{-H t}+1)^2-H^2 r^2\right)$, it follows the statement of the lemma.   \qed

\section{The semilinear equation with  large mass. Proof of Theorem~\ref{T_large_mass}}

\setcounter{equation}{0}
\renewcommand{\theequation}{\thesection.\arabic{equation}}
Let $M$ be a non-negative number such that 
$M^2:= \frac{m^2c^4}{h^2} - \frac{9H^2}{4 }\geq 0 $.

\subsection{The linear equation  without potential and source terms}

In this subsection we obtain 
decay estimates of solution of the linear equation 
\begin{eqnarray}
\label{setCP_WPST}
 &  &
    \frac{\partial^2 \psi }{\partial t^2}
+ 3H        \frac{\partial \psi }{\partial t}
-e^{-2Ht}{\mathcal A}(x,\partial_x) \psi +   \frac{m^2c^4}{h^2} \psi =0\,, \quad x \in \R^3,\,\, t \in [0,\infty)\,.
\end{eqnarray}

\begin{theorem}
\label{T13.2b}
For the solution of the Cauchy problem (\ref{setCP_WPST})\&(\ref{13.12})\&(\ref{13.13}), the following estimate holds:
\begin{eqnarray}
\label{LpLq}
\|\psi (x,t)\| _{H_{(s)}}
&  \leq &
 C_{H,\chi}  (1+t)^{1-\sgn (M)} e^{ -Ht}   \|\psi_0  (x)  \| _{H_{(s)}}\\
&  & 
+\,C_{H,\chi} (1+t)^{1-\sgn (M)} e^{-\frac{3H}{2}t}  (e^{H t}  - 1) (e^{H t}  + 1)^{-1 }\|  \psi_1 (x)\| _{H_{(s)}} \quad \mbox{  for all} \quad t>0.\nonumber
\end{eqnarray} 
\end{theorem}
\noindent
{\bf Proof.} 
Fix a cutoff function $\chi \in \C^\infty_0 (\R^3) $ such that  
$
 \chi (x)=1$  for all $(x,t)  \in \mbox{\rm supp}\, u $,  $t \in[0,\infty) .
$
For the function  $u= u(x,t)$   defined by (\ref{u}), when $ f=0$ and  $t \in[0,\infty)$, according to (\ref{Blplq}),   we obtain
\begin{eqnarray*}
\|u(x,t)\| _{H_{(s)}}
&  =  &
\|e ^{\frac{Ht}{2}}\chi (x)  v_{u_0}  (x, \phi (t))
+ \, 2\int_{ 0}^{\phi (t)}  K_0( s,t;-iM)\chi (x)v_{u_0}  (x, s)   ds  \nonumber \\
&  &
+\, 2\int_{0}^{\phi (t) } \chi (x)  v_{u_1 } (x,  s)
  K_1( s,t;-iM)   ds\| _{H_{(s)}} \,.
\end{eqnarray*}
We consider the case of   $s >[3/2]+1$, then
\begin{eqnarray*} 
\|u(x,t)\| _{H_{(s)}} 
&  \leq &
C_\chi e ^{\frac{Ht}{2}} \|v_{u_0}  (x, \phi (t))\| _{H_{(s)}}
+ \, C_\chi \int_{ 0}^{\phi (t)} | K_0( s,t;-iM) |\|v_{u_0}  (x, s)   ds \| _{H_{(s)}} \\
&  & 
+\,C_\chi\int_{0}^{\phi (t) }  | K_1( s,t;-iM) | \| v_{u_1 } (x,  s)\| _{H_{(s)}} 
 ds\\ 
&  \leq &
C_{H,\chi} e ^{\frac{Ht}{2}} \| u_0   (x)\| _{H_{(s)}}
+ \, C_{H,\chi} \|u_0  (x)  \| _{H_{(s)}}\int_{ 0}^{\phi (t)} | K_0( s,t;-iM) |   ds\\
&  & 
+\,C_{H,\chi}\| u_1 (x)\| _{H_{(s)}} \int_{0}^{\phi (t) }  | K_1( s,t;-iM) | 
 ds\,.
\end{eqnarray*}
Next we apply Lemma~\ref{L13.1} and Lemma~\ref{L13.2} (see below) and obtain 
\begin{eqnarray*}
\|u(x,t)\| _{H_{(s)}}
&  \leq &
C_{H,\chi} e ^{\frac{Ht}{2}} \| u_0   (x)\| _{H_{(s)}}
+ \, C_{H,\chi} \|u_0  (x)  \| _{H_{(s)}} (1+t)^{1-\sgn (M)} e^{- \frac{Ht}{2}} (e^{Ht}-1) \\
&  & 
+\,C_{H,\varphi}\| u_1 (x)\| _{H_{(s)}}    (1+t)^{1-\sgn (M)} (e^{H t}  - 1) (e^{H t}  + 1)^{-1 }\,.
\end{eqnarray*}

Thus, for the solution $\psi $ of the Cauchy problem (\ref{setCP_WPST}),  
due to the relations 
$
u = e^{\frac{3H}{2}t}\psi  $, $\,u_0=\psi_0$, and $ u_1=\frac{3H}{2}\psi_0+\psi_1$,  
we obtain the  estimate 
\begin{eqnarray*}
\|\psi (x,t)\| _{H_{(s)}}
&  \leq &
e^{-\frac{3H}{2}t}\Bigg[ C_{H,\chi} e ^{\frac{Ht}{2}} \| \psi_0   (x)\| _{H_{(s)}}
+ \, C_{H,\chi} \|\psi_0  (x)  \| _{H_{(s)}} (1+t)^{1-\sgn (M)} e^{- \frac{Ht}{2}} (e^{Ht}-1) \\
&  & 
+\,C_{H,\chi}\| \frac{3H}{2}\psi_0 (x)+\psi_1 (x)\| _{H_{(s)}}   (1+t)^{1-\sgn (M)} (e^{H t}  - 1) (e^{H t}  + 1)^{-1 }\Bigg]\\
&  \leq &
 C_{H,\chi}  (1+t)^{1-\sgn (M)} e^{ -Ht}   \|\psi_0  (x)  \| _{H_{(s)}}\\
&  & 
+\,C_{H,\chi} (1+t)^{1-\sgn (M)} e^{-\frac{3H}{2}t}  (e^{H t}  - 1) (e^{H t}  + 1)^{-1 }\|  \psi_1 (x)\| _{H_{(s)}}  
\end{eqnarray*}
for large $t$. The theorem is proved. 
\qed

\begin{lemma} 
\label{L13.1} Let $M\geq 0$ and $H $ be a positive number, then
\[
\int_{ 0}^{(1-e^{-Ht} )/H} | K_0( s,t;-iM) |   ds 
  \leq  
C_{M,H}(1+t)^{1-\sgn (M)} e^{- Ht/2} (e^{Ht}-1) \quad \mbox{for all}\quad t \in [0,\infty).
\]
\end{lemma}
\medskip

\noindent
{\bf Proof.} We use the $ K_0( r,t;-iM) $ given by   Lemma~\ref{L8.1}. 
Then, since $M$ is real, we obtain
\begin{eqnarray*}
&  &
\int_{ 0}^{(1-e^{-Ht} )/H} | K_0( r,t;-iM) |   dr\\
& \leq &
\int_{ 0}^{(1-e^{-Ht} )/H} \Bigg|\frac{ \left(\left(e^{-H t}+1\right)^2-H^2 r^2\right)^{-\frac{1}{2}}}{\left(1-e^{-H t}\right)^2-H^2 r^2}\Bigg[\left(H^2 r^2  iM- iM e^{-2 H t}+H e^{-H t}-H+ iM \right)  \\
&  &
\times F \left(\frac{1}{2}+\frac{ iM}{H},\frac{1}{2}+\frac{ iM}{H};1;\frac{\left(-1+e^{-H t}\right)^2-H^2 r^2}{\left(1+e^{-H t}\right)^2-H^2 r^2}\right)  \\
&  &
+ \left(\frac{H}{2}- iM\right) \left(H^2 r^2-e^{-2 H t}+1\right) \\
&  &
\times F\left(-\frac{1}{2}+\frac{iM}{H},\frac{1}{2}+\frac{iM}{H} ;1;\frac{\left(-1+e^{-H t}\right)^2-H^2 r^2}{\left(1+e^{-H t}\right)^2-H^2 r^2}\right)
\Bigg] \Bigg|   dr\,.
\end{eqnarray*}
Hence,  with    $z:=e^{Ht} $ and  
$r=y/(Hz)$, we derive
\begin{eqnarray*} 
&  &
\int_{ 0}^{(1-e^{-Ht} )/H} | K_0( r,t;-iM) |   dr \\
\!&\! \leq \!&\!
\int_{ 0}^{z-1 } \Bigg|\frac{ \left(\left(z+1\right)^2- y^2\right)^{-\frac{1}{2}}}{\left(z-1 \right)^2- y^2}\\
&  &
\times\Bigg[\left(y^2  i\frac{ M}{H }- i\frac{ M}{H }  + z-z^2+ z^2i\frac{ M}{H } \right)  F \left(\frac{1}{2}+\frac{ iM}{H},\frac{1}{2}+\frac{ iM}{H};1;\frac{\left(z-1 \right)^2- y ^2  }{\left(z+1 \right)^2-  y ^2}\right)  \\
&  &
+ \left(\frac{1}{2}- i\frac{ M}{H }\right) \left(y^2-1+z^2\right) F\left(-\frac{1}{2}+\frac{iM}{H},\frac{1}{2}+\frac{iM}{H} ;1;\frac{\left(z-1 \right)^2- y ^2  }{\left(z+1 \right)^2-  y ^2}\right)
\Bigg] \Bigg|   d y\\
\!&\! \leq \!&\!
C_{M/H}(z+1)^{-1/2}(z-1) \quad \mbox{\rm for all}\quad z \in [1,\infty)\,.
\end{eqnarray*}
In the last step, we used the  estimates, which are    proved  in \cite[Lemma 7.4]{Yag_Galst_CMP}.
The lemma is proved. \qed

Next we consider the kernel $K_1$.
\begin{lemma}
\label{L13.2} Let $M\geq 0$ and $H $ be a positive number,  then
\[
 \int_{0}^{\phi (t) }  | K_1( r,t;-iM) |  dr 
  \leq   
   C_{M } (1+t)^{1-\sgn (M)}\frac{1}{ H} (e^{H t}  - 1) (e^{H t}  + 1)^{-1 } \quad \mbox{for all} \quad t \in [0,\infty) \,.
\]
\end{lemma}
\medskip

\noindent
{\bf Proof.} We have
\begin{eqnarray*} 
 \int_{0}^{\phi (t) }  | K_1( r,t;-iM) |  dr 
 & \leq &
  \int_{0}^{\phi (t) } \Bigg| 4^{ \frac{iM}{H}} e^{-iM t} \left(\left(1+e^{-H t}\right)^2-(H r)^2\right)^{-i\frac{M}{H}-\frac{1}{2}} \\
&  &
\times  F \left(\frac{1}{2}+i\frac{M}{H},\frac{1}{2}+i\frac{M}{H};1;\frac{\left(1-e^{-H t} \right)^2-(r H)^2}{\left(1+e^{-H t} \right)^2-(r H)^2}\right)\Bigg|\, dr\,.
\end{eqnarray*}
Denote $ y:= e^{Ht}Hr$, $r=y/{zH}$ and $z:=e^{Ht} $, $ y = zHr$ and $y/ z =Hr $.  If $M$ is real, then
\[ 
 \int_{0}^{\phi (t) }  | K_1( r,t;-iM) |  dr 
  \leq  
 \frac{1}{ H} \int_{0}^{z-1}  \left( \left(z+1 \right)^2-y^2\right)^{-\frac{1}{2}}  \Bigg|F \left(\frac{1}{2}+i\frac{M}{H},\frac{1}{2}+i\frac{M}{H};1;\frac{\left(z-1 \right)^2- y ^2  }{\left(z+1 \right)^2-  y ^2}\right)\Bigg| \,dy\,.
\]
On the other hand, (see  \cite[Sec. 7]{Yag_Galst_CMP})
\begin{eqnarray}
\label{zeta}
\left| F \left(\frac{1}{2}+i\frac{M}{H},\frac{1}{2}+i\frac{M}{H};1;\zeta \right) \right| \leq  C_{M,H }   \left(1 - \ln  (1-\zeta  )     \right) ^{1- \sgn M} \quad \mbox{\rm for all}\quad \zeta  \in [0 ,1)  \,.
\end{eqnarray} 
According to Lemma 7.2~\cite{Yag_Galst_CMP} with $\rho =1$ if $M>0$, then
\[ 
 \int_{0}^{\phi (t) }  | K_1( r,t;-iM) |  dr 
  \leq  
  C_{M } \frac{1}{ H} \int_{0}^{z-1}  \left( \left(z+1 \right)^2-y^2\right)^{-\frac{1}{2}} dy \,.
\]
Then, for all $z>1$ the following inequality
\[ 
\int_{  0}^{ z  - 1} 
((z  + 1)^2  - r^2  )^{-\frac{1 }{2}}  
  d r  
    \leq  
C(z  - 1) (z  + 1)^{-1 }    
\]
implies
\[
 \int_{0}^{\phi (t) }  | K_1( r,t;-iM) | \, dr 
 \leq  
  C_{M } \frac{1}{ H} (z  - 1) (z  + 1)^{-1 }\,.
\]
If $M=0$ we obtain (see Lemma~7.2~\cite{Yag_Galst_CMP})
\[
\int_{0}^{z-1}  \left( \left(z+1 \right)^2-y^2\right)^{-\frac{1}{2}}   \Bigg|F \left(\frac{1}{2} ,\frac{1}{2} ;1;\frac{\left(z-1 \right)^2- y ^2  }{\left(z+1 \right)^2-  y ^2}\right)\Bigg|\, dy  
  \leq  
(1+\ln (z)) (z  - 1) (z  + 1)^{-1 }\,.
\]
The lemma is proved. \qed

\subsection{The linear equation  with   source and without potential}

Recall that in the case of large mass 
$
M^2:= \frac{m^2c^4}{h^2} - \frac{9H^2}{4 }\geq 0$ and  $M \geq 0 $\,.
 
\begin{theorem}
\label{T13.6}
For the solution of the problem
\[
\cases{
\psi_{tt} +3H \psi_t- e^{-2Ht}{\mathcal A}(x,\partial_x) \psi +\frac{m^2c^4}{h^2}\psi= \Psi, \quad   t >0,\cr
  u(x,0)=  0  \, , \quad u_t(x,0)=0\,,
}
\]
where supp\,$\Psi \subset \{(x,t) \in \R^3\times[0,\infty)\,|\,|x|>R_{ID}- c(1-e^{-tH})/H \,\}$ and $M \geq 0 $ one has
\[ 
\|\psi (x,t)\| _{H_{(s)}}  
   \leq 
 C_{M } e^{-\frac{3H}{2}t} \int_{ 0}^{t}\| \Psi(x, b)\|_{H_{(s)}} e^{\frac{3H}{2}b}  (1+ H (t-b  ))^{1- \sgn   M}\, db\quad \mbox{for all} \quad t >0\,.
\]
\end{theorem}
\medskip

\noindent
{\bf Proof.} 
The function  $u= u(x,t)$   defined by
\[ 
u(x,t) 
   =   
2   \int_{ 0}^{t} db
  \int_{ 0}^{\phi (t)- \phi (b)}  E(r,t;0,b;-iM )  v_f(x,r ;b) \, dr  
, \quad    t >0 ,
\]
where $\phi (t):= (1-e^{-Ht} )/H$,
 solves the problem (see, \cite{MN2015})
\begin{equation}
\label{13.17}
\cases{
u_{tt} - e^{-2Ht}{\mathcal A}(x,\partial_x)  u + M ^2 u= f, \quad  t >0,\cr
  u(x,0)=  0  \, , \quad u_t(x,0)=0\,.
}
\end{equation}

First, we prove that 
for the solution of the problem (\ref{13.17}), the following estimate
\[
\|u(x,t)\| _{H_{(s)}} \\
   \leq   
 C_{M } \int_{ 0}^{t}\| f(x, b)\|_{H_{(s)}} e^{-H (t-b  )}(e^{ H (t-b )}-1)(1+ H (t-b  ))^{1- \sgn   M}\, db 
\]
holds for all $t>0$. 
Indeed, it follows from (\ref{0.6}) that 
\begin{eqnarray*} 
E(r,t;0,t_0;-iM)
& :=  &
4^{ \frac{ iM}{H}} e^{ -iM  (  t_0+ t)} \left(\left(e^{-H t_0}+e^{-H t}\right)^2-(H r)^2\right)^{-i\frac{M}{H}-\frac{1}{2}}  \nonumber \\
&  &
\times  F \left(\frac{1}{2}+i\frac{ M}{H},\frac{1}{2}+i\frac{ M}{H};1;\frac{\left(-e^{-H t}+e^{-H t_0}\right)^2-(r H)^2}{\left(e^{-H t}+e^{-Ht_0}\right)^2-(r H)^2}\right)  \,.
\end{eqnarray*}
Then with the cutoff function $ \chi =\chi (x)$ we obtain
\[
\|u(x,t)\| _{H_{(s)}}
  =  
\|\chi u(x,t)\| _{H_{(s)}}  
 \, \leq  
\, 2 \|\chi\|_{H_{(s)}}  \int_{ 0}^{t} db
  \int_{ 0}^{\phi (t)- \phi (b)}  |E(r,t;0,b;M)|  \|v_f(x,r ;b)\|_{H_{(s)}} \, dr \,.
\]
Since (\ref{Blplq}),
we obtain
\begin{eqnarray*} 
\|u(x,t)\| _{H_{(s)}}  
&  \leq & 
C_\chi   \int_{ 0}^{t}\| f(x, b)\|_{H_{(s)}} \, db
  \int_{ 0}^{\phi (t)- \phi (b)}  \left(\left(e^{-H t_0}+e^{-H t}\right)^2-(H r)^2\right)^{-\frac{1}{2}}  \nonumber \\
&  &
\times  \left|F \left(\frac{1}{2}+i\frac{ M}{H},\frac{1}{2}+i\frac{ M}{H};1;\frac{\left(-e^{-H t}+e^{-H b}\right)^2-(r H)^2}{\left(e^{-H t}+e^{-Hb}\right)^2-(r H)^2}\right) \right|  \, dr \,.
\end{eqnarray*}
Consider the second integral with $z=e^{ H (t -b)} >1$, $y=e^{H t}Hr $,  and 
$\phi (t)- \phi (b)= (e^{-Hb}-e^{-Ht} )/{H}$. 
Then
\begin{eqnarray*}
&  &
 \int_{ 0}^{\frac{1}{H}(e^{-Hb}-e^{-Ht} )}  \left(\left(e^{-Hb}+e^{-H t}\right)^2-(H r)^2\right)^{-\frac{1}{2}}  \\
&  &
\times  \left|F \left(\frac{1}{2}+i\frac{ M}{H},\frac{1}{2}+i\frac{ M}{H};1;\frac{\left(-e^{-H t}+e^{-H b}\right)^2-(r H)^2}{\left(e^{-H t}+e^{-Hb}\right)^2-(r H)^2}\right) \right|  \, dr\\
& = &
 \frac{1}{H}\int_{ 0}^{z-1}  \left(\left(z+1\right)^2-y^2\right)^{-\frac{1}{2}}  \left|F \left(\frac{1}{2}+i\frac{ M}{H},\frac{1}{2}+i\frac{ M}{H};1;\frac{\left(z-1\right)^2-y^2}{\left(z+1\right)^2-y^2}\right) \right|  \,  dy\,.
\end{eqnarray*}
Next we use (\ref{zeta}) and obtain
\begin{eqnarray*}
&  &
 \int_{ 0}^{z}  \left(\left(e^{-Hb}+e^{-H t}\right)^2-(H r)^2\right)^{-\frac{1}{2}}   \left|F \left(\frac{1}{2}+i\frac{ M}{H},\frac{1}{2}+i\frac{ M}{H};1;\frac{\left(-e^{-H t}+e^{-H b}\right)^2-(r H)^2}{\left(e^{-H t}+e^{-Hb}\right)^2-(r H)^2}\right) \right|  \, dr\\
& \leq &
C_H  z^{-1}(z-1)(1+\ln z)^{1- \sgn M}\,.
\end{eqnarray*}
Hence,
\[ 
\|u(x,t)\| _{H_{(s)}}  
   \leq  
 C_{\chi, M } \int_{ 0}^{t}\| f(x, b)\|_{H_{(s)}} e^{-H (t-b  )}(e^{ H (t-b )}-1)(1+ H (t-b  ))^{1- \sgn M}\, db \,.
\]
Now we set
$u = e^{\frac{3H}{2}t}\psi  $    
and derive
\[
\|\psi (x,t)\| _{H_{(s)}} 
 \leq    
 C_{M }  \int_{ 0}^{t}\| \Psi(x, b)\|_{H_{(s)}} e^{-\frac{5H}{2}(t-b)} (e^{ H (t-b )}-1)(1+ H (t-b  ))^{1- \sgn M}\, db\,.
\]
The theorem is proved. 
\qed
\subsection{The integral equation. The global existence}

We study the Cauchy problem    (\ref{setCP})\&(\ref{CP})
    through the integral equation.
To define that integral equation, we  appeal to the operator 
\[
G:={\mathcal K}\circ {\mathcal EE}\,,
\]
where ${\mathcal EE}$ stands for the evolution (wave) equation in the exterior of BH in the universe without expansion as follows. For the function $f(x,t) $, we define
\[
v(x,t;b):= {\mathcal EE} [f](x,t;b)\,,
\]
where the function 
$v(x,t;b)$   
is a solution to the Cauchy problem 
\begin{eqnarray}
\label{1.6} 
&   &
\partial_t^2 v - {\mathcal A}(x,\partial_x)v =0, \quad x \in B^{ext}_{Sch} \subset {\mathbb R}^3, \quad t \geq 0, \\
\label{2.2a}
&  &
v(x,0;b)=f(x,b)\,, \quad v_t(x,0;b)= 0\,, \quad x \in B^{ext}_{Sch} \subset {\mathbb R}^3\,,  \quad b \geq 0,
\end{eqnarray} 
while ${\mathcal K}$ is introduced  by 
\begin{eqnarray}
\label{9.9}
{\mathcal K}[v]  (x,t) 
&  :=  &
2   e^{-\frac{3H}{2}t}\int_{ 0}^{t} db
  \int_{ 0}^{\phi(t)-\phi(b)} dr  \,  e^{\frac{3H}{2}b} v(x,r ;b) E(r,t; 0,b;-iM)  \,.
\end{eqnarray}
The kernel $ E(r,t; 0,b;M) $ is given by  (\ref{0.6}). Hence, 
\[
G[f]  (x,t)  
  =  
2   e^{-\frac{3H}{2}t}\int_{ 0}^{t} db
  \int_{ 0}^{ \phi(t)-\phi(b)} dr  \,  e^{\frac{3H}{2}b}\,{\mathcal EE} [f](x,r ;b) E(r,t; 0,b;-iM)  \,.
\]
Denote $\widetilde{C}^\ell ([0,T]; H_{(s)}) $ the complete subspace of $C^\ell([0,T]; H_{(s)})$ of all   functions $f=f(x,t)$ with supp\,$f \subset \{(x,t) \in \R^3\times[0,\infty)\,|\,|x|> R_{ID}- c(1-e^{-tH})/H \,\}$.  
According to  Section~\ref{S_FPS} and the theory of linear strictly hyperbolic equations with the smooth coefficients, for every $T>0$ the operator $G$ maps
\[
G \,:\, \widetilde{C}([0,T]; H_{(s)})\longrightarrow \widetilde{C}^2([0,T]; H_{(s)})
\]
 continuously.  
Thus, the Cauchy problem  (\ref{setCP})\&(\ref{13.12})  leads to the following integral equation
\begin{eqnarray} 
\label{5.1}
\psi  (x,t)
 = 
\psi  _{id}(x,t) + G[V  \psi ] (x,t)+
G[ F(\cdot)\Psi(\cdot ,\psi ) ] (x,t)    \,, 
\end{eqnarray}
where
\begin{eqnarray}
\label{psi0}
\psi_{id}  (x,t) 
 & =  & 
e^{-Ht} v_{\psi _0}  (x, \phi (t))
+ \,  e^{-\frac{3}{2}Ht}\int_{ 0}^{\phi (t)} \big(2  K_0( s,t;-iM)+ 3K_1( s,t;-iM)\big)v_{\psi _0}  (x, s) \,  ds  \nonumber \\
& &
+\, 2e^{-\frac{3}{2}Ht}\int_{0}^{\phi (t)}   v_{\psi _1 } (x,   s) 
  K_1(s,t;-iM)  \, ds
, \quad x \in B^{ext}_{Sch} \subset {\mathbb R}^n, \,\, t>0\, , 
\end{eqnarray}
and the function 
$v(x,t;b)$  of (\ref{9.9})  
is a solution to the Cauchy problem (\ref{1.6})\&(\ref{2.2a}),   
while $\phi (t):=  (1-e^{-Ht})/H $. 
Every solution to  the  Cauchy problem  (\ref{setCP})\&(\ref{13.12})  solves also the last integral equation with some function $\psi  _{id} (x,t)$, which is a solution to the problem for the linear equation without source and potential terms. We define a solution of the Cauchy problem (\ref{setCP})\&(\ref{13.12})  via integral equation (\ref{5.1}). 
\medskip

For the solution of the equation without self-interaction and potential terms, 
according to (\ref{psi0}) of Theorem~\ref{T13.2b}, we have
\[
\| \psi_{id}  (x,t)\| _{H_{(s)}} \\
  \leq   
C_{M,H} (1+t)^{1-\sgn (M)} e^{- Ht}\left(  \|  \psi _0 \|_{H_{(s)}}      
+\,  e^{-\frac{1}{2}Ht}\| \psi _1   \| _{H_{(s)}}    \right)
, \quad  \,\, t>0\, . 
\]
Consider the mapping $ S$ defined by the right-hand side of (\ref{5.1}):
\begin{equation}
\label{Sop}
S [\Phi] 
= \psi_{id}+G[V\Phi] +G[F\Psi(\Phi)]\,,
\end{equation}
where
\[
\psi_{id} \in X({R,H_{(s)},\gamma})\,.
\]
The operator $ S$ does not enlarge support of function $\Phi $ if supp$\,\Phi \subseteq  $ supp$\,\psi_{id}$.  
We claim that if $\Phi \in   X({R,H_{(s)},\gamma})$ with  $ \gamma \in (0,H)$,  and if supp~$\Phi \subseteq \{(x,t) \in \R^3\times[0,\infty)\,|\,|x|> R_{ID}- c(1-e^{-tH})/H \,\} $, then 
$S[\Phi] \in   X({R,H_{(s)},\gamma})$.     
 Moreover, $S$ is a contraction. Indeed, according to Theorem~\ref{T13.6} and condition
\[
\|V(x,t)  \Phi (t) \| _{H_{(s)}}\leq \varepsilon_0 \| \Phi (t) \| _{H_{(s)}}\,,
\] 
we obtain
\begin{eqnarray*} 
e^{ \gamma t} \|S [\Phi]\| _{H_{(s)}}  
&  \leq &  
e^{ \gamma t} \|\psi_{id}\| _{H_{(s)}} + e^{ \gamma t} \|G[V  \Phi] \| _{H_{(s)}}\\
&  &
+C_{M }  e^{ \gamma t}\int_{ 0}^{t}\left( \|  \Phi (x, b)\|_{H_{(s)}}\right)^{1+\alpha} e^{-\frac{3H}{2}(t-b)} (1+ H (t-b  ))^{1- \sgn M}\, db\,.
\end{eqnarray*}
First, we consider
\begin{eqnarray*}
&  & 
e^{ \gamma t}\int_{ 0}^{t}\left( \|  \Phi (x, b)\|_{H_{(s)}}\right)^{1+\alpha} e^{-\frac{3H}{2}(t-b)} (1+ H (t-b  ))^{1- \sgn M}\, db\\
&  \leq &  
\int_{ 0}^{t}\left(e^{ \gamma b} \|  \Phi (x, b)\|_{H_{(s)}}\right)^{1+\alpha} e^{ \gamma t- \gamma(1+\alpha)b} e^{-\frac{3H}{2}(t-b)}  (1+ H (t-b  ))^{1- \sgn M}\, db\,.
\end{eqnarray*}
If {$M>0$}, then
\begin{eqnarray*}
&  & 
e^{ \gamma t}\int_{ 0}^{t}\left( \|  \Phi (x, b)\|_{H_{(s)}}\right)^{1+\alpha} e^{-\frac{3H}{2}(t-b)} (1+ H (t-b  ))^{1- \sgn M}\, db\\
&  \leq & 
\left(\sup_{t \in [0,\infty) } e^{ \gamma t} \|  \Phi (x, t)\|_{H_{(s)}}\right)^{1+\alpha}\int_{ 0}^{t} e^{ \gamma t- \gamma(1+\alpha)b} e^{-\frac{3H}{2}(t-b)} \, db. 
\end{eqnarray*}
On the other hand, 
\begin{eqnarray*} 
\int_{ 0}^{t} e^{\gamma  t -(\alpha +1) b \gamma -\frac{3}{2} H (t-b)}  \, db  
& \leq &
\cases{\dsp \frac{-2}{ 3 H b -2(\alpha +1) b \gamma}e^{\gamma  t  -\frac{3}{2} Ht} \quad \mbox{if}\quad \gamma > \frac{3H}{2(\alpha +1)}\,,\cr
\dsp  \frac{2 e^{-\alpha  \gamma  t}-2 e^{\gamma  t-\frac{3 H t}{2}}}{3 H-2 (\alpha +1) \gamma  } < \frac{2 e^{-\alpha  \gamma  t}}{3 H-2 (\alpha +1) \gamma  }
\quad \mbox{if}\quad \gamma < \frac{3H}{2(\alpha +1)}\,,\cr
\dsp e^{\gamma  t  -\frac{3}{2} Ht}t\quad \mbox{if}\quad \gamma = \frac{3H}{2(\alpha +1)}<\frac{3H}{2}\,.}
\end{eqnarray*}
Hence, for $M>0$ we choose $ 0<\gamma  \leq \frac{3H}{2 }$ and  with $C(\gamma,H,\alpha)>0$ we obtain
\[
e^{ \gamma t}\int_{ 0}^{t}\left( \|  \Phi (x, b)\|_{H_{(s)}}\right)^{1+\alpha} e^{-\frac{3H}{2}(t-b)} (1+ H (t-b  ))^{1- \sgn M}\, db
   \leq   
C(\gamma,H,\alpha) \left(\sup_{t \in [0,\infty) } e^{ \gamma t} \|  \Phi (x, t)\|_{H_{(s)}}\right)^{1+\alpha} .
\]
If {$M=0$}, then
\begin{eqnarray*}
&  & 
e^{ \gamma t}\int_{ 0}^{t}\left( \|  \Phi (x, b)\|_{H_{(s)}}\right)^{1+\alpha} e^{-\frac{3H}{2}(t-b)} (1+ H (t-b  ))^{1- \sgn M}\, db\\
& \leq  & 
\left(\sup_{t \in [0,\infty) } e^{ \gamma t} \|  \Phi (x, t)\|_{H_{(s)}}\right)^{1+\alpha}e^{ \gamma t-  \frac{3H}{2}t} \int_{ 0}^{t} e^{ \frac{3H}{2}b- \gamma(1+\alpha)b  } (1+ H (t-b  )) \, db
\end{eqnarray*}
and we set  $ 0<\gamma  < \frac{3H}{2 }$  to obtain
\begin{eqnarray}
\label{13.18}  
e^{\gamma  t  -\frac{3}{2} Ht}\int_{ 0}^{t} e^{\frac{3}{2} H b -(\alpha +1) b \gamma }(1+ H (t-b  ))   \, db   
& \leq &
C_{\alpha,\gamma,H}   
\quad \mbox{\rm for all}\quad t \in [0,\infty) \,.  
\end{eqnarray}
Next, we consider the term with the potential $ V$ 
that  is analogous to the case of $\alpha=0$:
\begin{eqnarray*} 
 e^{ \gamma t} \|G[V  \Phi ] \| _{H_{(s)}} 
& \leq & 
\varepsilon_0 \left(\sup_{t \in [0,\infty) } e^{ \gamma t} \|  \Phi (x, t)\|_{H_{(s)}}\right) C_{\gamma,H}\cases{1 \quad \mbox{if}\quad \gamma <\frac{3}{2}H, \cr t^2 \quad \mbox{if}\quad \gamma =\frac{3}{2}H .}
\end{eqnarray*}
Thus, with any $\gamma$ such that $ 0<\gamma  \leq H$, we have
\begin{eqnarray*}
\left(\sup_{t \in [0,\infty) } e^{ \gamma t} \|S [\Phi]\| _{H_{(s)}} \right)
&  \leq & 
\varepsilon_0 C_{\gamma,H}\left(\sup_{t \in [0,\infty) } e^{ \gamma t} \| \Phi \| _{H_{(s)}} \right)  
+ \left(\sup_{t \in [0,\infty) } e^{ \gamma t} \|\Phi_{id}\| _{H_{(s)}} \right) \\
&  &
+C\left(\sup_{t \in [0,\infty) } e^{ \gamma t} \|  \Phi  \|_{H_{(s)}}\right)^{1+\alpha}
\quad \mbox{\rm for all}\quad   t \in [0,\infty) .
\end{eqnarray*}
If $ \psi_{id}(x,t)$ is generated by the initial data $\psi_0  (x) $ and   $\psi_1 (x) $, then
with $ \gamma \in(0, H ]$ we obtain
\begin{eqnarray*}  
e^{ \gamma t} \|S [\Phi]\| _{H_{(s)}}   
&  \leq &  
\varepsilon_0 C_{\gamma,H} \left(\sup_{t \in [0,\infty) } e^{ \gamma t} \|\Phi (t)\| _{H_{(s)}}\right)+C\left(\sup_{t \in [0,\infty) } e^{ \gamma t} \|  \Phi (x, t)\|_{H_{(s)}}\right)^{1+\alpha} \\
&  &
+ C_{H,\chi} (  (1+t)^{1-\sgn M}  \|\psi_0  (x)  \| _{H_{(s)}} +  e^{ -\frac{1}{2}H t}\|  \psi_1 (x)\| _{H_{(s)}}  )\quad \mbox{\rm for all}\quad   t \in [0,\infty) .
\end{eqnarray*}
For $ \varepsilon_0 C_{\gamma,H}<1$  and $M>0$ it follows 
\begin{eqnarray*} 
 \left(\sup_{t \in [0,\infty) } e^{ \gamma t} \|\Phi (t)\| _{H_{(s)}}\right) 
&  \leq &  
 \frac{1}{1-\varepsilon_0 C_{\gamma,H}  } C_{H,\chi}   (  \|\psi_0  (x)  \| _{H_{(s)}} +  e^{ -\frac{1}{2}H t}\|  \psi_1 (x)\| _{H_{(s)}}  ) \\
&  &
+\frac{1}{1-\varepsilon_0 C_{\gamma,H}  }C\left(\sup_{t \in [0,\infty) } e^{ \gamma t} \|  \Phi (x, t)\|_{H_{(s)}}\right)^{1+\alpha}\quad \mbox{\rm for all}\quad   t \in [0,\infty) .
\end{eqnarray*}
Then we choose initial data,  $\varepsilon_0 $, and $ R$ such that 
\[
 \frac{1}{1-\varepsilon_0 C_{\gamma,H}  }C_{H,\chi}   (  \|\psi_0  (x)  \| _{H_{(s)}} +  \|  \psi_1 (x)\| _{H_{(s)}}  )+ \frac{1}{1-\varepsilon_0 C_{\gamma,H}  }C     R ^{\alpha+1}  <R \,.
 \]
For $M=0$, we set  $ \gamma \in(0, H )$, appeal to (\ref{LpLq}),   and come to the same conclusion.  
 \medskip
 
To prove that $S$ is a contraction mapping, we obtain the contraction property from 
\[
\sup_{t \in [0,\infty) }e^{ \gamma t} \|S[\Phi_1 ](\cdot,t) -  S[\Phi_2  ](\cdot,t) \|_{H_{(s)}({\mathbb R}^n) }
 \leq  
CR(t) ^{\alpha } d(\Phi ,\Psi )\,,  
\]
where 
\begin{equation}
\label{Rt}
\displaystyle R(t):= \max\{ \sup_{0\leq \tau \leq t }   e^{ \gamma t}\| \Phi_1  (\cdot ,\tau ) \| _{H_{(s)}({\mathbb R}^n) }
, \sup_{0\leq \tau \leq t }  e^{ \gamma t} \| \Phi_2 (\cdot ,\tau ) \| _{H_{(s)}({\mathbb R}^n) } \}\leq R. 
\end{equation}
Indeed,  due to Theorem~\ref{T13.6}, we have
\begin{eqnarray*}
&  &
 e^{ \gamma t}\| S[\Phi_1 ](\cdot , t) -  S[\Phi_2 ](\cdot , t) \|_{H_{(s)}({\mathbb R}^n) } \\ 
& \leq &
 e^{ \gamma t}\| V(\cdot, t) [ 
   \Phi_1  -   \Phi_2   ](\cdot, t)
\|_{H_{(s)}({\mathbb R}^n) }   
+e^{ \gamma t}\| G[ \,
F( \Psi( \Phi_1  ) -  \Psi(\Phi_2   ) ) \, ](\cdot, t)
\|_{H_{(s)}({\mathbb R}^n) }  \,.
\end{eqnarray*}
Consider the term
\begin{eqnarray*}
&  &
e^{ \gamma t}\| G[ \,F
( \Psi( \Phi_1  ) -  \Psi(\Phi_2   ) ) \, ](\cdot, t)
\|_{H_{(s)}({\mathbb R}^n) }\\
& \le  &
C_{M,\alpha } e^{ \gamma t}\int_{ 0}^{t} e^{-\frac{5H}{2}(t-b)} (e^{ H (t-b )}-1)(1+ H (t-b  ))^{1- \sgn M}\\
&  &
\times \|\Phi_1 (\cdot , b) -\Phi_2 (\cdot , b) \|_{H_{(s)}({\mathbb R}^n) } 
\left( \| \Phi_1  (\cdot ,  b) \|_{H_{(s)}({\mathbb R}^n) } ^\alpha  
+ \| \Phi_2  (\cdot ,  b)  \|_{H_{(s)}({\mathbb R}^n) }^\alpha 
\right) \,db    \,.
\end{eqnarray*}
Thus,  taking into account  the last estimate  and a definition of the metric, we obtain  
\begin{eqnarray*}
&  &
e^{ \gamma t}\| G[ \,
F( \Psi( \Phi_1  ) -  \Psi(\Phi_2   ) ) \, ](\cdot, t)
\|_{H_{(s)}({\mathbb R}^n) }\\ 
& \le  &
C_{M,\alpha } e^{ \gamma t}\int_{ 0}^{t} e^{-\frac{5H}{2}(t-b)} (e^{ H (t-b )}-1)(1+ H (t-b  ))^{1- \sgn M}\\
&  &
\times e^{ -\gamma b}e^{ -\gamma \alpha  b}\Big( \max_{0 \le \tau  \leq b } e^{ \gamma \tau}\|\Phi_1 (\cdot , \tau ) -\Phi_2 (\cdot , \tau ) \|_{H_{(s)}({\mathbb R}^n) }  \Big) \\
&  &
\times 
\left( \left( e^{  \gamma    b}\| \Phi_1  (\cdot ,  b) \|_{H_{(s)}({\mathbb R}^n) } \right)^\alpha  
+ \left(e^{  \gamma    b}\| \Phi_2  (\cdot ,  b)  \|_{H_{(s)}({\mathbb R}^n) }\right)^\alpha 
\right) \,db  \\
& \le  &
C_{M,\alpha } d(\Phi ,\Psi )
R(t)^\alpha  \int_{ 0}^{t} e^{ \gamma t}e^{-\frac{3H}{2}(t-b)}   (1+ H (t-b  ))^{1- \sgn M}e^{ -\gamma b}e^{ -\gamma \alpha  b}\,db   \,,
\end{eqnarray*}
and, consequently, by (\ref{13.18}), we arrive at
\[ 
e^{ \gamma t}\| G[ \,
F( \Psi( \Phi_1  ) -  \Psi(\Phi_2   ) ) \, ](\cdot, t)
\|_{H_{(s)}({\mathbb R}^n) }  
  \le  
C_{M,\alpha } d(\Phi_1 ,\Phi_2 )
R^\alpha     \,.
\]
Similarly, for the term with potential, since $\alpha =0$,  we obtain 
\begin{eqnarray*}
e^{ \gamma t}\| V(\cdot, t) [ 
   \Phi_1  -   \Phi_2   ](\cdot, t)
\|_{H_{(s)}({\mathbb R}^n) }\leq \varepsilon_0 C_{M,\alpha,V }
   d(\Phi_1 ,\Phi_2  )\,.
\end{eqnarray*}
Finally,   
\[
\| S[\Phi_1 ](\cdot , t) -  S[\Phi_2 ](\cdot , t) \|_{H_{(s)}({\mathbb R}^n) }  
  \le  
 \varepsilon_0 C_{M,\alpha,V }  d(\Phi_1 ,\Phi_2 )+ C_{M,\alpha}      R(t)^\alpha   d(\Phi_1 ,\Phi_2 )  \,.
\]
Then we choose $\|\psi_{id}\| _{H_{(s)}} \leq \varepsilon $ and $ R$ such that $ \varepsilon_0C_{M,\alpha}  +C_{M,\alpha}    R ^\alpha  <1 $. 
Banach's fixed point theorem completes the proof of theorem.
\hfill $\square$

 \section{The semilinear equation with  small mass. Proof of Theorem~\ref{TIE}: Existence of global solution}

\setcounter{equation}{0}
\renewcommand{\theequation}{\thesection.\arabic{equation}}

 \subsection{Linear equation without source and potential terms}
\begin{theorem} 
\label{T13.2}
For every given $s \in {\mathbb R}$, 
the solution  $ \psi = \psi (x,t)$ of the Cauchy problem for the equation 
\[
    \frac{\partial^2 \psi }{\partial t^2}
+ 3H        \frac{\partial \psi }{\partial t}
-e^{-2Ht}{\mathcal A}(x,\partial_x) \psi+   \frac{m^2c^4}{h^2} \psi =0\,, \quad x \in \R^3,\,\, t \in [0,\infty)\nonumber
\]
with the initial conditions (\ref{13.12})
and $\Re M= \Re \left(\frac{9H^2 }{4 }-\frac{m^2c^4}{h^2} \right)^{1/2} \in \left(0,\frac{H}{2}\right)$  
satisfies the following  estimate
\[ 
 \|  \psi (x,t) \| _{H_{(s)}}  
  \leq   
C_{m,s}  e^{   -H t}\left(   \| \psi _0   
\|_{H_{(s)}}
+ (1- e^{-Ht}) \|\psi _1  
\|_{H_{(s)}} 
 \right) \quad \mbox{for all} \quad t \in (0,\infty)\,.
\]

If   
$\Re M= \Re \left(\frac{9H^2 }{4 }-\frac{m^2c^4}{h^2} \right)^{1/2}>  H/2$ or $M=H/2$,  
then    the solution  $ \psi  = \psi (x,t)$ of the Cauchy problem  
satisfies the following  estimate
\begin{eqnarray*}   
\| \psi  (x,t) \| _{H_{(s)}} 
   \leq 
C_{m,s}   e^{(\Re M-\frac{3H}{2})t}     \left( \|\psi _0 \|_{H_{(s)}} + (1- e^{-Ht}) \|\psi _1 \|_{H_{(s)}}\right) \quad \mbox{for all} \quad t \in (0,\infty)\,.
\end{eqnarray*}  
\end{theorem}
\medskip

\noindent
{\bf Proof.} {The case of $M=H/2 $} is an evident consequence of (\ref{Blplq}),(\ref{0.6}),(\ref{K0MH}),(\ref{K1MH}), and the representation (\ref{psi}), 
where 
\begin{equation}
\label{10.1}
 E(r, t ; 0, b ; H/2)  =  \frac{1}{2}   e^{\frac{H }{2}(t+b)}, \quad 
K_0(r,t;H/2)
  =  
-\frac{1}{4} H e^{\frac{H}{2} t}, \quad K_1(r,t;H/2)
 =  
\frac{1}{2} e^{\frac{H}{2} t}\,.
\end{equation}
Hence,
\[
\psi (x,t)
   =  
e ^{-Ht} v_{\psi_0}  (x, \phi (t))  
+   e ^{- Ht} H \int_{ 0}^{\phi (t)} v_{\psi_0 } (x,  s)
  ds   
+   e ^{- Ht}    \int_{0}^{\phi (t) }  v_{\psi_1 } (x,  s)
 ds\,. \nonumber 
\]
Then (\ref{Blplq}) 
and $ \phi (t):= (1-e^{-Ht})/H \leq 1/H$ imply
\begin{eqnarray*}
&  &
\|  \psi (x,t) \| _{H_{(s)}}\\
&  \leq  &
e ^{-Ht}\|  v_{\psi_0}  (x, \phi (t))\| _{H_{(s)}}  
+   e ^{- Ht} H \int_{ 0}^{\phi (t)}\|  v_{\psi_0 } (x,  s)\| _{H_{(s)}}
  ds   
+  e ^{- Ht}    \int_{0}^{\phi (t) } \|  v_{\psi_1 } (x,  s)\| _{H_{(s)}}
 ds\\
&  \leq  &
2e ^{-Ht}\|  {\psi_0}  (x )\| _{H_{(s)}}  
+  e ^{- Ht}  \|   {\psi_1 } (x )\| _{H_{(s)}}  (1-e ^{- Ht})/H\,.
\end{eqnarray*}
Now  we consider the {case of $M\not=H/2 $}. 
First 
we consider the case of $\psi _1=0 $. Then 
\begin{eqnarray*} 
\psi (x,t) 
&  =  &
 e ^{-Ht} v_{\psi_0}  (x, \phi (t))  
+   e ^{-\frac{3}{2}Ht}\int_{ 0}^{1}  \left[ 2K_0( \phi (t)s,t;M)    
+   3HK_1( \phi (t)s,t;M)  \right] v_{\psi_0 } (x,  \phi (t)s)
 \phi (t) ds 
\end{eqnarray*}
and, consequently, 
\begin{eqnarray} 
\label{2.9}
\|  \psi  (x,t) \| _{H_{(s)}}  
&   \leq &  
e^{-Ht} \|  v_{\psi _0}  (x, \phi (t)) \| _{H_{(s)}}   \\
&  &
+ \,  e^{-\frac{3}{2}Ht}\int_{ 0}^{1} \| v_{\psi _0}  (x, \phi (t)s)\| _{H_{(s)}} \big|2  K_0(\phi (t)s,t;M)+ 3K_1(\phi (t)s,t;M)\big|\phi (t)\,  ds\,. \nonumber
\end{eqnarray}
Further,  for the solution $v = v (x,t)$ of the Cauchy problem (\ref{2.1})\&(\ref{2.2}),  
 one has  the estimate (\ref{Blplq}).
Hence, 
\[  
e^{-Ht} \| v_{\psi _0}  (x, \phi (t)s) \| _{H_{(s)}}
  \leq  
C  e^{-Ht}\|\psi _0 \|_{H_{(s)}}  \quad \mbox{\rm for all} \,\,  t >0, \,\, s \in [0,1]\,.
\]  
For the second  term of (\ref{2.9}), we obtain
\begin{eqnarray*} 
&  &
e^{-\frac{3}{2}Ht}\int_{ 0}^{1} \|  v_{\psi _0}  (x, \phi (t)s) \| _{H_{(s)}} \big|2  K_0(\phi (t)s,t;M)+ 3K_1(\phi (t)s,t;M)\big|\phi (t)\,  ds\\
& \leq &
 \|\psi _0 \|_{H_{(s)}}e^{-\frac{3}{2}Ht} \int_{ 0}^{1} \left( \big|2  K_0(\phi (t)s,t;M)\big|+ 3\big|K_1(\phi (t)s,t;M)\big|\right) \phi (t)\,  ds \,.
\end{eqnarray*}
Next, we have to estimate the following two integrals of the last inequality:
\[ 
 \int_{ 0}^{1} 
\big|  K_i(\phi (t)s,t;M)\big|   \phi (t)\,  ds, \quad i=0,1 \,,
\]
where   $t>0$. To complete the estimate of the second term  of (\ref{2.9}), we are going to apply the following two lemmas
with $a=0$.
\begin{lemma}
\label{L2.3}
Let $  a>-1 $, $\Re M>0$,  and $\phi (t)= (1-e^{-Ht})/H$. Then 
\[ 
 \int_{ 0}^{1} \phi (t)^{a} s^{a}
\big|  K_1(\phi (t)s,t;M)\big|   \phi (t)\,  ds  
  \leq  
 C_M  e^{-a Ht}(e^{Ht }-1)^{a+1} (e^{Ht }+1)^{ \frac{\Re M}{H}-1}  \quad \mbox{  for all}  \quad t>0\,.
\]  
In particular,
\[ 
 \int_{ 0}^{1} \phi (t)^{a} s^{a}
\big|  K_1(\phi (t)s,t;M)\big|   \phi (t)\,  ds  
  \leq  
 C_{M,a}   e^{\Re M   t}  \quad \mbox{  for large}  \quad t\,.
\]
\end{lemma}
\medskip

\noindent
{\bf Proof.} By the definition of the kernel $K_1$, we obtain
\begin{eqnarray*} 
 \int_{ 0}^{1} \phi (t)^{a} s^{a}
\big|  K_1(\phi (t)s,t;M)\big|   \phi (t)\,  ds  
& \leq &
 4 ^{-\frac{\Re M}{H}} e^{ \Re Mt } \int_{ 0}^{(1-e^{-Ht})/H} r^{a} 
 \left(\left(1+e^{-H t}\right)^2-(H r)^2\right)^{\frac{\Re M}{H}-\frac{1}{2}} \\
&  &
\times  
\Bigg| F \left(\frac{1}{2}-\frac{M}{H},\frac{1}{2}-\frac{M}{H};1;\frac{\left(1-e^{-H t} \right)^2-(r H)^2}{\left(1+e^{-H t} \right)^2-(r H)^2}\right)\Bigg|    \,  dr \,.
\end{eqnarray*}
Then we use the substitution  $He^{Ht}r=y$, $r=e^{-Ht}y/H$ as follows  
\begin{eqnarray*}  
 \int_{ 0}^{1} \phi (t)^{a} s^{a}
\big|  K_1(\phi (t)s,t;M)\big|   \phi (t)\,  ds  
& \leq &
H^{-a-1} 4 ^{-\frac{\Re M}{H}} e^{-  \Re M t -aHt}  \int_{ 0}^{ e^{Ht}-1  }   y^a 
 \left(\left(e^{ H t}+1\right)^2-y^2\right)^{\frac{\Re M}{H}-\frac{1}{2}} \\
&  &
\times  
\Bigg| F \left(\frac{1}{2}-\frac{M}{H},\frac{1}{2}-\frac{M}{H};1;\frac{\left(e^{ H t}- 1\right)^2-y^2}{\left(e^{ H t}+ 1\right)^2-y^2}\right)\Bigg|    \,   dy \,.
\end{eqnarray*}
On the other hand, for $\Re M>0$, we have (see \cite[Section~A]{JMP2019})
\[  
\left|F\left(\frac{1}{2}-\frac{M}{H},\frac{1}{2}-\frac{M}{H} ;1; \zeta  \right) \right|     \leq C_M \quad \mbox{\rm for all} \quad \zeta  \in [0,1)\,,
\]
where 
\[ 
\zeta := \frac{ ( e^{Ht }-1)^2 -y^2 }{( e^{Ht }+1)^2 -y^2 }  \in [0,1) \quad \mbox{\rm for all} \quad y \in [0, e^{Ht }-1] \quad \mbox{\rm and all} \quad t>0.
\]
Hence,
\[
 \int_{ 0}^{1} \phi (t)^{a} s^{a}
\big|  K_1(\phi (t)s,t;M)\big|   \phi (t)\,  ds \\
  \leq  
 C_{H,M}4 ^{-\frac{\Re M}{H}} e^{-  \Re M t -aHt}  \int_{ 0}^{ e^{ Ht}-1  }   y^a 
 \left(\left(e^{ H t}+1\right)^2-y^2\right)^{\frac{\Re M}{H}-\frac{1}{2}}  \,   dy\,. 
\]
If we denote $z:= e^{Ht}$, then for $ M>0$ we have
\[
\int_{ 0}^{z-1}  y^{a} 
  \big((z+1)^2 -   y  ^2\big)^{-\frac{1}{2}+M    }     \,   dy  \\
  =  
\frac{1}{1+a}(z-1)^{1+a} (z+1)^{2 M-1} F\left(\frac{1+a}{2},\frac{1}{2}-M;\frac{3+a}{2};\frac{(z-1)^2}{(z+1)^2}\right) ,
 \]  
where $a>-1$ and $z\geq 1$. Hence, for $\Re M>0$ we have
\[  
\int_{ 0}^{1} \phi (t)^{a} s^{a}
\big|  K_1(\phi (t)s,t;M)\big|   \phi (t)\,  ds   
   \leq  
C_M  e^{-\Re Mt -aHt }(e^{Ht}-1)^{a+1} (e^{Ht}+1)^{2 \frac{\Re M}{H}-1}  
\]
for all  $ t>0 $. The lemma is proved. \hfill $\square$

\begin{lemma}
\label{L14.3}
Let $ a>-1$,  $\Re M>0$, and  $\phi (t)= (1-e^{-Ht})/H$. Then 
\[  
\int_{ 0}^{1} \phi (t)^{a} s^{a}
\big|  K_0(\phi (t)s,t;M)\big|   \phi (t)\,  ds 
   \leq  
  C_{M,a}  
(e^{Ht}-1)^{a+1}e^ { -aHt}  \times \cases{   (e^{Ht}+1)^{ -\frac{1}{2}} \quad \mbox{\rm if} \quad \Re M<H/2\,, \cr
 e^ { \Re M t}  (e^{Ht}+1)^{ -1} \quad \mbox{\rm if} \quad \Re M > H/2,} 
\] 
for all   $t>0$. In particular,
\[
 \int_{ 0}^{1} \phi (t)^{a} s^{a}
\big|  K_0(\phi (t)s,t;M)\big|   \phi (t)\,  ds  
  \leq   C_{M,a}  
\times \cases{   e^ {  \frac{1}{2}Ht}  \quad \mbox{\rm if} \quad \Re M<H/2\,, \cr
 e^ { (\Re  M- H) t}    \quad \mbox{\rm if} \quad \Re M> H/2\,,}  
\]
for large  $t$.
\end{lemma}
\medskip

\noindent
{\bf Proof.} By definition of $K_0$, we obtain
\begin{eqnarray*} 
&  & 
 \int_{ 0}^{1} \phi (t)^{a} s^{a}
\big|  K_0(\phi (t)s,t;M)\big|   \phi (t)\,  ds \\
& \leq  &
4^{-\frac{\Re M}{H}} e^{ t \Re M} \int_{ 0}^{(1-e^{-Ht})/H} r^{a}
\frac{   \left(\left(e^{-H t}+1\right)^2-H^2 r^2\right)^{\frac{ \Re M}{H}-\frac{1}{2}} }{\left(1-e^{-H t}\right)^2-H^2 r^2}\\
&  &
\times\Bigg|\Bigg[ \left(H e^{-H t}-H+M e^{-2 H t}-M-H^2 M r^2\right)  F \left(\frac{1}{2}-\frac{M}{H},\frac{1}{2}-\frac{M}{H};1;\frac{\left(1 -e^{-H t}\right)^2-H^2 r^2}{\left(1+e^{-H t}\right)^2-H^2 r^2}\right)\\
&  &
+\left(\frac{H}{2}+M\right) \left(H^2r^2- e^{-2 H t}+1\right)  F \left(-\frac{1}{2}-\frac{M}{H},\frac{1}{2}-\frac{M}{H};1;\frac{\left(1- e^{-H t}\right)^2-H^2 r^2}{\left(1+e^{-H t}\right)^2-H^2 r^2}\right)\Bigg]\Bigg|    \,  d r \,.
\end{eqnarray*} 
Now we   make the change    $r=e^{-Ht}y H^{-1}$ in the last integral and obtain
\begin{eqnarray*} 
&  & 
 \int_{ 0}^{1} \phi (t)^{a} s^{a}
\big|  K_0(\phi (t)s,t;M)\big|   \phi (t)\,  ds \\
& \leq  &
C e^{ -t \Re M-aHt}  \int_{ 0}^{ e^{Ht}-1 } y^{a}
  \frac{   \left(\left(e^{H t}+1\right)^2-  y  ^2\right)^{\frac{ \Re M}{H}-\frac{1}{2}} }{\left( e^{H t}-1\right)^2-   y  ^2}\\
&  &
\times\Bigg|\Bigg[ \left(  H e^{ H t}- e^{2Ht}H+  M  - e^{2Ht}M-  M y ^2\right)  F \left(\frac{1}{2}-\frac{M}{H},\frac{1}{2}-\frac{M}{H};1;\frac{\left(e^{ H t}- 1\right)^2-y^2}{\left(e^{ H t}+ 1\right)^2-y^2}\right)\\
&  &
+\left(\frac{H}{2}+M\right) \left(  y ^2- 1+ e^{2Ht}\right)  F \left(-\frac{1}{2}-\frac{M}{H},\frac{1}{2}-\frac{M}{H};1;\frac{\left(e^{ H t}- 1\right)^2-y^2}{\left(e^{ H t}+ 1\right)^2-y^2}\right)\Bigg]\Bigg|    \,  d y  \,.
\end{eqnarray*}  
Then we denote $z=e^{Ht}$ and derive
\begin{eqnarray*} 
&  & 
 \int_{ 0}^{1} \phi (t)^{a} s^{a}
\big|  K_0(\phi (t)s,t;M)\big|   \phi (t)\,  ds \\
& \leq  &
C z^{ -t \frac{\Re M}{H}  -a t}  \int_{ 0}^{ z-1  } y^{a}
  \frac{   \left(\left(z+1\right)^2-  y  ^2\right)^{\frac{ \Re M}{H}-\frac{1}{2}} }{\left( z-1\right)^2-   y  ^2}\\
&  &
\times\Bigg|\Bigg[ \left(  H z- z^{2}H+  M  - z^{2 }M-  M y ^2\right)  F \left(\frac{1}{2}-\frac{M}{H},\frac{1}{2}-\frac{M}{H};1;\frac{\left(z- 1\right)^2-y^2}{\left(z+ 1\right)^2-y^2}\right)\\
&  &
+\left(\frac{H}{2}+M\right) \left(  y ^2- 1+ z^{2 }\right)  F \left(-\frac{1}{2}-\frac{M}{H},\frac{1}{2}-\frac{M}{H};1;\frac{\left(z- 1\right)^2-y^2}{\left(z+ 1\right)^2-y^2}\right)\Bigg]\Bigg|    \,  d y  \,.
\end{eqnarray*}  
To complete the proof of Lemma~\ref{L14.3}, we apply the following statement.
\begin{proposition}
If\, $a >-1 $ and $\Re M>0$, then 
\begin{eqnarray} 
&  &    
\label{14.24}
\int_{ 0}^{ z-1  } y^{a}
  \frac{   \left(\left(z+1\right)^2-  y  ^2\right)^{\frac{ \Re M}{H}-\frac{1}{2}} }{\left( z-1\right)^2-   y  ^2}\\
&  &
\times\Bigg|\Bigg[ \left(  H z- z^{2}H+  M  - z^{2 }M-  M y ^2\right)  F \left(\frac{1}{2}-\frac{M}{H},\frac{1}{2}-\frac{M}{H};1;\frac{\left(z- 1\right)^2-y^2}{\left(z+ 1\right)^2-y^2}\right)\nonumber\\
&  &
+\left(\frac{H}{2}+M\right) \left(  y ^2- 1+ z^{2 }\right)  F \left(-\frac{1}{2}-\frac{M}{H},\frac{1}{2}-\frac{M}{H};1;\frac{\left(z- 1\right)^2-y^2}{\left(z+ 1\right)^2-y^2}\right)\Bigg]\Bigg|    \,  d y  \nonumber\\
& \leq &
 C_{M,H,a}(z-1)^{1+a}\times \cases{  (z+1)^{\frac{ \Re M}{H}-\frac{1}{2}} \quad \mbox{\rm if} \quad \Re M<H/2\,,\cr
 (z+1)^{2\frac{ \Re M}{H}-1} \quad \mbox{\rm if}  \quad \Re M> H/2\,.} \nonumber
\end{eqnarray}  
\end{proposition}
\medskip

\noindent
{\bf Proof.}   We follow the arguments that have been used in the proof of Lemma~7.4~\cite{Yag_Galst_CMP}. For $\Re M>0$,   both hypergeometric functions are bounded.  We divide the domain of integration into  two zones, 
\begin{eqnarray*} 
Z_1(\varepsilon, z) 
& := &
\left\{ (z,y) \,\Big|\, \frac{ (z-1)^2 -y^2   }{ (z+1)^2 -y^2 } \leq \varepsilon, \quad  0 \leq y \leq z-1,\,\, z \geq 1 \right\} \,,\\   
Z_2(\varepsilon, z) 
& := &
\left\{ (z,y) \,\Big|\, \varepsilon \leq  \frac{ (z-1)^2 -y^2   }{ (z+1)^2 -y^2 }\,,\,\,  \quad  0 \leq y \leq z-1,\,\, z \geq 1  \right\} \,,
\end{eqnarray*}
and then  split the integral into  two parts,
\[
\int_{ 0}^{z-1} \star  \, dr 
  =  
\int_{ (z,r) \in Z_1(\varepsilon, z)   }  \star  \, dr 
+ \int_{ (z,r) \in Z_2(\varepsilon, z)  }  \star   \, dr \,,
\] 
where $ \star $ denotes the integrand in (\ref{14.24}). In the first zone $Z_1(\varepsilon, z) $, we have
\begin{eqnarray*}
F\Big(\frac{1}{2}-\frac{M}{H},\frac{1}{2}-\frac{M}{H};1; \frac{ (z-1)^2 -y^2   }{ (z+1)^2 -y^2 }   \Big)  
 &  =  &
 1 + \left( \frac{1}{2}-\frac{M}{H} \right)^2\frac{ (z-1)^2 -y^2   }{ (z+1)^2 -y^2 } \\
&  &  
+ O\left(\left( \frac{ (z-1)^2 -y^2   }{ (z+1)^2 -y^2 }\right)^2\right),  \\  
F\Big(-\frac{1}{2}-\frac{M}{H},\frac{1}{2}-\frac{M}{H};1; \frac{ (z-1)^2 -y^2   }{ (z+1)^2 -y^2 }   \Big)  
 &  =  & 
 1 - \left( \frac{1}{4}-\left(\frac{M}{H} \right)^2 \right)\frac{ (z-1)^2 -y^2   }{ (z+1)^2 -y^2 }  \\
&  &   
+ O\left(\left( \frac{ (z-1)^2 -y^2   }{ (z+1)^2 -y^2 }\right)^2\right)   . 
\end{eqnarray*} 
 We use the last formulas to estimate the term containing hypergeometric functions: 
\begin{eqnarray*} 
&  &    
\Bigg|\Bigg[ \left(  H z- z^{2}H+  M  - z^{2 }M-  M y ^2\right)  F \left(\frac{1}{2}-\frac{M}{H},\frac{1}{2}-\frac{M}{H};1;\frac{\left(z- 1\right)^2-y^2}{\left(z+ 1\right)^2-y^2}\right)\\
&  &
+\left(\frac{H}{2}+M\right) \left(  y ^2- 1+ z^{2 }\right)  F \left(-\frac{1}{2}-\frac{M}{H},\frac{1}{2}-\frac{M}{H};1;\frac{\left(z- 1\right)^2-y^2}{\left(z+ 1\right)^2-y^2}\right)\Bigg]\Bigg|      \\
& \leq &  
\Bigg|\Bigg[ \left(  H z- z^{2}H+  M  - z^{2 }M-  M y ^2\right) \\
&  &
\times \Bigg[  1 + \left( \frac{1}{2}-\frac{M}{H} \right)^2\frac{ (z-1)^2 -y^2   }{ (z+1)^2 -y^2 }   
+ O\left(\left( \frac{ (z-1)^2 -y^2   }{ (z+1)^2 -y^2 }\right)^2\right)\Bigg]  \\
&  &
+\left(\frac{H}{2}+M\right) \left(  y ^2- 1+ z^{2 }\right)\Bigg[  1 - \left( \frac{1}{4}-\left(\frac{M}{H} \right)^2 \right)\frac{ (z-1)^2 -y^2   }{ (z+1)^2 -y^2 }  + O\left(\left( \frac{ (z-1)^2 -y^2   }{ (z+1)^2 -y^2 }\right)^2\right) \Bigg]\Bigg|  \\
& \leq &  
\Bigg|\Bigg[\frac{1}{2}  H \left(y^2-(z-1)^2\right)\\
&  &
-\frac{1}{8} H \left(1-\frac{2 M}{H}\right) \left( 2 \frac{M}{H} \left(3 y^2+z^2+2 z-3\right) +\left(y^2+3 z^2-2 z-1\right)\right)\left( \frac{ (z-1)^2 -y^2   }{ (z+1)^2 -y^2 }\right)\\
&  &
+\frac{1}{2}  H \left(y^2-(z-1)^2\right)O\left(\left( \frac{ (z-1)^2 -y^2   }{ (z+1)^2 -y^2 }\right)^2\right) \Bigg]\Bigg| \,.      
\end{eqnarray*} 
Thus, on the left-hand side of (\ref{14.24}), we have to consider the following two integrals, which can be easily  estimated,
\[
A_1
   :=  
\int_{ (z,y) \in Z_1(\varepsilon, z)  }  y^{a}  \big((z+1)^2 - y^2\big)^{\frac{\Re M}{H} -\frac{1}{2}   }\,  dy   
  \,, \quad
A_2
  :=  
z^2  \int_{  (z,y) \in Z_1(\varepsilon, z) }  y^{a}  \big((z+1)^2 - y^2\big)^{ \frac{\Re M}{H} -\frac{3}{2}    }  dy 
   ,
\] 
 for all  $z \in [1,\infty)$. Indeed, for $A_1$ we obtain 
 \begin{eqnarray*}
A_1
& \leq   &
\int_{ 0 }^{z-1}  y^{a}  \big((z+1)^2 - y^2\big)^{\frac{\Re M}{H} -\frac{1}{2}   }\,  dy \\
& = &
\frac{1}{1+a}(z-1)^{1+a} (z+1)^{2 \frac{\Re M}{H}-1} F\Big(\frac{1+a}{2},\frac{1}{2}-\Re M;\frac{3+a}{2};\frac{(z-1)^2}{(z+1)^2}\Big) \\
& \leq &
C_{M,a} (z-1)^{1+a} (z+1)^{2 \frac{\Re M}{H}-1}\,.
\end{eqnarray*} 
 Similarly, if $\Re M>0$,  then
  \begin{eqnarray}
  \label{1.9}
A_2
& \leq   &
z^2  \int_{ 0 }^{z-1}   y^{a}  \big((z+1)^2 - y^2\big)^{ \frac{\Re M}{H} -\frac{3}{2}    }  dy \nonumber \\
& = &
z^2\frac{1}{1+a}(z-1)^{1+a} (z+1)^{2 \frac{\Re M}{H}-3} F\Big(\frac{1+a}{2},\frac{3}{2}-\frac{\Re M}{H};\frac{3+a}{2};\frac{(z-1)^2}{(z+1)^2}\Big) \,.
\end{eqnarray}

Here and henceforth, if $A$ and $B$ are two non-negative quantities, we use $A \lesssim  B$ to denote the statement that $A\leq CB $ for some absolute constant $C>0$. 

It suffices to consider the case of real valued $M$. Then  \cite[(A5)]{JMP2019} and  (\ref{1.9}), in the case of $\Re M<H/2$, imply 
 \[
A_2
  \lesssim   
z^2\frac{1}{1+a}(z-1)^{1+a} (z+1)^{2 \frac{\Re M}{H}-3}  z^{\frac{1}{2}-\frac{\Re M}{H}}  
  \lesssim    
  (z-1)^{1+a} (z+1)^{\frac{\Re M}{H}-\frac{1}{2} }\,.
\] 
In the case of $M\geq  H/2$ due to ~\cite[(A5)]{JMP2019},  
we derive
 \[
A_2
  \lesssim  
z^2(z-1)^{1+a} (z+1)^{2\frac{\Re M}{H}-3}   
  \lesssim   
 (z-1)^{1+a} (z+1)^{2 \frac{\Re M}{H}-1}  \,.
\] 
 Finally, for the integral over  the first zone $Z_1(\varepsilon, z) $, we  obtain
\[
\int_{ (z,r) \in Z_1(\varepsilon, z)   }  \star\,  dr
  \lesssim   
(z-1)^{1+a}\times \cases{  (z+1)^{\frac{\Re M}{H}-\frac{1}{2}} \quad \mbox{\rm if} \quad \Re M<H/2\,,\cr
 (z+1)^{2\frac{\Re M}{H}-1} \quad \mbox{\rm if}  \quad \Re M> H/2\,.}
\] 
In the second zone, we have
\[
0< \varepsilon \leq  \frac{ (z-1)^2 -r^2   }{ (z+1)^2 -r^2 } < 1 \quad \mbox{\rm and}  \quad 
\frac{ 1  }{ (z-1)^2 -r^2 }  \leq  \frac{ 1   }{ \varepsilon[(z+1)^2 -r^2] }\,.
\]
Then, the hypergeometric functions for $\Re M>0$ obey the estimates
\[
\left| F\Big(-\frac{1}{2}-\frac{  M}{H},\frac{1}{2}-\frac{  M}{H};1; \zeta      \Big) \right|  \leq C \,\,  \mbox{\rm and}  \,\,   
\left| F\Big(\frac{1}{2}-\frac{ M}{H},\frac{1}{2}-\frac{  M}{H};1; \zeta   \Big) \right|  \leq C   \,\, \, \mbox{\rm for all}\,\, \,  \zeta  \in [\varepsilon ,1) .
\]
This allows us to  estimate  the integral over the second zone as follows: 
\begin{eqnarray*} 
&  &     
\int_{  (z,y) \in Z_2(\varepsilon, z) } y^{a}
  \frac{   \left(\left(z+1\right)^2-  y  ^2\right)^{\frac{ \Re M}{H}-\frac{1}{2}} }{\left( z-1\right)^2-   y  ^2}\\
&  &
\times\Bigg|\Bigg[ \left(  H z- z^{2}H+  M  - z^{2 }M-  M y ^2\right)  F \left(\frac{1}{2}-\frac{M}{H},\frac{1}{2}-\frac{M}{H};1;\frac{\left(z- 1\right)^2-y^2}{\left(z+ 1\right)^2-y^2}\right)\nonumber\\
&  &
+\left(\frac{H}{2}+M\right) \left(  y ^2- 1+ z^{2 }\right)  F \left(-\frac{1}{2}-\frac{M}{H},\frac{1}{2}-\frac{M}{H};1;\frac{\left(z- 1\right)^2-y^2}{\left(z+ 1\right)^2-y^2}\right)\Bigg]\Bigg|    \,  d y  \nonumber\\
& \lesssim &
z^2   \int_{  (z,y) \in Z_2(\varepsilon, z) }  y^{a}
    \big((z+1)^2 - y^2\big)^{  \frac{\Re M}{H}  -\frac{3}{2}  } \,  d  y  \\ 
& \lesssim  &
 z^2   \int_{ 0}^{z-1}  y^{a}
    \big((z+1)^2 - y^2\big)^{ \frac{\Re M}{H}  -\frac{3}{2}  } \,  d  y  \,. 
\end{eqnarray*}  
Then we apply (\ref{1.9}) and Lemma~A.1\cite{JMP2019}:    
\[
z^2   \int_{  (z,y) \in Z_2(\varepsilon, z) }  y^{a}
    \big((z+1)^2 - y^2\big)^{ \frac{\Re M}{H} -\frac{3}{2}  } \,  d  y    
   \lesssim  
 (z-1)^{1+a}\times \cases{  (z+1)^{\frac{\Re M}{H}-\frac{1}{2}} \quad \mbox{\rm if} \quad \Re M<H/2\,,\cr
 (z+1)^{2\frac{\Re M}{H}-1} \quad \mbox{\rm if}  \quad \Re M> H/2}  
\] 
for all  $z \in [1,\infty)$.  Finally, for the integral over  the second zone $Z_2(\varepsilon, z) $, we obtain
\[
\int_{ (z,r) \in Z_2(\varepsilon, z)   }  \star\,  dr
  \lesssim  
(z-1)^{1+a}\times \cases{  (z+1)^{\frac{\Re M}{H}-\frac{1}{2}} \quad \mbox{\rm if} \quad \Re M<H/2\,,\cr
 (z+1)^{2\frac{\Re M}{H}-1} \quad \mbox{\rm if}  \quad \Re M> H/2\,.}
\]
The rest of the proof is a repetition of the above-used arguments.   
Thus, the  proposition is proved.    \qed 
\medskip 

\noindent
{\bf Completion of the proof of Theorem~\ref{T13.2}.} Thus, if $\psi _1=0 $, then from (\ref{2.9}) we derive
\begin{eqnarray*} 
\| \psi  (x,t) \| _{H_{(s)}}  
&   \leq &  
e^{-Ht} \|  v_{\psi _0}  (x, \phi (t)) \| _{H_{(s)}}   \\
&  &
+ \,  e^{-\frac{3}{2}Ht}\int_{ 0}^{1} \| v_{\psi _0}  (x, \phi (t)s)\| _{H_{(s)}}  \big|2  K_0(\phi (t)s,t;M)+ 3K_1(\phi (t)s,t;M)\big|\phi (t)\,  ds\\
&   \lesssim  & 
  e^{-Ht}\|\psi _0 \|_{H_{(s)}}  \\
&  &
+  \|\psi _0 \|_{H_{(s)}}e^{-\frac{3}{2}Ht}(e^{Ht }-1) 
 \left(   (e^{Ht }+1)^{ \frac{\Re M}{H}-1}    +    
 \cases{    (e^{Ht}+1)^{ -\frac{1}{2}} \quad \mbox{\rm if} \quad \Re M<H/2\,, \cr
 e^ { \frac{\Re M}{H}t}  (e^{Ht}+1)^{ -1} \quad \mbox{\rm if} \quad \Re M>  H/2}  \right) .  
\end{eqnarray*}
 In particular, for large $t$, we obtain
\[  
\|  \psi  (x,t) \| _{H_{(s)}}  
   \lesssim  
    \|\psi _0 \|_{H_{(s)}}\left(  e^{-Ht} 
+ e^{-\frac{1}{2}Ht}     \left[   e^{( \Re M -H)t }     +    
 \cases{   e^{ -\frac{1}{2}Ht} \quad \mbox{\rm if} \quad \Re M<H/2\,, \cr
 e^ { \Re M t}  e^{-Ht} \quad \mbox{\rm if} \quad \Re M>H/2}  \right]\right) .  
\]
In the case of $\psi _0=0 $, we have
\begin{eqnarray*}
\|\psi (x,t) \| _{H_{(s)}}
& \leq  & 
2
e ^{-\frac{3}{2}Ht}\int_{0}^{\phi (t) }  \|v_{\psi_1 } (x,  s)\| _{H_{(s)}}
 | K_1( s,t;M)|   ds \\
& \leq  & 
C_{s}e ^{-\frac{3}{2}Ht} \|\psi_1(x)\| _{H_{(s)}}\int_{0}^{\phi (t) } 
 | K_1( s,t;M)|   ds ,\quad  x \in \R^3 \,,\,\, t >0. \nonumber 
\end{eqnarray*}
Due to Lemma~\ref{L2.3},   we obtain
\[
\|\psi (x,t) \| _{H_{(s)}}
   \lesssim    
e ^{(\frac{\Re M} {H}-\frac{3}{2})Ht}  (1-e^{-Ht }) \|\psi_1(x)\| _{H_{(s)}}\quad   t >0\,. \nonumber 
\]
The theorem is proved. \hfill $\square$

\subsection{The linear equation with source term and without potential}

We consider equations with $m \in {\mathbb C} $. 
This is why in   this section we focus on the cases of $\Re M >0 $ and complex valued $M$. Thus, we are also interested in the Higgs boson equation,   in  massive scalar fields, and in the tachyons having $m^2 < 0$.

\begin{theorem} 
\label{T14.4} 
Let $ \psi = \psi (x,t)$ be a solution of the Cauchy problem for the equation
\[
    \frac{\partial^2 \psi }{\partial t^2}
+ 3H        \frac{\partial \psi }{\partial t}
-e^{-2Ht}{\mathcal A}(x,\partial_x)\psi+   \frac{m^2c^4}{h^2} \psi =f\,, \quad x \in \R^3,\,\, t \in [0,\infty)\nonumber
\]
with the initial conditions
$
\psi (x,0)=0$, $\,\partial_t \psi (x,0)=0$, 
where supp\,$f \subset \{(x,t) \in \R^3\times[0,\infty)\,|\,|x|> R_{ID}- c(1-e^{-tH})/H \,\}$  and 
$
\Re M= \Re \left(\frac{9H^2 }{4 }-\frac{m^2c^4}{h^2} \right)^{1/2} \,.
$ 
  
Then  the solution $ \psi = \psi (x,t)$ for $0< \Re M <H/2 $    
satisfies the following  estimate:
\[  
 \|    \psi (x ,t) \| _{H_{(s)} }
  \leq  
Ce^{-Ht}   
 \int_{ 0}^{t} e^{Hb} 
   \|f(x,b)\|_{H_{(s)} }\,db   \quad \mbox{for all}\quad   t>0.
\]  
If either $\Re M> H/2$ or $M=H/2$, then 
\[
  \| \psi (x ,t)  \| _{H_{(s)} }
  \leq   
  e^{ (\Re M   -\frac{3}{2}H)t} \int_{ 0}^{t} e ^{-(\Re M -\frac{3H}{2})   b}\|  f(x, b)  \| _{H_{(s)} }db    \quad \mbox{for all}\quad   t>0\,.   
\]
\end{theorem}
\medskip

\noindent
{\bf Proof.}   The case of $M=H/2 $  is an evident consequence of  the representations (\ref{psi}) 
 and (\ref{10.1}).   Indeed,
\begin{eqnarray*} 
\| \psi (x,t)\|_{H_{(s)} } 
&  \leq  &
e ^{- Ht}   \int_{ 0}^{t}e ^{2H b}  db
  \int_{ 0}^{\phi (t)- \phi (b)}   \| v_f(x,r ;b)\|_{H_{(s)} } \, dr \\ 
&   \lesssim    &
e ^{- Ht}   \int_{ 0}^{t}e ^{2H b}  db
  \int_{ 0}^{\phi (t)- \phi (b)}   \|  f(x, b)\|_{H_{(s)} } \, dr \\ 
&   \lesssim    &
e ^{- Ht}   \int_{ 0}^{t}e ^{ H b}   \|  f(x, b)\|_{H_{(s)} }
  db
\quad \mbox{for all}\quad   t>0. \nonumber 
\end{eqnarray*}
For the {case of $M\not=H/2 $}, 
  according to (\ref{Blplq}), we can write
\[ 
\| v(x,r ;b)   \| _{H_{(s)} }  
\leq C  \|f(x,b)\|_{H_{(s)} }  \quad \mbox{\rm for all} \quad  r \in [0,1/H], \quad b \geq 0\,.
\]
Hence, from (\ref{psi}), due to (\ref{0.6}), we derive
\begin{eqnarray*} 
\| \psi (x,t)  \| _{H_{(s)} }
&  \lesssim  &
e ^{-\frac{3}{2}Ht}  e^{ \Re M  t } \int_{ 0}^{t} e ^{\frac{3H}{2}b}  e^{ \Re M  b}\|  f(x, b)  \| _{H_{(s)} }db
  \int_{ 0}^{\phi (t)- \phi (b)}  \left(\left(e^{-H b}+e^{-H t}\right)^2-(H r)^2\right)^{\frac{\Re M}{H}-\frac{1}{2}}  \nonumber \\
&  &
\times  \left| F \left(\frac{1}{2}-\frac{M}{H},\frac{1}{2}-\frac{M}{H};1;\frac{\left(-e^{-H t}+e^{-H b}\right)^2-(r H)^2}{\left(e^{-H t}+e^{-Hb}\right)^2-(r H)^2}\right) \right | \, dr\,.
\end{eqnarray*}
Following the outline of the proof of Lemma~\ref{L2.3}, we set     $r=e^{-Ht}y H^{-1}$  and,  from the last inequality,  obtain 
\begin{eqnarray} 
&  &
\| \psi (x,t)  \| _{H_{(s)} } \nonumber\\
&  \lesssim  &
\label{10.5b}
e ^{-\frac{3}{2}Ht}  e^{ \Re M  t } e^{-Ht}e^{-2H(\frac{\Re M}{H}-\frac{1}{2}) t}\int_{ 0}^{t} e ^{\frac{3H}{2}b}  e^{ \Re M  b}\|  f(x, b)  \| _{H_{(s)} }db   \\
&  &
 \times  \int_{ 0}^{e^{H(t-b)}-1}  \left(\left(e^{H(t -b)}+1\right)^2-y^2\right)^{\frac{\Re M}{H}-\frac{1}{2}}    \left| F \left(\frac{1}{2}-\frac{M}{H},\frac{1}{2}-\frac{M}{H};1;\frac{\left( e^{H(t- b)}- 1\right)^2-y^2}{\left(e^{H(t- b)}+ 1\right)^2-y^2}\right) \right |  dy  . \nonumber 
\end{eqnarray}
In order to estimate the second integral, we apply Lemma~A.5~\cite{JMP2019}  with $z=e^{H(t-b)} >1$ and $a=0$.
Hence, the   estimate \cite[(A.7)]{JMP2019} 
\[
\int_{ 0}^{ z- 1} \,     \Big(( z+1)^2 - y^2\Big)^{-\frac{1}{2}+\Re\frac{M}{H}    } \left| F\Big(\frac{1}{2}-\frac{M}{H}   ,\frac{1}{2}-\frac{M}{H}  ;1; 
\frac{ (z-1)^2 -y^2 }{(z+1)^2 -y^2 } \Big) \right| dy 
   \leq   
C_M (z-1)   z^{\Re\frac{M}{H}-\frac{1}{2} } 
\] 
implies 
\[ 
\| \psi (x,t)  \| _{H_{(s)} }  
   \lesssim   
e ^{-\frac{3}{2}Ht}  e^{ \Re M  t } e^{-Ht}e^{-2H(\frac{\Re M}{H}-\frac{1}{2}) t}\int_{ 0}^{t} e ^{\frac{3H}{2}b}  e^{ \Re M  b}\|  f(x, b)  \| _{H_{(s)} }(e^{H(t-b)}-1)   e^{H(t-b) (\Re\frac{M}{H}-\frac{1}{2} )}\,db , 
\]
that is,  the following estimate
\[ 
\| \psi (x,t)  \| _{H_{(s)} }  
  \lesssim  
e ^{-Ht}    \int_{ 0}^{t} e ^{Hb}  \|  f(x, b)  \| _{H_{(s)} }db, \quad 0< \Re M<H/2\,.
\]

\ndt
For the case of  $\Re M>H/2$,  we apply \cite[(A.8)]{JMP2019}
\[
\int_{ 0}^{ z- 1} \, y^{a}   \Big(( z+1)^2 - y^2\Big)^{-\frac{1}{2}+\frac{\Re M}{H}    } \left| F\Big(\frac{1}{2}-\frac{M}{H}   ,\frac{1}{2}-\frac{ M}{H}  ;1; 
\frac{ (z-1)^2 -y^2 }{(z+1)^2 -y^2 } \Big) \right| dy 
  \lesssim  
 (z-1)^{1+a} (z+1)^{2\frac{\Re M}{H}  -1} ,
\]
 and from  (\ref{10.5b}), we obtain
\[ 
\| \psi (x,t)  \| _{H_{(s)} }  
   \lesssim   
e ^{-\frac{3}{2}Ht}  e^{ \Re M  t } e^{-Ht}e^{-2H(\frac{\Re M}{H}-\frac{1}{2}) t}\int_{ 0}^{t} e ^{\frac{3H}{2}b}  e^{ \Re M  b}\|  f(x, b)  \| _{H_{(s)} }  e^{H(t -b) 2 \frac{\Re M}{H} }\,db, 
\]
that is,  the estimate 
\[ 
\| \psi (x,t)  \| _{H_{(s)} }  
   \lesssim  
  e^{ (\Re M   -\frac{3}{2}H)t} \int_{ 0}^{t} e ^{-(\Re M -\frac{3H}{2})   b}\|  f(x, b)  \| _{H_{(s)} }db, \quad  \Re M>H/2\,.
\] 
The theorem is proved. \qed

\subsection{Global solution   to semilinear equation. Proof of Theorem~\ref{TIE}}

We are going to apply  Banach's fixed-point theorem.  
In order to estimate nonlinear terms, we use the   Lipschitz condition (${\mathcal L}$).
First, we consider the integral equation  (\ref{5.1}),  
where the function $ \psi _{id} (x,t) \in C([0,\infty);H_{(s)})$ is given. 
The operator $G$ and the structure of the nonlinear term determine the solvability of the integral equation (\ref{5.1}). 
\medskip

\noindent
 (i) In this case   $0< \Re M<H/2$. Consider the mapping (\ref{Sop}), 
where the function $
\psi_{id}  
$ 
is generated by initial data, that is, by  (\ref{psi}). We have 
\begin{eqnarray*} 
\psi_{id}  (x,t) 
&  =  &
 e ^{-Ht} v_{\psi_0}  (x, \phi (t)) 
+ \, e ^{-\frac{3}{2}Ht}\int_{ 0}^{\phi (t)}  \left[ 2K_0( s,t;M)    
+   3HK_1( s,t;M)  \right] v_{\psi_0 } (x,  s)
  ds  \nonumber \\
&  &
+\, 2e ^{-\frac{3}{2}Ht}\int_{0}^{\phi (t) }  v_{\psi_1 } (x,  s)
  K_1( s,t;M)   ds,\quad    t >0\,.
\end{eqnarray*}
The operator $S$ does not enlarge the support of function $\Phi $ if supp$\,\Phi \subseteq  $ supp$\,\psi_{id}$. We claim that if $ \Phi \in X({R,H_{(s)},\gamma})$ with   $\gamma \in [0,H] $ and if  supp~$\Phi \subseteq \{(x,t) \in \R^3\times[0,\infty)\,|\,|x|> R_{ID}- c(1-e^{-tH})/H \,\} $, then 
$S[\Phi] \in   X({R,H_{(s)},\gamma})$.  Moreover, $S$
is a contraction, provided that   $\varepsilon  $, $\varepsilon_0  $, and $R$ are sufficiently small. 

Consider the case of  $\Re M <H/2$. 
First, we note that due to Theorem~\ref{T13.2}
\[ 
e^{\gamma t}  \|  \psi_{id}  (x,t) \| _{H_{(s)}}  
  \leq   
C_{m,s} e^{(\gamma   -H )t}\Big(  \| \psi _0   
\|_{H_{(s)}}
+  \|\psi _1  
\|_{H_{(s)}} 
 \Big) \leq   
\varepsilon C_{m,s} e^{(\gamma   -H )t} \quad \mbox{for all} \quad t >0\,.
\]
Further,  due to Theorem~\ref{T14.4}, we obtain
\[
\|S [\Phi](t)\| _{H_{(s)}} 
   \leq   
\|\psi_{id}\| _{H_{(s)}}+\|G[V  \Phi (t)] \| _{H_{(s)}}
+C_{M } e^{-Ht} \int_{ 0}^{t}e^{Hb}\left( \|  \Phi (x, b)\|_{H_{(s)}}\right)^{1+\alpha}  \, db\,.
\]
Then, for $\gamma \in {\mathbb R} $, we have 
\begin{eqnarray*} 
e^{\gamma t} \|S [\Phi](t)\| _{H_{(s)}}  
&  \leq & 
e^{\gamma t} \|\psi_{id}\| _{H_{(s)}}+e^{\gamma t} \|G[V  \Phi (t)] \| _{H_{(s)}}\\
&  &
+C_{M } \left( \sup_{\tau  \in [0,\infty)} e^{\gamma  \tau  }     \| \Phi (\cdot ,\tau  ) \| _{H_{(s)}}   \right)^{\alpha  +1} e^{\gamma t-Ht} \int_{ 0}^{t}e^{Hb}e^{ - \gamma (\alpha +1)b}  \, db\,.
\end{eqnarray*}
For $   \gamma \in [0,   H]  $ and $\alpha >0 $, the following function is bounded  
\begin{eqnarray} 
\label{Int}
&  &
 e^{ \gamma  t- H t} \int_{ 0}^{t} e^{H  b}  e^{ - \gamma (\alpha +1)b} \,db \leq C_{\alpha,\gamma,H}  \quad \mbox{\rm for all}\quad  t \in [0,\infty)\,.
\end{eqnarray}
Consequently,
\[
e^{\gamma t} \|S [\Phi](t)\| _{H_{(s)}} 
   \leq  
e^{\gamma t} \|\psi_{id}\| _{H_{(s)}}+e^{\gamma t} \|G[V  \Phi (t)] \| _{H_{(s)}}
+C_{M }C_{\alpha,\gamma,H}  \left( \sup_{\tau  \in [0,\infty)} e^{\gamma  \tau  }     \| \Phi (\cdot ,\tau  ) \| _{H_{(s)}}   \right)^{\alpha  +1}  .
\]
Further, for $\gamma \in[0,H)$, according to condition (${\mathcal V})$ and the finite propagation speed property, 
 we have
\[ 
e^{\gamma t} \|G[V  \Phi (t)] \| _{H_{(s)}} 
  \leq  
 \varepsilon_0 C_{M }\frac{1}{H-\gamma}\left( \sup_{\tau  \in [0,\infty)} e^{\gamma  \tau  }     \| \Phi (\cdot ,\tau  ) \| _{H_{(s)}}   \right)  \,,
\]
and, consequently,  
\begin{eqnarray*}
e^{\gamma t} \|S [\Phi](t)\| _{H_{(s)}} 
&  \leq & 
e^{\gamma t} \|\psi_{id}\| _{H_{(s)}}+ \varepsilon_0 C_{M }\frac{1}{H-\gamma}\left( \sup_{\tau  \in [0,\infty)} e^{\gamma  \tau  }     \| \Phi (\cdot ,\tau  ) \| _{H_{(s)}}   \right)\\
&  &
+C_{M }C_{\alpha,\gamma,H} \left( \sup_{\tau  \in [0,\infty)} e^{\gamma  \tau  }     \| \Phi (\cdot ,\tau  ) \| _{H_{(s)}}   \right)^{\alpha  +1} \,.
\end{eqnarray*}
By Theorem~\ref{T13.2}, for $\varepsilon_0 C_{M } <H-\gamma  $  it follows 
\begin{eqnarray*} 
&  &
 \left( \sup_{\tau  \in [0,\infty)} e^{\gamma  \tau  }     \| \Phi (\cdot ,\tau  ) \| _{H_{(s)}}   \right)\\
&  \leq & 
\frac{1}{(1- \varepsilon_0 C_{M }\frac{1}{H-\gamma})} C_{m,s} \left(         \| \psi _0   
\|_{H_{(s)}}
+   \|\psi _1  
\|_{H_{(s)}} 
  \right) 
+ C_{M } \frac{1}{(1- \varepsilon_0 C_{M }\frac{1}{H-\gamma})}\left( \sup_{\tau  \in [0,\infty)} e^{\gamma  \tau  }     \| \Phi (\cdot ,\tau  ) \| _{H_{(s)}}   \right)^{\alpha  +1} \,.
\end{eqnarray*}
Thus, the last inequality proves that the operator $S$ maps $X({R,s,\gamma})$ into itself if $\varepsilon_0  $, $\varepsilon  $, and $R$ are sufficiently small, namely, if
\[
\frac{1}{(1- \varepsilon_0 C_{M }\frac{1}{H-\gamma})} C_{m,s} 
\varepsilon  +C_{M } \frac{1}{(1- \varepsilon_0 C_{M }\frac{1}{H-\gamma})} R^{\alpha +1} < R  .
\]
To prove that $S$ is a contraction mapping, we derive the contraction property from 
\[
\sup_{t \in [0,\infty) }e^{ \gamma t} \|S[\Phi_1 ](\cdot,t) -  S[\Phi_2  ](\cdot,t) \|_{H_{(s)}({\mathbb R}^n) }
 \leq  
CR(t) ^{\alpha } d(\Phi ,\Psi )\,,  
\]
where $R(t) $ is defined in (\ref{Rt}). 
Indeed,    we have
\begin{eqnarray*}
&  &
 e^{ \gamma t}\| S[\Phi_1 ](\cdot , t) -  S[\Phi_2 ](\cdot , t) \|_{H_{(s)}({\mathbb R}^n) } \\ 
& \leq &
 e^{ \gamma t}\| V(\cdot, t) [ 
   \Phi_1  -   \Phi_2   ](\cdot, t)
\|_{H_{(s)}({\mathbb R}^n) }   
+e^{ \gamma t}\| G[ \,
F( \Psi( \Phi_1  ) -  \Psi(\Phi_2   ) ) \, ](\cdot, t)
\|_{H_{(s)}({\mathbb R}^n) }  \,.
\end{eqnarray*}
For the second term, due to Theorem~\ref{T14.4}, we obtain
\begin{eqnarray*}
&  &
e^{ \gamma t}\| G[ \,
F( \Psi( \Phi_1  ) -  \Psi(\Phi_2   ) ) \, ](\cdot, t)
\|_{H_{(s)}({\mathbb R}^n) } \\
& \leq &
e^{ \gamma t} e^{ -Ht} \int_{ 0}^{t} e ^{Hb}  \|\Phi_1 (\cdot , b) -\Phi_2 (\cdot , b) \|_{H_{(s)}({\mathbb R}^n) } 
\left( \| \Phi_1  (\cdot ,  b) \|_{H_{(s)}({\mathbb R}^n) } ^\alpha  
+ \| \Phi_2  (\cdot ,  b)  \|_{H_{(s)}({\mathbb R}^n) }^\alpha 
\right) \,db\,.
\end{eqnarray*}
Thus,  taking into account  the last estimate  and the definition of the metric, we obtain 
\begin{eqnarray*}
&  &
e^{ \gamma t}\| G[ \,
F( \Psi( \Phi_1  ) -  \Psi(\Phi_2   ) ) \, ](\cdot, t)
\|_{H_{(s)}({\mathbb R}^n) } \\
& \leq &
e^{ \gamma t} e^{ -Ht} \int_{ 0}^{t} e ^{H b}  e^{ -\gamma b}e^{ -\gamma \alpha  b}\Big( \max_{0 \le \tau  \leq b } e^{ \gamma \tau}\|\Phi_1 (\cdot , \tau ) -\Phi_2 (\cdot , \tau ) \|_{H_{(s)}({\mathbb R}^n) }  \Big) \\
&  &
\times 
\left( \left( e^{  \gamma    b}\| \Phi_1  (\cdot ,  b) \|_{H_{(s)}({\mathbb R}^n) } \right)^\alpha  
+ \left(e^{  \gamma    b}\| \Phi_2  (\cdot ,  b)  \|_{H_{(s)}({\mathbb R}^n) }\right)^\alpha 
\right) \,db\\
& \leq &
C_{M,\alpha } d(\Phi ,\Psi )
R(t)^\alpha  e^{ (\gamma   -H)t} \int_{ 0}^{t} e ^{(H  -\gamma    -\gamma \alpha ) b}\, db\,.
\end{eqnarray*} 
Consequently, by (\ref{Int}), the following inequality holds:
\[ 
e^{ \gamma t}\| G[ \,
F( \Psi( \Phi_1  ) -  \Psi(\Phi_2   ) ) \, ](\cdot, t)
\|_{H_{(s)}({\mathbb R}^n) }  
  \leq  
C_{\alpha,H, M} d(\Phi ,\Psi )
R(t)^\alpha   \,.
\] 
Similarly, for the term with potential, since  $\gamma \in [0,H) $,  we obtain 
\[
e^{ \gamma t}\| V(\cdot, t) [ 
   \Phi_1  -   \Phi_2   ](\cdot, t)
\|_{H_{(s)}({\mathbb R}^n) } \leq \varepsilon_0 C_{M,\alpha,V }
   d(\Phi_1 ,\Phi_2  )\,.
\]
Finally,   
\[
\| S[\Phi_1 ](\cdot , t) -  S[\Phi_2 ](\cdot , t) \|_{H_{(s)}({\mathbb R}^n) }  
  \le  
 \varepsilon_0 C_{M,\alpha,V }  d(\Phi_1 ,\Phi_2 )+ C_{\alpha,H, M}      R(t)^\alpha   d(\Phi_1 ,\Phi_2 )  \,.
\]
Then we choose $\|\Phi_{id}\| _{H_{(s)}} \leq \varepsilon $ and $ R$ such that $ \varepsilon_0C_{M,\alpha,V}  +C_{\alpha,H,M}    R ^\alpha  <1 $. 
Banach's fixed point theorem completes the proof of 
the case of (i).
\medskip

\ndt
(ii) We claim that if $\Re M \in [H/2,3H/2)$, then the operator $ S \,: \, X({R,H_{(s)},\gamma})\longrightarrow X({R,H_{(s)},\gamma})$ of (\ref{Sop}) 
with   $\gamma =\frac{1}{\alpha +1}(\frac{3}{2}H- \Re M -\delta )>0$ is a contraction, provided that   $\varepsilon_0  $, $\varepsilon  $, and $R$ are sufficiently small.  By Theorem~\ref{T14.4}, we obtain
\[
\|S [\Phi](t)\| _{H_{(s)}} 
  \leq  
\|\psi_{id}\| _{H_{(s)}}+\|G[V  \Phi (t)] \| _{H_{(s)}}
+C_{M }   e^{ (\Re M   -\frac{3}{2}H)t} \int_{ 0}^{t} e ^{-(\Re M -\frac{3H}{2})   b}\left( \|  \Phi (x, b)\|_{H_{(s)}}\right)^{1+\alpha}  \, db .
\]
Then, for $\gamma \in {\mathbb R} $, we have 
\begin{eqnarray*} 
e^{\gamma t} \|S [\Phi](t)\| _{H_{(s)}}  
&  \leq & 
e^{\gamma t} \|\psi_{id}\| _{H_{(s)}}+e^{\gamma t} \|G[V  \Phi (t)] \| _{H_{(s)}}\\
&  &
+C_{M } \left( \sup_{\tau  \in [0,\infty)} e^{\gamma  \tau  }     \| \Phi (\cdot ,\tau  ) \| _{H_{(s)}}   \right)^{\alpha  +1} e^{\gamma t}   e^{ (\Re M   -\frac{3}{2}H)t} \int_{ 0}^{t} e ^{-(\Re M -\frac{3H}{2})   b}e^{ - \gamma (\alpha +1)b}  \, db\,.
\end{eqnarray*}
Further,
\[
 e^{\gamma t}   e^{ (\Re M   -\frac{3}{2}H)t}\int_{ 0}^{t} e ^{-(\Re M -\frac{3H}{2})   b}e^{ - \gamma (\alpha +1)b}  \, db 
  \leq  
C
\cases{ 
 e ^{   - \gamma \alpha t}  
\quad \mbox{if}\quad   \Re M -\frac{3H}{2}     + \gamma (\alpha +1)<0,\cr
 e^{[\gamma  + (\Re M   -\frac{3}{2}H)]t}  
\quad \mbox{if}\quad   \Re M -\frac{3H}{2}     + \gamma (\alpha +1)>0,\cr
t e^{-\gamma  \alpha  t} \quad \mbox{if}\quad   \Re M -\frac{3H}{2}  + \gamma (\alpha +1)=0\,. }
\]
Hence, for $\Re M \in [H/2,3H/2)$, $\alpha >0$, and $\gamma   \leq   \frac{1}{\alpha +1}(\frac{3}{2}H- \Re M )$, we have
\[
 e^{\gamma t}   e^{ (\Re M   -\frac{3}{2}H)t}\int_{ 0}^{t} e ^{-(\Re M -\frac{3H}{2})   b}e^{ - \gamma (\alpha +1)b}  \, db 
\leq C\,.
\] 

Further, for $\Re M \in [H/2,3H/2)$, according to condition
$\|V(x,t)  \Phi (t) \| _{H_{(s)}}\leq \varepsilon_0 \| \Phi (t) \| _{H_{(s)}}$, 
 we have
\begin{eqnarray*} 
e^{\gamma t} \|G[V  \Phi (t)] \| _{H_{(s)}} 
& \leq &
 C_{M } e^{\gamma t} e^{ (\Re M   -\frac{3}{2}H)t} \int_{ 0}^{t} e ^{-(\Re M -\frac{3H}{2})   b}\| V  \Phi (t)(x, b)\|_{H_{(s)}}  \, db \\
& \leq &
 \varepsilon_0 C_{M }C_{\alpha, \gamma, H}\left( \sup_{\tau  \in [0,\infty)} e^{\gamma  \tau  }     \| \Phi (\cdot ,\tau  ) \| _{H_{(s)}}   \right)  \,.
\end{eqnarray*} 
Thus,
\begin{eqnarray*}
\left( \sup_{\tau  \in [0,\infty)} e^{\gamma t} \|S [\Phi](t)\| _{H_{(s)}}  \right)
&  \leq & 
\left( \sup_{\tau  \in [0,\infty)}e^{\gamma t} \|\psi_{id}\| _{H_{(s)}}\right)+\varepsilon_0 C_{M }C_{\alpha, \gamma, H}\left( \sup_{\tau  \in [0,\infty)} e^{\gamma  \tau  }     \| \Phi (\cdot ,\tau  ) \| _{H_{(s)}}   \right) \\
&  &
+C_{M } \left( \sup_{\tau  \in [0,\infty)} e^{\gamma  \tau  }     \| \Phi (\cdot ,\tau  ) \| _{H_{(s)}}   \right)^{\alpha  +1}  \,.
\end{eqnarray*}
By Theorem~\ref{T13.2},
\[   
e^{\gamma t} \|\psi_{id}\| _{H_{(s)}} \|  
   \leq 
C_{m,s} e^{\gamma t}  e^{(\Re M-\frac{3H}{2})t}     \left(\|\psi _0 \|_{H_{(s)}} +   \|\psi _1 \|_{H_{(s)}}\right) \quad \mbox{for all} \quad t >0\,.
\] 
It  follows $ \psi  \in X(R,s,\gamma  ) $, provided that $ R$, $\varepsilon_0  $, and $\varepsilon  $ are sufficiently small. 
We skip the remaining part of the proof since it is similar to  case (i).
\medskip

\ndt
(iii) If  $\Re M > 3H/2 $ and $\varepsilon_0$ is sufficiently small, then according to the estimate of Theorem~\ref{T13.2}, 
we have  
$ \Phi_{id} (x,t) \in X(R,s,\gamma )$ with $\gamma <\frac{1}{\alpha+1}( 3H/2 -\Re M )<0$ for some $R>0$. 
On the other hand,
\begin{eqnarray*}
 e^{\gamma t} \|S [\Phi](t)\| _{H_{(s)}} 
&  \leq & 
e^{\gamma t} \|\psi_{id}\| _{H_{(s)}}+\varepsilon_0 C_{M }C_{\alpha, \gamma, H}\left( \max_{\tau  \in [0,t]} e^{\gamma  \tau  }     \| \Phi (\cdot ,\tau  ) \| _{H_{(s)}}   \right) \\
&  &
+C_{M } \left( \max_{\tau  \in [0,t]} e^{\gamma  \tau  }     \| \Phi (\cdot ,\tau  ) \| _{H_{(s)}}   \right)^{\alpha  +1} e^{\gamma t}   e^{ (\Re M   -\frac{3}{2}H)t} \int_{ 0}^{t} e ^{(- \Re M +\frac{3H}{2}     - \gamma (\alpha +1))b}  \, db\\
&  \leq & 
e^{\gamma t} \|\psi_{id}\| _{H_{(s)}}+\varepsilon_0 C_{M }C_{\alpha, \gamma, H}\left( \max_{\tau  \in [0,t]} e^{\gamma  \tau  }     \| \Phi (\cdot ,\tau  ) \| _{H_{(s)}}   \right) \\
&  &
+C_{M } \left( \max_{\tau  \in [0,t]} e^{\gamma  \tau  }     \| \Phi (\cdot ,\tau  ) \| _{H_{(s)}}   \right)^{\alpha  +1} 
\frac{      e ^{ - \gamma \alpha  t}  -  e^{\gamma t+(\Re M   -\frac{3}{2}H)t} }{ \frac{3H}{2}- \Re M      - \gamma (\alpha +1) }\,.
\end{eqnarray*}
Next we define  
\[
 T_\varepsilon :=\inf \{ T\,:\,\max_{\tau  \in [0,T]}e^{\gamma  \tau }  \|  \psi (x ,\tau ) \|_{ H_{(s)}({\mathbb R}^n)  } \geq 2\varepsilon \}\,, \quad \varepsilon :=   \max_{\tau  \in [0,\infty)} e^{\gamma  \tau } \| \Phi_{id}(\cdot ,\tau ) \|_{H_{(s)}({\mathbb R}^n)  }  \,.
 \] 
Then  
\[ 
2 \varepsilon  
   \leq   
\varepsilon +\varepsilon_0 2\varepsilon
 +  C_M  
 \varepsilon ^{\alpha  +1}    
\frac{e ^{ - \gamma \alpha T_\varepsilon }}{ \frac{3H}{2} -  \Re M  -\gamma (\alpha +1))}  
\]
implies \,
$
 T_\varepsilon \geq  -\frac{1}{\gamma}\ln \left(    \varepsilon  \right)- C(\alpha,\gamma,\varepsilon _0,H, M )
$. The global existence in Theorem~\ref{TIE} is proved.

\section{Proof of Theorem~\ref{TIE}: Decay of time derivative of solution }    
\subsection{Estimate of derivative of solution to linear equation. No source term}

\setcounter{equation}{0}
\renewcommand{\theequation}{\thesection.\arabic{equation}}

According to (\ref{T_EC}), for the equation (\ref{Lin5.1})   the energy $E(t)$ 
is conserved, that is, for all times of existence of the solution,   
$
\frac{d}{dt}E(t) =0
$. 
We state a global in time    ``energy estimate''  as folows.
\begin{theorem} 
\label{T14.5}
Consider the Cauchy problem
\begin{eqnarray}
\label{14.25}
\hspace{-1.6cm} &  &
  \psi _{tt}  - e^{-2H t} {\mathcal A}(x,\partial_x)  \psi  +   3H \psi  _t   +  \frac{m^2c^4}{h^2}  \psi +V(r)\psi  =  0 ,
   \\
\label{14.26}
\hspace{-1.6cm} &  &
\psi (x,0)=  \psi _0 (x)  , \,\,    \psi _t(x,0)= \psi _1 (x) ,\,\, \mbox{\rm supp}\, \psi _0, \mbox{\rm supp} \, \psi _1\subset \left\{ x \in \R^3\,\left|\, |x|> R_{ID}> {c}/{H}+R_{Sch} \right\},\right.
 \end{eqnarray} 
where $ {\mathcal A}(x,\partial_x)  $ is defined  in (\ref{OpA}),  $m^2   \in {\mathbb R} $,   and the potential $V$ is real-valued and bounded, $V(r) \in {\mathcal B}^\infty (\R^3)$.  
Then  there is 
a number $C>0 $
such that
\begin{equation}
\label{10.12} 
   \|  \psi _t (t) \|_{L^2(\R^3)}  + e^{-Ht}  \| \psi  (t)\|_{H_{(1)}}  
   \leq    
 C\left(    \|  \psi (t)\|_{L^2(\R^3)}  +   e^{-\frac{3}{2}Ht} \| \psi  _1\|_{L^2(\R^3)}  +     e^{-\frac{3}{2}Ht} \|  \psi _0 \|_{H_{(1)}}     \right) \,\, \mbox{ for all}\,\, t>0.  
\end{equation}
\end{theorem}
\medskip

\noindent
{\bf Proof.}
After application of the  Liouville transform  
$
 \psi =e^{-\frac{3}{2}Ht}\sqrt{F(r)}u$,  
we  arrive at the problem  
\[ 
u_{tt} - e^{-2Ht}{\mathcal A}_{3/2}(x,\partial_x) u -M^2 u +V(r)u= 0 ,\quad \quad u(x,0)= u_0 (x)\, , \quad u_t(x,0)=u_1 (x)\,, 
\]
with  smooth initial functions $ u_0 (x)$ and $ u_1 (x) $. Here $M^2 =  \frac{9 H^2}{4 } -\frac{m^2c^4}{h^2} $. The equation leads to the  identity
\begin{eqnarray*} 
&  &
\frac{1}{2} \frac{d}{dt} \left[ (u_{t},u_{t})_{L^2(\R^3)} - e^{-2Ht}  ( {\mathcal A}_{3/2}(x,\partial_x) u,  u )_{L^2(\R^3)}  -    M^2     (u,u )_{L^2(\R^3)} -  (Vu,u)_{L^2(\R^3)}\right] \\
&  &
\hspace{1cm} - He^{-2Ht}( {\mathcal A}_{3/2}(x,\partial_x) u,  u )_{L^2(\R^3)}=  0 \,.
\end{eqnarray*}
Here $ (u,u)_{L^2(\R^3)}$ denotes the scalar product in $L^2(\R^3) $. Since the operator ${\mathcal A}_{3/2}(x,\partial_x) $ is self-adjoint and non-positive, it follows 
\[
\frac{1}{2} \frac{d}{dt} \left[ (u_{t},u_{t})_{L^2(\R^3)} - e^{-2Ht}  ({\mathcal A}_{3/2}(x,\partial_x) u,  u )_{L^2(\R^3)}  -  \Re M^2     (u,u )_{L^2(\R^3)}-\Re (Vu,u)_{L^2(\R^3)} \right] \leq   0 \,.
\]
The integration in time gives 
\begin{eqnarray*} 
&  &
  (u_{t},u_{t})_{L^2(\R^3)} - e^{-2Ht}  ({\mathcal A}_{3/2}(x,\partial_x) u,  u )_{L^2(\R^3)}  -  \Re M^2     (u,u )_{L^2(\R^3)}  -\Re (Vu,u)_{L^2(\R^3)} \\
& \leq &
  (u_1,u_1)_{L^2(\R^3)} 
- e^{-2Ht}  ( {\mathcal A}_{3/2}(x,\partial_x) u_0 ,  u_0  )_{L^2(\R^3)}  -  \Re M^2     (u_0 ,u_0  )_{L^2(\R^3)} -\Re (Vu_0,u_0)_{L^2(\R^3)} ,
\end{eqnarray*}
and, consequently,
\[
 \|u_{t}(t)\|_{L^2(\R^3)}  + e^{-Ht}  \| u  (t)\|_{H_{(1)}} \leq  C_s  (\|u(t)\|_{L^2(\R^3)} + \|u_1\|_{L^2(\R^3)} +  \|u_0 \|_{H_{(1)}} )  \,.
\]
Hence, for the function $ \psi $, we have 
\[
 e^{-Ht}  \|\psi (t)\|_{H_{(1)}}  \leq  C_s ( \|\psi(t)\|_{L^2(\R^3)} + e^{-\frac{3}{2}Ht}\|u_{1} \|_{L^2(\R^3)} +  e^{-\frac{3}{2}Ht}\|u_0 \|_{H_{(1)}} )   
\]
and
\[
    \left\|e^{\frac{3}{2}Ht}\frac{3H}{\sqrt{F(r)}} \psi  (t) + e^{\frac{3}{2}Ht}\frac{2}{\sqrt{F(r)}}\psi _t (t) \right\|_{L^2(\R^3)}
  \leq   
 C_s\left(    e^{\frac{3}{2}Ht}\|  \psi (t)\|_{L^2(\R^3)}  +   \| \psi _1\|_{L^2(\R^3)}  +     \|  \psi _0 \|_{H_{(1)}}     \right)   .
\]
Then
\[
    \left\| \frac{ 2}{\sqrt{F(r)}}\psi _t (t) \right\|_{L^2(\R^3)}    
  \leq  
   \left\| \frac{3H}{\sqrt{F(r)}} \psi  (t)   \right\|_{L^2(\R^3)}    + e^{-\frac{3}{2}Ht}C_s\left(    e^{\frac{3}{2}Ht}\|  \psi (t)\|_{L^2(\R^3)}  +   \| \psi _1\|_{L^2(\R^3)}  +     \|  \psi _0 \|_{H_{(1)}}     \right) .
\]
Thus, the  theorem is proved. \hfill \qed

\begin{theorem}
\label{T11.2}
For   $ s \in {\mathbb N}\cup\{0\}$ and $   V =0 $,  the solution    $ \psi = \psi (x,t)$ of the Cauchy problem (\ref{14.25}) \&(\ref{14.26})
for    $\Re M \in (0, H/2)$  
satisfies the following estimate: 
\[
    \|  \psi _t (t) \|_{H_{(s)}}   
  \leq  
 C_s    e^{-Ht}\left(   \| \psi _1\|_{H_{(s+1)}}  +      \|  \psi _0 \|_{H_{(s+1)}}     \right) \quad \mbox{for all}\quad t \geq 0 \,.
\]
If   and  $\Re M >   {H}/{2}$ or $M=H/2 $, then 
\[
    \|  \psi _t (t) \|_{H_{(s)}}   
  \leq  
 C_s      e^{(\Re M-\frac{3}{2})Ht}      \left( \|\psi _0 \|_{H_{(s+1)}} +  \|\psi _1 \|_{H_{(s+1)}}\right) \quad \mbox{for all}\quad t \geq 0  \,.
\]
\end{theorem}
\medskip

\noindent
{\bf Proof.} 
According to Theorem~\ref{T13.2},   if $\Re M \in (0, H/2)$, then
\[
\|  \psi  (t) \|_{H_{(s)}}   
  \leq  
C_s   e^{   -H t}\Big(   \| \psi _0   \|_{H_{(s)}} 
+ \|\psi _1  \|_{H_{(s)}}  
 \Big) \,.
\]
Hence,  the inequality (\ref{10.12}) 
of Theorem~\ref{T14.5} implies
\begin{eqnarray*}
    \|  \psi _t (t) \|_{L^2(\R^3)}   
& \lesssim  &
  \|  \psi (t)   \|_{L^2(\R^3)} +   e^{-\frac{3}{2}Ht}\left(     \|  \psi _1\|_{L^2(\R^3)}  +     \|   \psi _0 \|_{H_{(1)}}     \right) \\   
& \lesssim  &
   e^{   -H t}\left(   \|\psi _0   \|_{H_{(1)}} 
+  \|\psi _1  \|_{L^2(\R^3)}  \right)   \,.
\end{eqnarray*}
For the case of $s >0$, we use induction. Indeed, if $  \partial_x $ is a first-order differential  operator, then the function $w=  \partial_x \psi$ solves equation 
\[ 
  w _{tt}  - e^{-2H t} {\mathcal A}(x,\partial_x)  w  +   3H w  _t   +  \frac{m^2c^4}{h^2}  w  =  e^{-2H t} [  \partial_x ,{\mathcal A}(x,\partial_x)]\psi \,,  
 \] 
where the commutator $ [  \partial_x ,{\mathcal A}(x,\partial_x)]$ is the second-order operator. We write $w= \tilde w+\tilde{\tilde w} $, where
\begin{eqnarray*}
\hspace{-1.6cm} &  &
  \tilde w _{tt}  - e^{-2H t} {\mathcal A}(x,\partial_x) \tilde w +   3H \tilde w _t   +  \frac{m^2c^4}{h^2} \tilde w   =  0 , \quad
\tilde w (x,0)=  \partial_x\psi _0 (x)  , \,\,    \tilde w _t(x,0)= \partial_x\psi _1 (x)\,,\\
\hspace{-1.6cm} &  &
  \tilde{\tilde w} _{tt}  - e^{-2H t} {\mathcal A}(x,\partial_x) \tilde{\tilde w} +   3H \tilde{\tilde w} _t   +  \frac{m^2c^4}{h^2} \tilde{\tilde w}   =   e^{-2H t} [\partial_x ,{\mathcal A}(x,\partial_x)]\psi , \quad
\tilde{\tilde w} (x,0)=  0  , \,\,    \tilde{\tilde w} _t(x,0)= 0 .
 \end{eqnarray*} 
An application of Theorem~\ref{T14.5} and Theorem~\ref{T13.2}  leads to
\begin{eqnarray*}
   \|  \tilde w _t (t) \|_{L^2(\R^3)}  + e^{-Ht}  \| \tilde w(t)\|_{H_{(1)}}  
  & \lesssim   &  
 C\left(    \|  \tilde w (t)\|_{L^2(\R^3)}  +   e^{-\frac{3}{2}Ht} \|   w _1\|_{L^2(\R^3)}  +     e^{-\frac{3}{2}Ht} \|   w_0 \|_{H_{(1)}}     \right)  \\
  & \lesssim   &  
 C\left(    e^{-Ht}  (\|   \psi _1\|_{H_{(1)}} +     \|   \psi_0 \|_{H_{(1)}}    +   e^{-\frac{3}{2}Ht} \|   \psi _1\|_{H_{(1)}}  +     e^{-\frac{3}{2}Ht} \|   \psi_0 \|_{H_{(2)}}     \right) \\
  & \lesssim  &  
 C    e^{-Ht}  \left(\|   \psi _1\|_{H_{(1)}} +     \|   \psi_0 \|_{H_{(2)}}  \right)   
\end{eqnarray*} 
while Theorem~\ref{TdecG} and Theorem~\ref{T13.2} lead to 
 \begin{eqnarray*} 
 \| \tilde{ \tilde w}_t(t)\|_{L^2(\R^3)}   
  & \lesssim   &
e ^{- Ht}     \int_{ 0}^{t} e ^{ Hb} \| e^{-2H b} [ \partial_x ,{\mathcal A}(x,\partial_x)] \psi (x, b) \| _{L^2(\R^3)}   db  \\
  & \lesssim   &
e ^{- Ht}     \int_{ 0}^{t} e ^{ -Hb} \| \psi (x, b) \| _{H_{(2)}}   db  \\
  & \lesssim   &
e ^{- Ht}      \left(  \| \psi_0   \| _{H_{(2)}} + \| \psi_1  \| _{H_{(2)}}\right)    \,,
\end{eqnarray*}
respectively. Hence
\[ 
 \| \psi_t(t)\|_{H_{(1)}}   
 \lesssim   
e ^{- Ht}      \left(  \| \psi_0   \| _{H_{(2)}} + \| \psi_1  \| _{H_{(2)}}\right)    \,.
\]
The induction completes the proof of the case of $ s \in {\mathbb N}  $.

Next, we  can consider the case of  $\Re M >   H/2$. According to Theorems~\ref{T13.2},  if $\Re M >  H/2$ or $ M = H/2$, then 
\[
\|  \psi  (t) \|_{H_{(s)}}   
  \leq   
C_{m,s}   e^{(\Re M-\frac{3}{2})Ht}     \left( \|\psi _0 \|_{H_{(s)}} +  \|\psi _1 \|_{H_{(s)}}\right)  
\]
while also holds (\ref{10.12})  
by Theorem~\ref{T14.5}. It follows
\begin{eqnarray*}
   \|  \psi _t (t) \|_{L^2(\R^3)} 
 & \leq &  
 C\left(   C_{m,s}   e^{(\Re M-\frac{3}{2}H)t}     \left( \|\psi _0 \|_{L^2(\R^3)} +  \|\psi _1 \|_{L^2(\R^3)} \right)   +   e^{-\frac{3}{2}Ht} \| \psi  _1\|_{L^2(\R^3)}  +     e^{-\frac{3}{2}Ht} \|  \psi _0 \|_{H_{(1)}}     \right)  \\
 & \leq &  
     C_{m,s}   e^{(\Re M-\frac{3}{2}H)t}     \left( \|\psi _0 \|_{H_{(1)}} +  \|\psi _1 \|_{L^2(\R^3)}\right)       \,\, \mbox{ for all}\,\, t>0. 
\end{eqnarray*}
The remaining part of the proof with $s >0$  is similar to the previous case. 
The  theorem is proved. \hfill $\square$

\subsection*{Estimate of derivative of solution to linear equation. Vanishing initial functions} 

\begin{theorem}
\label{TdecG}
The operator $G$ has the following     property :
\begin{eqnarray}
\label{11.18} 
\hspace{-0.3cm} \mbox{\rm (i)} \quad \| \partial_t G[f](t,x)\|_{H_{(s)} }  
  & \lesssim   &
e ^{- Ht}     \int_{ 0}^{t} e ^{ Hb} \| f(x, b) \| _{H_{(s)}}   db  ,
  \quad \mbox{\rm if} \,\, 0< \Re M<\frac{H}{2} \,\, or \,\, M=\frac{H}{2}\,,\\
\label{11.19} 
\hspace{-0.3cm}\mbox{\rm (ii)}\quad \| \partial_t G[f](t,x)\|_{H_{(s)} } 
  & \lesssim   &
 e ^{(\Re M-\frac{1}{2}H)t}     \int_{ 0}^{t} e ^{-(\Re M-\frac{1}{2}H)b} \| f(x, b) \| _{H_{(s)}} 
   db       ,
   \,\,  \mbox{\rm if} \quad\frac{H}{2}<\Re M< \frac{3}{2}H \,,
\end{eqnarray}
for all $t \geq 0$, where supp\,$f \subset \{(x,t) \in \R^3\times[0,\infty)\,|\,|x|> R_{ID}- c(1-e^{-tH})/H \,\}$. 
\end{theorem}
\medskip

\ndt
{\bf Proof.} 
 In   {the case of $M=H/2$ }  one has 
$
E(r,t;0,b;H/2)
  :=  \frac{1}{2} e^{\frac{1}{2} H (b+t)} 
$. 
For $\psi=G[f] $, the representation 
\[ 
\psi (x,t) 
   =  
e ^{-Ht}     \int_{ 0}^{t}    e^{2 H b} db
  \int_{ 0}^{\phi (t)- \phi (b)} v_f(x,r ;b) \, dr
  \] 
  implies
\[
\partial_t \psi (x,t) 
  =   
 -H \psi (x,t) +
    e ^{-2Ht}      \int_{ 0}^{t}    e^{2 H b}   
   v_f(x,\phi (t)- \phi (b) ;b)    db\,.  
  \]
Consequently, 
\begin{eqnarray*} 
\|\partial_t \psi (x,t) \| _{H_{(s)} }
&  \lesssim  &
 \| \psi (x,t)\| _{H_{(s)} } +
  e ^{-2Ht}      \int_{ 0}^{t}    e^{2 H b}   
   \| v_f(x,\phi (t)- \phi (b) ;b) \| _{H_{(s)} }     db\\
&  \lesssim  &
 \| \psi (x,t)\| _{H_{(s)} } +
  e ^{-2Ht}      \int_{ 0}^{t}    e^{2 H b}   
   \|  f(x,  b) \| _{H_{(s)} }     db\\
&  \lesssim  &
  e ^{- Ht}      \int_{ 0}^{t}    e^{  H b}   
   \|  f(x,  b) \| _{H_{(s)} }     db +
  e ^{-2Ht}      \int_{ 0}^{t}    e^{2 H b}   
   \|  f(x,  b) \| _{H_{(s)} }     db\\
&  \lesssim  &
  e ^{- Ht}      \int_{ 0}^{t}    e^{  H b}   
   \|  f(x,  b) \| _{H_{(s)} }     db \,. 
  \end{eqnarray*}
Hence, (\ref{11.18}) with $M=H/2$ is proved. 

If { $\Re M\not= H/2$}, then, in order to estimate the time  derivative  of the function $ \psi $ if $\Re M\not= H/2$,  we write
\[ 
\partial_t \psi (x,t) 
   =   
A_1+A_2+A_3  ,
\]
where, with $E(\phi (t)- \phi (b),t; 0,b;M) 
  = 
\frac{1}{2}   e^{ \frac{1}{2}H (b+t)}$, we denoted 
\begin{eqnarray*} 
A_1 & := & 
-\frac{3}{2}H \psi (x,t),\quad 
A_2   :=  
 e ^{-2Ht}    \int_{ 0}^{t}    e^{2H b} v_f(x,\phi (t)- \phi (b) ;b) \, db,\\
A_3 & := &e ^{-\frac{3}{2}Ht}  2   \int_{ 0}^{t} db
  \int_{ 0}^{\phi (t)- \phi (b)} e ^{\frac{3H}{2}b}   v_f(x,r ;b)  \partial_t E(r,t;0,b;M)   \, dr \,.
\end{eqnarray*}
Due to  (\ref{Blplq})   we have  
$
\|  v(x,r ;b) \| _{H_{(s)} } \leq   C \|f(x,b)\|_{H_{(s)} }$  for all $r \in(0, \phi (t)- \phi (b)] 
$. 
 Hence,  for $A_2$, we obtain 
\begin{eqnarray} 
\label{Es_A_2} 
\|A_2  \| _{H_{(s)} } 
  & \lesssim &  
  e ^{-2Ht}    \int_{ 0}^{t}    e^{2H b}\|f(x,b)\|_{H_{(s)} }\, db\,.  
\end{eqnarray}
For the   term $A_3$ of the derivative $\partial_t \psi  $,   we have
\[
 \|A_3\| _{H_{(s)}}
  \lesssim  
   e ^{-\frac{3}{2}Ht}   \int_{ 0}^{t} db
  \int_{ 0}^{\phi (t)- \phi (b)} e ^{\frac{3H}{2}b}  \| v_f(x,r ;b) \| _{H_{(s)}} | \partial_t E(r,t;0,b;M) |  \, dr  \,,
\]
that is,
\begin{equation}
\label{14.29}  
 \|A_3\| _{H_{(s)}} 
   \lesssim   
 e ^{-\frac{3}{2}Ht}     \int_{ 0}^{t} e ^{\frac{3H}{2}b} \| f(x, b) \| _{H_{(s)}} db
  \int_{ 0}^{\phi (t)- \phi (b)}  | \partial_t E(r,t;0,b;M) |  \, dr \,.
\end{equation}
We apply the following estimate for the time derivative of the kernel $E(r,t; 0,b;M) $.
\begin{proposition}
\label{P6.1}
If $\Re M>0 $, then 
\[
\int_{ 0}^{ (e^{-Hb}- e^{-Ht})/H} | \partial_t   E(r,t; 0,b;M) |\,dr  
  \lesssim   
 \cases{  e^{ \frac{1}{2}H(t-b)} ,   \quad  \mbox{\rm for}\quad \Re M<H/2, \cr
e^{(\Re M +H)(t- b)} ,   \quad  \mbox{\rm for}\quad \Re M>H/2, }   
\] 
for all $t\geq 0 $ and $b\geq 0 $ such that $b<t $.
\end{proposition}
\medskip

\noindent
{\bf Proof.} 
We have from (\ref{0.6}) the expression 
\begin{eqnarray} 
\label{14.30}
&  &
\partial_t E(r,t; 0,b;M) 
   = 
I_1(b,t,r)+I_2(b,t,r)\,,
\end{eqnarray}
where
\begin{eqnarray*}
I_1(b,t,r)
& := &  
\left( \partial_t 4^{-\frac{M}{H}} e^{ M  (  b+ t)} \left(\left(e^{-H b}+e^{-H t}\right)^2-(H r)^2\right)^{\frac{M}{H}-\frac{1}{2}} \right) \nonumber \\
&  &
\times  F \left(\frac{1}{2}-\frac{M}{H},\frac{1}{2}-\frac{M}{H};1;\frac{\left(-e^{-H t}+e^{-H b}\right)^2-(r H)^2}{\left(e^{-H t}+e^{-Hb}\right)^2-(r H)^2}\right)\,, \nonumber  \\
I_2(b,t,r)
& := &  
4^{-\frac{M}{H}} e^{ M  (  b+ t)} \left(\left(e^{-H b}+e^{-H t}\right)^2-(H r)^2\right)^{\frac{M}{H}-\frac{1}{2}}  \nonumber \\
&  &
\times \partial_t  F \left(\frac{1}{2}-\frac{M}{H},\frac{1}{2}-\frac{M}{H};1;\frac{\left(-e^{-H t}+e^{-H b}\right)^2-(r H)^2}{\left(e^{-H t}+e^{-Hb}\right)^2-(r H)^2}\right) \,. \nonumber  
\end{eqnarray*}
For $I_1$ we have
\begin{eqnarray*} 
I_1(b,t,r) 
& = &  
-4^{-\frac{M}{H}} e^{M (b+t)} \left(\left(e^{-b H}+e^{-H t}\right)^2-H^2 r^2\right)^{\frac{M}{H}-\frac{3}{2}} \nonumber \\
&  &
\times\left(-M e^{-2 b H}-H e^{-b H-H t}+H^2 M r^2+M e^{-2 H t}-H e^{-2 H t}\right)    \nonumber \\
&  &
\times  F \left(\frac{1}{2}-\frac{M}{H},\frac{1}{2}-\frac{M}{H};1;\frac{\left(-e^{-H t}+e^{-H b}\right)^2-(r H)^2}{\left(e^{-H t}+e^{-Hb}\right)^2-(r H)^2}\right) \,.
\end{eqnarray*}
Since $\Re M>0 $ and $b\leq t$, we obtain
\[
|I_1(b,t,r)| 
  \lesssim  
  e^{\Re M (b+t)} e^{-2 b H}   \left(\left(e^{-b H}+e^{-H t}\right)^2-H^2 r^2\right)^{\frac{\Re M}{H}-\frac{3}{2}} \,. 
\]
To estimate the integral 
\[
\int_0^{\frac{1}{H}(e^{-Hb}-e^{-H t})}|I_1(b,t,r)| dr \\
  \lesssim   
e^{\Re M (b+t)} e^{-2 b H} \int_0^{\frac{1}{H}(e^{-Hb}-e^{-H t})}  \left(\left(e^{-b H}+e^{-H t}\right)^2-H^2 r^2\right)^{\frac{\Re M}{H}-\frac{3}{2}}   dr  
\]
 we set     $r=e^{-Ht}y H^{-1}$  and $z :=e^{H(t-b)} \in [1,\infty)$  in the last integral   
\begin{eqnarray*}  
&  &
\int_0^{\frac{1}{H}(e^{-Hb}-e^{-H t})}\left(\left(e^{-b H}+e^{-H t}\right)^2-H^2 r^2\right)^{\frac{\Re M}{H}-\frac{3}{2}} dr \\
& = &
e^{-Ht} H^{-1}e^{-( 2\Re M -3H)t} (z-1) (z+1)^{2 ({\frac{\Re M}{H}-\frac{3}{2}})} F\left(\frac{1}{2},\frac{3}{2}-\frac{\Re M}{H};\frac{3}{2};\frac{(z-1)^2}{(z+1)^2}\right)\,.
\end{eqnarray*} 
 Hence,     
\begin{eqnarray*}  
&  &
\int_0^{\frac{1}{H}(e^{-Hb}-e^{-H t})}\left(\left(e^{-b H}+e^{-H t}\right)^2-H^2 r^2\right)^{\frac{\Re M}{H}-\frac{3}{2}} dr \\
& \lesssim &
 e^{- 2(\Re M - H)t}  (e^{H(t-b)}-1) (e^{H(t-b)}+1)^{({2\frac{\Re M}{H}-3})} F\left(\frac{1}{2},\frac{3}{2}- \frac{\Re M}{H} ;\frac{3}{2};\frac{(z-1)^2}{(z+1)^2}\right)\,.
\end{eqnarray*}
We estimate the last hypergeometric function in the next lemma.
\begin{lemma}
The following is true
\begin{equation}
\label{14.31}  
\int_0^{\frac{1}{H}(e^{-Hb }-e^{-H t})}\left(\left(e^{-b H}+e^{-H t}\right)^2-H^2 r^2\right)^{\frac{\Re M}{H}-\frac{3}{2}} dr  
\lesssim 
 \cases{     e^{-t ( \Re M  -\frac{1}{2}H)}     e^{  -b ({ \Re M -  \frac{3}{2}H })}  \,\, \mbox{  if}\,\,\, \Re M<H/2,\cr
e^{- 2(\Re M - H)b}\,\, \mbox{  if}\,\,\, \Re M>H/2,}
\end{equation}
for all $t\geq 0 $ and $b\geq 0 $ such that $b<t $.
\end{lemma} 
\medskip

\noindent
{\bf Proof.}
For the case of $\Re M>H/2 $, we write    
\[  
\int_0^{\frac{1}{H}(e^{-Hb }-e^{-H t})}\left(\left(e^{-b H}+e^{-H t}\right)^2-H^2 r^2\right)^{\frac{\Re M}{H}-\frac{3}{2}} dr 
  \lesssim  
e^{- 2(\Re M - H)b}\,.
\] 
This proves the second case of (\ref{14.31}).  
For  the case of $\Re M<H/2 $,  we write    
\begin{eqnarray*}  
&  &
\int_0^{\frac{1}{H}(e^{-Hb}-e^{-H t})}\left(\left(e^{-b H}+e^{-H t}\right)^2-H^2 r^2\right)^{\frac{\Re M}{H}-\frac{3}{2}} dr \\
& \lesssim &
 e^{- 2(\Re M - H)t}  (e^{H(t-b)}-1) (e^{H(t-b)}+1)^{({2\frac{\Re M}{H}-3})} F\left(\frac{1}{2},\frac{3}{2}- \frac{\Re M}{H} ;\frac{3}{2};\frac{(z-1)^2}{(z+1)^2}\right)\\
& \lesssim &
  e^{- 2(\Re M - H)t} (e^{H(t-b)}-1) (e^{H(t-b)}+1)^{({2\frac{\Re M}{H}-3})}  \left(1-\frac{(z-1)^2}{(z+1)^2} \right)^{\frac{\Re M}{H} -\frac{1}{2}} F\left(1,  \frac{\Re M}{H} ;\frac{3}{2};\frac{(z-1)^2}{(z+1)^2}\right)\\
& \lesssim &
 e^{- ( \Re M  -\frac{1}{2}H)t}     e^{  - ({ \Re M -  \frac{3}{2}H })b}\,.   
\end{eqnarray*} 
Thus,   
(\ref{14.31}) and the lemma are proved. \qed \\

Thus,  for  $\Re M<H/2$, we obtain
 \[ 
\int_0^{\frac{1}{H}(e^{-Hb}-e^{-H t})}|I_1(b,t,r)| dr
  \lesssim   
 e^{\Re M (b+t)} e^{-2 b H}    e^{-t ( \Re M  -\frac{1}{2}H)}     e^{  -b ({ \Re M -  \frac{3}{2}H })} \quad  \mbox{\rm for}\quad \Re M<H/2\,,
\] 
that is,
 \[  
\int_0^{\frac{1}{H}(e^{-Hb}-e^{-H t})}|I_1(b,t,r)| dr 
  \lesssim   
 e^{ \frac{1}{2}H(t-b)}      ,   \quad  \mbox{\rm for}\quad \Re M<H/2\,.
\] 
For the case of $\Re M>H/2$, we have 
 \[ 
\int_0^{\frac{1}{H}(e^{-Hb}-e^{-H t})}|I_1(b,t,r)|\, dr 
  \lesssim   
 e^{\Re M (t- b)} ,   \quad  \mbox{\rm for}\quad \Re M>H/2\,.
\]
Finally, for $I_1$, we have obtained
 \[ 
\int_0^{\frac{1}{H}(e^{-Hb}-e^{-H t})}|I_1(b,t,r)| dr 
  \lesssim   
 \cases{  e^{ \frac{1}{2}H(t-b)} ,   \quad  \mbox{\rm for}\quad \Re M<H/2\,,\cr
e^{\Re M (t- b)} ,   \quad  \mbox{\rm for}\quad \Re M>H/2\,.}
\]
Next we consider the   term $I_2$ of  (\ref{14.30}). If $  \Re M> H/2$, then 
\begin{eqnarray*} 
|I_2(b,t,r)|
& = &
\Bigg| 4^{-\frac{M}{H}} e^{ M  (  b+ t)} \left(\left(e^{-H b}+e^{-H t}\right)^2-(H r)^2\right)^{\frac{M}{H}-\frac{1}{2}}  \nonumber \\
&  &
\times 
\Bigg[ -\frac{(H-2 M)^2 e^{H (b+t)} \left(H^2 r^2 e^{2 H (b+t)}+e^{2 b H}-e^{2 H t}\right)  }{H \left(H^2 r^2 \left(-e^{2 H (b+t)}\right)+2 e^{H (b+t)}+e^{2 b H}+e^{2 H t}\right)^2} \\
&  &
\times F \left(\frac{3}{2}-\frac{M}{H},\frac{3}{2}-\frac{M}{H};2;\frac{\left(e^{-b H}-e^{-H t}\right)^2-H^2 r^2}{\left(e^{-b H}+e^{-H t}\right)^2-H^2 r^2}\right)\Bigg]\Bigg| \\
 &  \lesssim &
\Bigg|  e^{ M  (  b+ t)} \left(\left(e^{-H b}+e^{-H t}\right)^2-(H r)^2\right)^{\frac{M}{H}-\frac{1}{2}} 
\Bigg[ \frac{e^{H (b+t)} \left(H^2 r^2 e^{2 H (b+t)}+e^{2 b H}-e^{2 H t}\right)  }{ e^{  4H (b+t)} \left[     e^{  -H b}  e^{ - H t} \right]^2} \Bigg]\Bigg| \,.
\end{eqnarray*}
Thus, 
\begin{eqnarray*} 
|I_2(b,t,r)|
&  \lesssim &
  e^{  \Re M  (  b+ t)}e^{ H (b+t)} \left(\left(e^{-H b}+e^{-H t}\right)^2-(H r)^2\right)^{\frac{\Re M}{H}-\frac{1}{2}}  \nonumber \\
&  &
\times 
 \left|H^2 r^2 -e^{-2 b H}+e^{-2 H t}\right|   \quad  \mbox{\rm for}\quad \Re M>H/2\,.
\end{eqnarray*}
On the other hand, 
we denote $ B:=e^{-  b H}$ and $T:=e^{-  t H} $, then
\begin{eqnarray*} 
&  &
\int_0^{\phi (t)- \phi (b)}\left(\left(e^{-H b}+e^{-H t}\right)^2-(H r)^2\right)^{\frac{\Re M}{H}-\frac{1}{2}}  
\left|H^2 r^2 -e^{-2 b H}+e^{-2 H t}\right|dr \\
& = &
\frac{(B-T)^2 (B+T)^{\frac{2 \Re M}{H}-1}}{3 H}
\Bigg\{3 (B+T)  F \left(\frac{1}{2},\frac{1}{2}-\frac{\Re M}{H};\frac{3}{2};\frac{(B-T)^2}{(B+T)^2}\right)\\
&  &
+(T-B)  F \left(\frac{3}{2},\frac{1}{2}-\frac{\Re M}{H};\frac{5}{2};\frac{(B-T)^2}{(B+T)^2}\right)\Bigg\}\,.
\end{eqnarray*}
Since  $\frac{3}{2} -\frac{1}{2}-\left(\frac{1}{2}-\frac{\Re M}{H}\right)>0$, we have
 $ F \left(\frac{1}{2},\frac{1}{2}-\frac{\Re M}{H};\frac{3}{2};\frac{(B-T)^2}{(B+T)^2}\right) \lesssim  1 $
 and $F \left(\frac{3}{2},\frac{1}{2}-\frac{\Re M}{H};\frac{5}{2};\frac{(B-T)^2}{(B+T)^2}\right) \lesssim  1$ and 
\[ 
\int_0^{\phi (t)- \phi (b)}\left(\left(e^{-H b}+e^{-H t}\right)^2-(H r)^2\right)^{\frac{\Re M}{H}-\frac{1}{2}}  
\left|H^2 r^2 -e^{-2 b H}+e^{-2 H t}\right|dr   
 \lesssim  
 (B-T)^2 B^{\frac{2 \Re M}{H} }\,.   
\]
 Hence, we obtain the estimate   
 \begin{eqnarray*} 
    \int_{ 0}^{\phi (t)- \phi (b)}|I_2(b,t,r)| \, dr 
    & \lesssim & 
e^{  \Re M  (  b+ t)}e^{ H (b+t)}  \int_{ 0}^{\phi (t)- \phi (b)}   \left(\left(e^{-H b}+e^{-H t}\right)^2-(H r)^2\right)^{\frac{\Re M}{H}-\frac{1}{2}}  \nonumber \\
&  &
\times 
 \left|H^2 r^2 -e^{-2 b H}+e^{-2 H t}\right|  \, dr\\
    & \lesssim & 
e^{  (\Re M +H) (t-  b )} 
  \quad  \mbox{\rm for}\quad \Re M>H/2\,.
\end{eqnarray*} 
It follows
\begin{eqnarray*}  
\int_0^{\frac{1}{H}(e^{-Hb}-e^{-H t})} |\partial_t E(r,t; 0,b;M)| \, d r 
 & \leq &
\int_0^{\frac{1}{H}(e^{-Hb}-e^{-H t})} |I_1(b,t,r)|\, d r+\int_0^{\frac{1}{H}(e^{-Hb}-e^{-H t})}|I_2(b,t,r)|\, d r\\
 &  \lesssim  &
 \cases{  e^{ \frac{1}{2}H(t-b)} ,   \quad  \mbox{\rm for}\quad \Re M<H/2 \,,\cr
e^{(\Re M +H) (t- b)}   ,   \quad  \mbox{\rm for}\quad \Re M>H/2 \,.}  
\end{eqnarray*}
 For  $\Re M> {H}/{2}$,  we derive
\[
\| A_3 \| _{H_{(s)}}  
 \lesssim    
   e^{(\Re M -\frac{1}{2}H) t}   \int_{ 0}^{t}   e^{ -(\Re M-\frac{1}{2}H)b  }  \| f(x, b) \| _{H_{(s)}}    db\quad \mbox{\rm if} \quad \Re M>\frac{H}{2}\,.
\]
In the case of  $\Re M<H/2 $, we use \cite[(23)~Sec 2.1.4]{B-E} and obtain 
\begin{eqnarray*} 
&  &
\Bigg|  e^{ M  (  b+ t)} \left(\left(e^{-H b}+e^{-H t}\right)^2-(H r)^2\right)^{\frac{M}{H}-\frac{1}{2}}  \partial_t  F \left(\frac{1}{2}-\frac{M}{H},\frac{1}{2}-\frac{M}{H};1;\frac{\left(e^{-H b}-e^{-H t}\right)^2-(r H)^2}{\left(e^{-H b}+e^{-Ht}\right)^2-(r H)^2}\right)  \Bigg|\\
& = &
\Bigg|  (H-2 M)^2 e^{(b+t) (H+M)} \left(H^2 r^2 e^{2 H (b+t)}+e^{2 b H}-e^{2 H t}\right) \left(\left(e^{-b H}+e^{-H t}\right)^2-H^2 r^2\right)^{\frac{M}{H}-\frac{5}{2}}  \nonumber \\
&  & 
\times H^{-1} e^{-4 H (b+t)}F \left(\frac{3}{2}-\frac{M}{H},\frac{3}{2}-\frac{M}{H};2;\frac{\left(e^{-b H}-e^{-H t}\right)^2-H^2 r^2}{\left(e^{-b H}+e^{-H t}\right)^2-H^2 r^2}\right)\Bigg|\,.
\end{eqnarray*}
On the other hand,
\begin{eqnarray*}
&  &
F \left(\frac{3}{2}-\frac{M}{H},\frac{3}{2}-\frac{M}{H};2;\frac{\left(e^{-b H}-e^{-H t}\right)^2-H^2 r^2}{\left(e^{-b H}+e^{-H t}\right)^2-H^2 r^2}\right)\\
& = &
4^{\frac{2 M}{H}-1}   e^{- (t+b)( 2 M - H)}  \left( \left(e^{-b H}+e^{-H t}\right)^2-H^2 r^2 \right)^{1-\frac{2 M}{H} }\\
&  &
\times F \left(\frac{M}{H}+\frac{1}{2},\frac{M}{H}+\frac{1}{2};2;\frac{\left(e^{-b H}-e^{-H t}\right)^2-H^2 r^2}{\left(e^{-b H}+e^{-H t}\right)^2-H^2 r^2}\right)\,.
\end{eqnarray*}
Consequently,
\begin{eqnarray*} 
&  &
\Bigg|  e^{ M  (  b+ t)} \left(\left(e^{-H b}+e^{-H t}\right)^2-(H r)^2\right)^{\frac{M}{H}-\frac{1}{2}} \partial_t  F \left(\frac{1}{2}-\frac{M}{H},\frac{1}{2}-\frac{M}{H};1;\frac{\left(-e^{-H t}+e^{-H b}\right)^2-(r H)^2}{\left(e^{-H t}+e^{-Hb}\right)^2-(r H)^2}\right)  \Bigg|\\
& \lesssim &
 \left|  e^{-M(b+t) }  \left(H^2 r^2 -e^{-2 b H}+e^{-2 H t}\right) \left(\left(e^{-b H}+e^{-H t}\right)^2-H^2 r^2\right)^{-\frac{M}{H}-\frac{3}{2}} \right|  \,.
\end{eqnarray*} 
\begin{lemma} If $\Re M<H/2 $, then
\[
   \int_{ 0}^{\phi (t)- \phi (b)}\left|  e^{-M(b+t) }  \left(H^2 r^2 -e^{-2 b H}+e^{-2 H t}\right) \left(\left(e^{-b H}+e^{-H t}\right)^2-H^2 r^2\right)^{-\frac{M}{H}-\frac{3}{2}} \right|    \, dr \\
  \lesssim   
    e^{  \frac{1}{2}(t - b) H}   
\]
for all $ t\geq b \geq 0$.
\end{lemma}
\medskip

\noindent
{\bf Proof.} Indeed, if we denote $ B:=e^{-  b H}$ and $T:=e^{-  t H} $, then
\begin{eqnarray*} 
 &  &
   \int_{ 0}^{\phi (t)- \phi (b)}\left|  e^{-M(b+t) }  \left(H^2 r^2 -e^{-2 b H}+e^{-2 H t}\right) \left(\left(e^{-b H}+e^{-H t}\right)^2-H^2 r^2\right)^{-\frac{M}{H}-\frac{3}{2}} \right|    \, dr \\
 & \lesssim  &
  e^{-\Re M(b+t) }  e^{-2 b H}  \int_{ 0}^{\phi (t)- \phi (b)}   \left(\left(e^{-b H}+e^{-H t}\right)^2-H^2 r^2\right)^{-\frac{\Re M}{H}-\frac{3}{2}}     \, dr \\
    & \lesssim & 
 e^{-\Re M(b+t) }  e^{-2 b H}  \left((B-T) (B+T)^{-\frac{2 \Re M}{H}-3}\right) \left(1-\frac{(B-T)^2}{(B+T)^2}\right)^{-\frac{\Re M}{H}-\frac{1}{2}}    F \left(1,-\frac{\Re M}{H};\frac{3}{2};\frac{(B-T)^2}{(B+T)^2}\right)\\
    & \lesssim &
  e^{  \frac{1}{2}(t - b) H} \,.  
\end{eqnarray*}
The lemma is proved. 
\qed
\medskip

Finally
\[ 
  \int_{ 0}^{\phi (t)- \phi (b)} |I_2 (b,t,r)|\, d r 
  \lesssim  
  e^{  \frac{1}{2}(t - b) H} \quad \mbox{\rm if} \quad \Re M< {H}/{2}  \,. 
\] 
It follows
\[   
\int_0^{\frac{1}{H}(e^{-Hb}-e^{-H t})} |\partial_t E(r,t; 0,b;M)| \, d r  
  \lesssim  
 \cases{  e^{ \frac{1}{2}H(t-b)} ,   \quad  \mbox{\rm for}\quad \Re M<H/2 \,,\cr
e^{(\Re M +H) (t- b)}  ,   \quad  \mbox{\rm for}\quad \Re M>H/2 \,.}  
\]
Proposition \ref{P6.1} is proved \qed
\medskip

Now we estimate the norm of  $A_3$ . 
We use (\ref{14.29}) for $0\leq\Re M<\frac{H}{2}$ 
and obtain   
\[
\| A_3 \| _{H_{(s)}}  
   \lesssim    
 e ^{- Ht}     \int_{ 0}^{t} e ^{Hb} \| f(x, b) \| _{H_{(s)}} 
   db    \quad \mbox{\rm if} \quad 0\leq\Re M<\frac{H}{2}\,.
\]
For  $\Re M> {H}/{2}$,  we derive
\[
\| A_3 \| _{H_{(s)}}  
  \lesssim  
 e ^{(\Re M-\frac{1}{2}H)t}     \int_{ 0}^{t} e ^{-(\Re M-\frac{1}{2}H)b} \| f(x, b) \| _{H_{(s)}} 
   db  \quad \mbox{\rm if} \quad \Re M>\frac{H}{2}\,.
\]
If we collect estimates for $A_1$, $A_2$ (\ref{Es_A_2}), and $A_3$, in the case of $0<\Re M<\frac{H}{2} $, then
\begin{eqnarray*} 
\| \partial_t \psi(t,x)\|_{H_{(s)} }  
  & \lesssim   &
e^{-Ht}   
 \int_{ 0}^{t} e^{Hb} 
   \|f(x,b)\|_{H_{(s)} }\,db  +   e ^{-2Ht}    \int_{ 0}^{t}    e^{2H b}\|f(x,b)\|_{H_{(s)} }\, db \\
&  &
+   e ^{- Ht}     \int_{ 0}^{t} e ^{Hb} \| f(x, b) \| _{H_{(s)}}   db  ,
  \quad \mbox{\rm if} \quad\Re M<\frac{H}{2} .  
\end{eqnarray*} 
Thus, for  $\Re M<{H}/{2}$,  we obtain (\ref{11.18}).

For $\Re M>{H}/{2}$,  we have
\begin{eqnarray*} 
\| \partial_t \psi(t,x)\|_{H_{(s)} }  
  & \lesssim   &
e^{-Ht}   
 \int_{ 0}^{t} e^{Hb} 
   \|f(x,b)\|_{H_{(s)} }\,db  +   e ^{-2Ht}    \int_{ 0}^{t}    e^{2H b}\|f(x,b)\|_{H_{(s)} }\, db \\
&  &
+   e ^{(\Re M-\frac{1}{2}H)t}     \int_{ 0}^{t} e ^{-(\Re M-\frac{1}{2}H)b} \| f(x, b) \| _{H_{(s)}} 
   db    \\
  & \lesssim   &e ^{(\Re M-\frac{1}{2}H)t}     \int_{ 0}^{t} e ^{-(\Re M-\frac{1}{2}H)b} \| f(x, b) \| _{H_{(s)}} 
   db         ,
  \quad \mbox{\rm if} \quad\Re M>\frac{H}{2} \,.
\end{eqnarray*}
Hence, for  $\Re M>{H}/{2}$,  we obtain (\ref{11.19}). 
Thus, we have proved Theorem~\ref{TdecG}.\qed

\subsection{Estimate of time-derivative of solution to semilinear equation}

If the function $\psi=\psi(t,x) $ solves the equation (\ref{5.1}), 
then
\[
\partial_t\psi= \partial_t\psi_{id}+\partial_tG[V\psi] +\partial_tG[F\Psi(\psi)]\,.
\]
According to (\ref{3.2}) and Theorem~\ref{T11.2}, for  $\Re M< {H}/{2}$,  we have
\[
    \|  \partial_t\psi_{id}(t) \|_{H_{(s)}}   
  \leq   
 C    e^{-Ht}\left(   \| \psi _1\|_{H_{(s)}}  +      \|  \psi _0 \|_{H_{(s+1)}}     \right)  \,.
\] 
Further, according to (\ref{3.2}) and  (\ref{11.18}),  Theorem~\ref{TdecG}, with $\gamma <H/2 $, $\gamma <(3H/2 -\Re   M)/(\alpha+1)$,   we have
\begin{eqnarray*} 
\| \partial_tG[F\Psi(\psi)](t,x)\|_{H_{(s)} }  
  & \lesssim   &
e ^{- Ht}     \int_{ 0}^{t} e ^{ Hb} \|  \psi (b,x) \|^{1+\alpha} _{H_{(s)}}  db\\
  & \lesssim   &
e ^{- Ht}     \int_{ 0}^{t} e ^{ (H -\gamma (1+\alpha)) b} (e^{\gamma b }\|  \psi (b,x) \|_{H_{(s)}})^{1+\alpha}   db\\
  & \lesssim   &
  \varepsilon e ^{- Ht}     \int_{ 0}^{t} e ^{ (H -\gamma (1+\alpha)) b}   db\\
  & \lesssim   &
2 \varepsilon e ^{-\gamma (1+\alpha)  }     
  \quad \mbox{\rm if} \quad\Re M<\frac{H}{2}   \quad \mbox{\rm and} \quad \gamma <\frac{H}{1+\alpha}  \,.
\end{eqnarray*}
Similarly,   we have
\[ 
\| \partial_tG[V\psi](t,x)\|_{H_{(s)} }  
  \lesssim   
  \varepsilon e ^{-\gamma  t }      
  \quad \mbox{\rm if} \quad\Re M<\frac{H}{2} \quad \mbox{\rm and} \quad \gamma  <H \,.
\]
Thus, for $\Re M< {H}/{2} $,  the estimate  (\ref{der_t_i}) for the time derivative is proved. 
 
For    $ {H}/{2}<\Re M< {3H}/{2}  $, according to Theorem~\ref{T11.2}, we have
\[
    \|  \partial_t\psi_{id}(t) \|_{H_{(s)}}   
 \leq  
 C      e^{(\Re M-\frac{3}{2}H)t}      \left( \|\psi _0 \|_{H_{(s+1)}} +  \|\psi _1 \|_{H_{(s)}}\right)  \,.
\]
From (\ref{11.19})  with $\gamma  < \left(\frac{3}{2}H -\Re M\right)/(1+\alpha) <0$   we derive 
\begin{eqnarray*} 
\|  \partial_tG[F\Psi(\psi)](t,x)\|_{H_{(s)} } 
  & \lesssim   &
 e ^{(\Re M-\frac{1}{2}H)t}     \int_{ 0}^{t} e ^{-(\Re M-\frac{1}{2}H)b} \|  \psi (b,x) \|^{1+\alpha} _{H_{(s)}}  
   db     \\
  & \lesssim   &
 e ^{(\Re M-\frac{1}{2}H)t}     \int_{ 0}^{t} e ^{\delta b-(\Re M-\frac{1}{2}H+ \gamma  (1+\alpha) )b} \left(e^{\gamma  b}\|  \psi (b,x) \| _{H_{(s)}} \right)^{1+\alpha}  
   db     \,.
\end{eqnarray*}
Since
$
\delta-(\Re M-\frac{1}{2}H+ \gamma  (1+\alpha) )>0
$, 
we obtain
 \[ 
\|  \partial_tG[F\Psi(\psi)](t,x)\|_{H_{(s)} }
  \lesssim   
  \varepsilon e ^{( \delta-\gamma  (1+\alpha ) ) t}        ,
  \quad \mbox{\rm if} \quad\Re M>\frac{H}{2}  \quad \mbox{\rm and} \quad \gamma   (1+\alpha) <\left(\frac{3}{2}H -\Re M\right)\,.
\]
Similarly,  with $\gamma  <\left( { 3H}/{2} -\Re M\right)/(1+\alpha)  $ we obtain
\begin{eqnarray*} 
\|  \partial_tG[V\psi ](t,x)\|_{H_{(s)} } 
  & \lesssim   &
 e ^{(\Re M-\frac{1}{2}H)t}     \int_{ 0}^{t} e ^{-(\Re M-\frac{1}{2})Hb} \|  \psi (b,\cdot)  \| _{H_{(s)}} 
   db  \\
  & \lesssim   &
 e ^{(\Re M-\frac{1}{2}H)t}     \int_{ 0}^{t} e ^{\delta b-(\Re M-\frac{1}{2})Hb-\gamma  b} (e^{\gamma b}\|  \psi (b,\cdot)  \| _{H_{(s)}}) 
   db  \\
  & \lesssim   &
 \varepsilon e ^{(\delta-\gamma)    t}
  \quad \mbox{\rm if} \quad\Re M>\frac{H}{2}   \,.
\end{eqnarray*}
Thus, for    $ {H}/{2}<\Re M< {3H}/{2}  $, we have proved (\ref{der_t_ii}). This completes the proof of Theorem~\ref{TIE}. \qed


\begin{thebibliography}{99}

 
\bibitem{Andersson-Blue-Wang}
L.Andersson,  P.Blue, and J.Wang,   ``Morawetz estimate for linearized gravity in Schwarzschild,'' Ann. Henri Poincar\'e 21 , no. 3, 761--813 (2020).

\bibitem{Bachelot-Nicolas}
A.Bachelot,   J.-P.Nicolas,   
``\'Equation non lin\'eaire de Klein-Gordon dans des m\'etriques de type Schwarzschild,'' 
C. R. Acad. Sci. Paris S\'er. I Math. 316, no. 10, 1047--1050 (1993).
 

\bibitem{B-E}
H.~Bateman, A.~ Erdelyi, Higher  Transcendental  Functions, Vol. 1,2, McGraw-Hill, New York, 1954.

\bibitem{Birrell}
N.D.~Birrell,     P.C.W.~Davies,  
 { Quantum Fields in Curved Space},  Cambridge, New York, Cambridge University Press, 1984.



\bibitem{Catania-Georgiev} 
D.Catania, V.Georgiev,   ``Blow-up for the semilinear wave equation in the Schwarzschild metric,'' Differential Integral Equations 19, no. 7, 799--830 (2006). 




\bibitem{Chandrasekhar}
S.Chandrasekhar,   The mathematical theory of black holes. International Series of Monographs on Physics, 69. Oxford Science Publications. The Clarendon Press, Oxford University Press, New York, 1983.

\bibitem{Christ-Klainerman}
D.Christodoulou,   S.Klainerman,   
The global nonlinear stability of the Minkowski space.  
Princeton Mathematical Series, 41. Princeton University Press, Princeton, NJ, 1993.   

\bibitem{Choquet-Bruhat_book}
Y.Choquet-Bruhat, General relativity and the Einstein equations. Oxford Mathematical Monographs. Oxford University Press, Oxford~(2009)

\bibitem{Daf_Rod_2005b}
M.Dafermos,    I.Rodnianski,   ``A proof of Price's law for the collapse of a self--gravitating scalar field,'' 
Invent. Math. 162, no. 2, 381-457 (2005).

\bibitem{Daf_Rod_2005}
M.Dafermos,    I.Rodnianski,  ``Small-amplitude nonlinear waves on a black hole background,'' J. Math. Pures Appl. (9) 84 , no. 9, 1147--1172 (2005). 



\bibitem{Dafermos-Holzegel-Rodnianski}
 M.Dafermos,   G.Holzegel, and I.Rodnianski, ``The linear stability of the Schwarzschild solution to gravitational perturbations,'' Acta Math. 222, no. 1, 1--214 (2019).



\bibitem{Dimock}
J.Dimock, ``Scattering for the wave equation on the Schwarzschild metric,'' Gen.Relativ. Gravitation, 17,
numero 4, 1985, p.353-369.






\bibitem{Epstein-Moschella}
 H.Epstein,   U.Moschella,   ``de Sitter tachyons and related topics,'' Comm. Math. Phys. 336, no. 1, 381-430 (2015).

\bibitem{Faraoni}
V.Faraoni,   Cosmological and black hole apparent horizons. Lecture Notes in Physics, 907. Springer, Cham, 2015. 
 
\bibitem{Farrah}
D.Farrah  et all, ``Observational evidence for cosmological coupling of black holes and its implications for an astrophysical source of dark energy,'' The Astrophysical Journal Letters, 944:L31 (9pp), 2023 February 20. 
	arXiv:2302.07878  
 
\bibitem{Frankel}
T.Frankel,  
Gravitational curvature.
An introduction to Einstein's theory. W. H. Freeman and Co., San Francisco, Calif., 1979. 


\bibitem{NA2015}
  A.Galstian, K.Yagdjian,   ``Global solutions for semilinear Klein-Gordon equations in FLRW spacetimes,'' Nonlinear Anal. 113, 339-356 (2015).
  


\bibitem{Giorgi}
E.Giorgi,   ``The mathematics of stable black holes,'' Notices Amer. Math. Soc. 70, no. 4, 552--563 (2023).


\bibitem{Griffiths}
D.J.Griffiths,   (2017). Introduction to Quantum Mechanics. Cambridge, United Kingdom: Cambridge University Press. p. 415. 



\bibitem{Hawking}
  W. Hawking and G. F. R. Ellis, The Large Scale Structure of Space--Time, Cambridge Monographs on Mathematical Physics (Cambridge University Press, 1973).


\bibitem{Hintz_2021} 
P.Hintz,    
``Black hole gluing in de Sitter space,''  Comm. Partial Differential Equations 46, no. 7, 1280--1318 (2021).
 
 \bibitem{Hintz_2022_CMP} 
 P.Hintz,   ``A sharp version of Price's law for wave decay on asymptotically flat spacetimes,'' Comm. Math. Phys. 389, no. 1, 491--542 (2022).  
 
\bibitem{Hintz_2022}
  P.Hintz, Y.Xie,   ``Quasinormal modes of small Schwarzschild--de Sitter black holes,'' J. Math. Phys. 63, no. 1, Paper No. 011509, 26 pp (2022).



\bibitem{Kastor}
D.Kastor,   J.Traschen,   ``Cosmological multi-black-hole solutions,'' Phys. Rev. D (3) 47, no. 12, 5370--5375 (1993).

\bibitem{Lai-Zhou}
Ning-An Lai,  Y. Zhou,   ``Blow-up and lifespan estimate to a nonlinear wave equation in Schwarzschild spacetime,'' J. Math. Pures Appl. (9) 173, 172--194 (2023).
  
\bibitem{McVittie} 
G.C. McVittie, ``THE MASS-PARTICLE IN AN EXPANDING UNIVERSE,'' Mon. Not. R. Astron. Soc. 93, 325
(1933).  
  
 
\bibitem{Metcalfe} 
 J.Metcalfe,   G.Wang,   ``The Strauss conjecture on asymptotically flat space-times,'' SIAM J. Math. Anal. 49, no. 6, 4579--4594 (2017). 

\bibitem{Mizohata} 
S.Mizohata,   The theory of partial differential equations.  Cambridge University Press, New York, 1973. xii+490 pp.



\bibitem{Nakamura_JMAA2014}
 M.Nakamura,   ``The Cauchy problem for semi-linear Klein-Gordon equations in de Sitter spacetime,'' J. Math. Anal. Appl. 410 , no. 1, 445--454 (2014). 




\bibitem{Nicolas}
 J.-P.Nicolas,   ``Nonlinear Klein-Gordon equation on Schwarzschild--like metrics,'' J. Math. Pures Appl. (9) 74, no. 1, 35--58 (1995). 

\bibitem{Ohanian-Ruffini}    
H.~Ohanian, R.~Ruffini, Gravitation and Spacetime, Norton, New York, 1994.   






\bibitem{Parker}
L.~E.~Parker,   D.~J.~Toms, Quantum Field Theory in
Curved Spacetime, 
Quantized fields and gravity, Cambridge Monographs on Mathematical Physics, Cambridge University Press, Cambridge, 2009. 


\bibitem{Perlick}
  V.Perlick,   O.Y.Tsupko, and G.S.Bisnovatyi-Kogan,  ``Black hole shadow in an expanding universe with a cosmological constant,'' Phys. Rev. D 97, no. 10, 104062, 11 pp (2018).

 
\bibitem{Taylor}   
M.Taylor,   Partial differential equations III. Nonlinear equations. Second edition. Applied Mathematical Sciences, 117. Springer, New York, 2011.


\bibitem{YagBook}
K.Yagdjian,   
The Cauchy problem for hyperbolic operators.  
Multiple characteristics. Micro-local approach. Mathematical Topics, 12. Akademie Verlag, Berlin, 1997. 



\bibitem{Yag_Galst_CMP} 
K.~Yagdjian, A.~Galstian,  ``Fundamental solutions for the Klein-Gordon equation in de Sitter spacetime,'' 
Comm. Math. Phys. {285} 293--344 (2009). 

\bibitem{CPDE_2010}
K.Yagdjian, ``On the global solutions of the Higgs boson equation.'' Comm. Partial Differential Equations 37, no. 3, 447-478 (2012). 

\bibitem{Macao} 
K.Yagdjian, ``Integral transform approach to time--dependent partial differential equations,'' Mathematical analysis, probability and applications--plenary lectures.
Papers from the 10th International Congress (ISAAC 2015)   Macau, August 3--8, 2015. Edited by Tao Qian and Luigi G. Radino. 
 281--336, Springer Proc. Math. Stat., 177, Springer,  Cham, 2016.

\bibitem{JMP2019}
K.Yagdjian, ``Global existence of the self--interacting scalar field in the de Sitter universe,'' J. Math. Phys. 60, no. 5, 051503, 29 pp (2019).
 


\bibitem{MN2015}
K.~Yagdjian, ``Integral transform approach to solving Klein-Gordon equation with variable coefficients,'' Math. Nachr. 288, no. 17--18  2129--2152 (2015). 

 
\end{thebibliography}
\end{document}